\documentclass{article}
\usepackage[T1]{fontenc} 
\usepackage{newunicodechar}
\newunicodechar{‑}{\mbox{-}}
\newunicodechar{–}{--}

\usepackage{microtype}
\usepackage{graphicx}
\usepackage{subcaption}
\usepackage{booktabs} 
\usepackage[dvipsnames]{xcolor}
\usepackage{tikz} 

\usepackage{hyperref}

\newcommand{\theHalgorithm}{\arabic{algorithm}}

\usepackage[accepted]{icml2026}
\makeatletter
\renewcommand{\Notice@String}{}
\makeatother
\hypersetup{pdfsubject={Constrained combinatorial optimization with tensor trains}}
\icmltitlerunning{Exploration and Exploitation with U(1)-Symmetric Tensor Trains}

\usepackage{amsmath}
\usepackage{amssymb}
\usepackage{mathtools}
\usepackage{amsthm}

\usepackage[capitalize,noabbrev]{cleveref}

\theoremstyle{plain}

\theoremstyle{definition}

\theoremstyle{remark}

\usepackage[textsize=tiny]{todonotes}

\usepackage{placeins}
\begin{document}

\setcounter{topnumber}{4}
\setcounter{bottomnumber}{3}
\setcounter{totalnumber}{6}
\setcounter{dbltopnumber}{4}
\renewcommand{\topfraction}{0.92}
\renewcommand{\bottomfraction}{0.85}
\renewcommand{\textfraction}{0.08}
\renewcommand{\floatpagefraction}{0.75}
\renewcommand{\dbltopfraction}{0.92}
\renewcommand{\dblfloatpagefraction}{0.75}
\setlength{\floatsep}{8pt plus 2pt minus 2pt}
\setlength{\textfloatsep}{10pt plus 2pt minus 2pt}
\setlength{\dblfloatsep}{8pt plus 2pt minus 2pt}
\setlength{\dbltextfloatsep}{10pt plus 2pt minus 2pt}
\setlength{\intextsep}{8pt plus 2pt minus 2pt}
\raggedbottom

\twocolumn[
  \icmltitle{Balancing Exploration and Exploitation in Constrained Combinatorial Optimization with U(1)-Symmetric Tensor Trains}



  \icmlsetsymbol{equal}{*}

  \begin{icmlauthorlist}
    \icmlauthor{Sergei Iudin}{MIRIAI}
    \icmlauthor{Mikhail Podobrii}{MIRIAI}
    \icmlauthor{Dmitry Zheltkov}{INM}
    \icmlauthor{Stanislav Moiseev}{T}
  \end{icmlauthorlist}

  \icmlaffiliation{MIRIAI}{Moscow Independent Research Institute of Artificial Intelligence (MIRAI), Russian Federation}
  \icmlaffiliation{INM}{Institute of Numerical Mathematics RAS, Russian Federation}
  \icmlaffiliation{T}{T-Technologies}

  \icmlcorrespondingauthor{Sergei Iudin}{iudin.s@miriai.org}

  \icmlkeywords{Tensor Train, Constrained Optimization}

  \vskip 0.3in
]


\printAffiliationsAndNotice{}  

\begin{abstract}
  Tensor Train (TT) decomposition has emerged as a powerful tool for mitigating the curse of dimensionality in a variety of machine learning and combinatorial optimization tasks.  Recent advances show that TT‑based generative models can efficiently represent and explore high‑dimensional discrete spaces, even when the cost function structure is unknown.  However, incorporating hard linear constraints into TT‑based optimization remains challenging, as existing interpolation and sampling approaches often fail when feasible points are rare.  In this work, we develop a general TT‑based framework for constrained combinatorial optimization built upon U(1)‑symmetric tensor networks.  Our method introduces an expressibility‑enhancement mechanism for constructing sparse symmetric tensors and a graph‑based diversity measure that guides sample selection during training.  On assignment benchmarks, expressibility enhancement and diversity-aware selection improve mean solution quality over their respective baselines. The assignment-selected diversifying strategy also extends to nonpolynomial Set Partitioning, with larger relative gains observed in the higher-tensor-rank data families. Separate runtime benchmarks show substantial speedups over ConstrainTNet and PROTES on non-degenerate assignment workloads; comparisons with simulated annealing depend on constraint information and computational budget, but, overall, support the practicality of the current implementation.
\end{abstract}

\section{Introduction}


Tensor Train (TT) decomposition has become a powerful tool in the domains of Artificial Intelligence and Combinatorial Optimization in recent years. 
The idea of squeezing high-dimensional data via tensor train decomposition originated in the field of quantum mechanics \cite{affleck1987rigorous, white1992density}, where many-dimensional systems naturally arise from particle interactions.  In 2011, Oseledets formalized the approach and brought these ideas into computational mathematics \cite{oseledets2011tensor}. 
To date, TT decomposition was applied to low‑rank network layer approximation \cite{novikov2015tensorizing}, LLM fine‑tuning \cite{yang2024loretta}, optimization \cite{antil2023ttrisk, mugel2022dynamic}, options pricing \cite{arenstein2025full}, and some others \cite{dolgov2020parallel, iudin2025multistart,  sulimov2017evaluation, sozykin2026global}

This paper focuses on applications of TT in combinatorial optimization. 
In comparison to classical mixed-integer linear programming (MILP)
TT is not limited to linear cost functions and does not use any integrality relaxation. In comparison to Constraint Programming (CP) solvers 
\cite{biere2009handbook} 
and QAOA \cite{farhi2014quantum}, the cost function structure may be unknown. Compared with Bayesian optimization \cite{jones1998efficient}, TT can offer favorable scaling when the required ranks remain moderate 
. 

We know two paradigms for TT-based optimization. First, we may request cost function values at specific grid points, then interpolate the entire landscape using a low TT-rank assumption \cite{zheltkov2019global, sozykin2022ttopt}. 
Second, we may apply a generative approach \cite{batsheva2023protes}: based on some data, learn a probability distribution that encourages low-cost values, and store it as TT. Then we may sample new data from it and use the data for further training until a stopping criterion is reached. TT-format
speeds up both the training and sampling phases, making high-dimensional distribution tractable.

The challenge of applying TT to combinatorial optimization appears when variables are constrained to a feasible region (e.g., by linear equalities). In such a case, building a TT surrogate for the cost function is difficult due to a large fraction of infeasible grid points. 
The probability of finding any feasible point is small, and interpolation based on infeasible points is not informative.
On the contrary, even a straightforward generative approach may be pre-trained on feasible data, and in a sampling phase, all infeasible samples may be eliminated. However, if the constraints are hard to satisfy, the algorithm's efficiency may decrease significantly. To address the challenge of constraint embedding, the authors in \cite{nakada2025quick} used the method of nilpotent matrices. While effective for some practical cases, it is not generalizable to arbitrary linear equations because of exponential scaling with the number of constraints. In \cite{Ryzhakov2023Constructive}, the authors represent the method of exact constraints encoding, but only for constraints that have a particular structure, which is consistent with the tensor train linear topology and locality of interactions. The paper \cite{LopezPiqueres2023} suggests using $U(1)$- symmetric tensor networks for problems with linear constraints. The available feasible dataset was used to build a sparse probability distribution tensor, assigning zero probability to all infeasible vectors. Later, an exact framework for incorporating any linear inequalities into a symmetric TT structure was developed \cite{10.21468/SciPostPhys.18.6.192}, with a code available on git \cite{ConstrainTNet2024}. 

\textbf{Main contributions.} We developed a library for constrained optimization, using tensor trains. Our approach is based on a generative $U(1)$-symmetric paradigm, which accommodates arbitrary integer linear equalities, although computational cost depends on the resulting charge structure. Two major contributions expand the paradigm: i) We suggested the \textit{expressibility enhancement} step, increasing feasible dataset exploitation. Thus, instead of only two strategies of building a sparse $U(1)$-symmetric tensor from the feasible data, as in the pioneer paper, we suggest a family of such strategies, controlled by a parameter (exploitation power). The tested enhancement improves mean assignment solution quality over the original construction. ii) We represented a sparse symmetric tensor train as a graph of $U(1)$ charges, and suggested a graph diversity measure, which helps to choose which samples to use for training on the next algorithm iteration (increasing exploration).
The assignment-selected diversifying strategy also improves mean quality on entropic Set Partitioning, with benefits concentrated in higher-tensor-rank data families. Our Julia implementation achieves $77$--$448\times$ speedups over ConstrainTNet \cite{ConstrainTNet2024} on the specified non-degenerate assignment kernel benchmark, depending on size and sampling backend.

\section{Problem Statement}
We consider a similar class of optimization problems as in \cite{LopezPiqueres2023}: 

\begin{equation}
\left\{
\begin{aligned}
    & C(x) \;\to\; \min, \\[4pt]
    & A x = b, \\[4pt]
    & x \in \{0,1\}^N ,
\end{aligned}
\right.
\label{eq:problem_statement}
\end{equation}

where C : $\{0,1\}^N \to \mathbb{R} $ is a black-box objective function, $A$ and $b$ are integer-valued. We emphasize that the domain of the $x$ variable may be, in principle, any finite integer domain; also, $A$ and $b$ may be rational, if they are obtained from an integer by dividing by some scaling factor.

Generic black-box solvers can struggle with hard constraints. E.g., Simulated Annealing (SA) \cite{kirkpatrick1983optimization} requires a random step implementation that maintains feasibility. Or constraints may be implemented as penalties of the form $\lambda \|Ax-b\|_2^2$, but this may be even more problematic, because of the necessity to choose a value for $\lambda$ and probably a low fraction of feasible solutions evaluated. In this work, we assume that some feasible data for hard constraints is available. This may be due to one of the following reasons: i) There is a feasible dataset sampler available, which may produce diverse feasible samples using some knowledge about the matrix $A$ and vector $b$. The optimization itself may still be very challenging due to the black-box function $C(x)$, which, in the worst case, takes a long time to evaluate. ii) Feasible space may be explored in some other way, e.g., with a random search / MILP / CP solver, etc.

\section{Tensor Train Generative Modeling }
\label{sec:generative}

Algorithm~\ref{alg:tt_generative} describes a family of generative optimization
procedures based on tensor-train (TT) probability models. The method maintains
a distribution $P_\theta(x)$ for an integer vector $x=(x_1,\ldots,x_N)$:
\begin{equation}
\begin{gathered}
\pi_\theta(x)=\sum_{\alpha_1,\ldots,\alpha_{N-1}}
G^{(1)}_{1,x_1,\alpha_1}G^{(2)}_{\alpha_1,x_2,\alpha_2}\cdots
G^{(N)}_{\alpha_{N-1},x_N,1},\\
P_\theta(x)=|\pi_\theta(x)|^2/Z_\theta
\quad\text{or}\quad P_\theta(x)=\pi_\theta(x).
\label{eq:tt-def}
\end{gathered}
\end{equation}
The first choice is a Born machine; the second requires a nonnegative,
normalized TT. The implementations considered below use the Born rule.
An initial feasible dataset supplies a starting model and a best-so-far
solution. At each iteration, the model generates candidates, a selection
policy constructs the next training set, and a weighted likelihood update
biases the model toward promising solutions. The incumbent is maintained
separately from training data and is never discarded by diversification.

\textbf{Selection and memory.} The candidate pool may contain only current
TT samples, or also externally generated feasible points and retained data.
The selection policy may favor low cost, charge diversity, or their combination;
it also specifies whether initial points or a history participate. Neither
elite-only selection nor cumulative retention is intrinsic to the framework.
For finite computational cost, experiments specify a selection quota and a
memory policy explicitly.

\textbf{Optional diversification.} Before fitting the selected data, an
operator $\mathcal D_t$ may leave the TT unchanged, reset its amplitudes, or apply a preliminary
sweep with equal sample weights. The last choice flattens the empirical
training target; a single sweep does not guarantee a uniform model over its
entire feasible support. Diversification can be periodic or triggered by
stagnation. It must preserve feasibility and provide nonzero support for
positive-weight training points. These options define different algorithms,
not interchangeable implementations of the same experiment.
With $\mathcal D_t$ equal to the identity, learned amplitudes persist across
iterations. The comparison in \cref{subsec:eels_optimization} uses this
choice; In \cref{subsec:add}, algorithm rebuilds the model after selection, while the SA comparison in
\cref{subsec:sa_comparison} uses the explicitly stated reset/rebuild variant.
The schematic loop can be cyclically organized as train--sample--select
when the initial dataset is trained before the first sampling phase.

For the constrained problem \eqref{eq:problem_statement}, it is much more effective to initialize the tensor with some sparse structure, with zero probabilities for vectors that violate constraints. It helps avoid generating such samples that violate the problem's hard constraints. For this purpose, we exploit properties of $U(1)$-symmetry. Let us define it carefully, but first, give some preliminary definitions. We consider a tensor $\pi_{x_1\dots x_N}$ for which the following holds: 
\begin{itemize}
    \item Each index $x_k$ has one of two possible orientations (directions): $x_k \in I_{\text{in}}$ or $x_k \in I_{\text{out}}$. We call $I_{\text{in}}$ a set of incoming indices and $I_{\text{out}}$ a set of outgoing indices.
    
    \item Each value of index $x_k$ may be decomposed into the pair of values $(c_{x_k}, t_{c_{x_{k}}})$, where we call $c_{x_k}$ a \textit{charge} index, and  $t_{c_{x_{k}}}$ a \textit{degeneracy} index. In this paper, the charges may be integers or integer vectors. More formally, for each $k$ linear space $V^{(k)}$ (associated with index $x_k$) is a direct sum of charge subspaces:
    \begin{equation}
    \label{eq:direct_sum}
        V^{(k)} = \bigoplus_{c_k}V^{(c_k)}.
    \end{equation}
\end{itemize}

Let $C_{\text{in}} = \sum_{x_i \in I_{\text{in}}}c_{x_i}$, and $C_{\text{out}}=\sum_{x_i\in I_{\text{out}}}c_{x_i}$.
We say that the set of index values $(x_1,\dots x_N)$ \textit{conserves the charge} $c$ if $C_{\text{in}} - C_{\text{out}}=c$.
Finally, we call the tensor $\pi_{x_1\dots x_d}$ $U(1)$-symmetric if the following holds for some charge $c$ and all possible values of indices $(x_1,\dots x_N)$:

\begin{equation}
    \pi_{x_1\dots x_N} = \left(\mathcal{P}_{c_{1}\dots c_{N}}\right)_{t_{c_1}\dots t_{c_N}} \delta_{C_{\text{in}},  C_{\text{out}} + c},
\label{eq:charge_conservation}
\end{equation}

where $\delta_{i,j}$ is the Kronecker delta symbol, $c_k=c_{x_k} \quad \forall k$ and $\mathcal{P}_{c_{1}\dots c_{N}}$ is a tensor indexed only by charges. Informally, that means that $\pi$ becomes block-sparse. If $\pi_{x_1,\dots, x_N}$ is $U(1)$-symmetric, then $\pi_{x_1,\dots, x_N}=0$ for all sets of values $(x_1,\dots, x_N)$ that do not conserve the charge $c$.  In such cases, $c$ is called a $flux$ of the tensor. We will also call \eqref{eq:charge_conservation} a \textit{charge conservation law}.

\begin{algorithm}
\caption{Generative optimization with a tensor train}
\label{alg:tt_generative}
\begin{algorithmic}[1]
\REQUIRE Initial feasible data $\mathcal T_0$, sampling budget $K$, iteration budget $T$, selection/memory policy $\mathcal S_t$, optional diversification $\mathcal D_t$
\STATE Initialize $P_\theta$ and incumbent $x^\star$ from $\mathcal T_0$; initialize hystory state $H$
\FOR{$t=1,\ldots,T$}
    \STATE Sample $K$ candidates from $P_\theta$; optionally augment the pool from available feasible-data sources
    \STATE Evaluate candidates as prescribed; update $x^\star$ using evaluated feasible points
    \STATE $(\mathcal T_{\mathrm{new}},H)\gets\mathcal S_t(\text{candidates},\mathcal T_0,H)$
    \STATE $\theta\gets\mathcal D_t(\theta,\mathcal T_{\mathrm{new}},H)$ \COMMENT{Identity if omitted}
    \STATE Assign normalized weights $w_i$ to $x^{(i)}\in\mathcal T_{\mathrm{new}}$
    \STATE Update $\theta$ by minimizing $\displaystyle -\sum_{x^{(i)}\in\mathcal T_{\mathrm{new}}}w_i\log P_\theta(x^{(i)})$
\ENDFOR
\STATE \textbf{return} $x^\star$
\end{algorithmic}
\end{algorithm}

\subsection{Generative modeling with $U(1)$-symmetry}
\label{subsec:U1_generative}

The rigorous description of a generative algorithm based on $U(1)$-symmetric tensor trains may be found in \cite{LopezPiqueres2023}. All crucial steps from \cref{alg:tt_generative} remain the same, but require less memory for storing $P_{\theta}$, as only nonzero blocks are stored. Also, the big saving of computational resources takes place while calculating the gradients of $\mathcal{L}(\theta)$ and applying gradient descent steps, because all these operations are done on sparse tensors.

In detail, we describe only initialization from the training data. The idea is to associate each variable $x_i$ in problem \eqref{eq:problem_statement} with a charge $x_i \vec{A_i}$, where $\vec{A_{i}}$ is the $i$-th column of the constraints matrix $A$. Then, the probability amplitude $\pi(x)$ is given as a symmetric TT with the flux $b$, and each TT-core $G^{(k)}$ is a $U(1)$-symmetric tensor of rank 3. Then, for each TT-core, we choose the direction of the indices. Left (virtual) link index $\alpha_{k-1} \in I_{\text{in}}$, (site) index $x_k \in I_{\text{out}}$, right (virtual) link index $\alpha_{k}\in I_{\text{out}}$:

\[
G^{(k)}(\alpha_{k-1}, x_k, \alpha_{k}) \;=\;
\vcenter{\hbox{\begin{tikzpicture}[baseline=(G.base), every node/.style={circle, draw, minimum size=5mm}]
  \node[fill=cyan] (G) {};
  \draw[<-] (G.west) -- ++(-0.5,0) node[left, draw=none] {$\alpha_{k-1}$};
  \draw[->] (G.east) -- ++(0.5,0) node[right, draw=none] {$\alpha_{k}$};
  \draw[->] (G.south) -- ++(0,-0.5) node[below, draw=none] {$x_k$};
  \fill (G.west) ++(-0.2,0) circle (0.6pt);
  \fill (G.east) ++(0.2,0) circle (0.6pt);
  \fill (G.south) ++(0,-0.2) circle (0.6pt);
\end{tikzpicture}.}}
\]

By convention, we put the flux on the $1$-st TT-core, on the first virtual index of $G^{(1)}$. Thus, the example TT, which encodes the probability distribution for the problem with 4 binary variables, is:

\begin{center}
\begin{tikzpicture}[every node/.style={circle, draw, minimum size=5mm}]
    \node[fill=cyan] (x1) at (-2.25, -5) {};
    \node[fill=cyan] (x2) at (-0.75, -5){};
    \node[fill=cyan] (x3) at (0.75, -5) {};
    \node[fill=cyan] (x4) at (2.25, -5) {};

    \node[fill=none, draw=none] (x1_) at (-2.25, -4.5) {$G^{(1)}$};
    \node[fill=none, draw=none] (x2_) at (-0.75, -4.5) {$G^{(2)}$};
    \node[fill=none, draw=none] (x3_) at (0.75, -4.5) {$G^{(3)}$};
    \node[fill=none, draw=none] (x4_) at (2.25, -4.5) {$G^{(4)}$};

    \node[fill=none, draw=none] (c1_) at (-2.25, -6) {$\{\vec{A_1}, \vec{0}\}$};
    \node[fill=none, draw=none] (c2_) at (-0.75, -6) {$\{\vec{A_2}, \vec{0}\}$};
    \node[fill=none, draw=none] (c3_) at (0.75, -6) {$\{\vec{A_3}, \vec{0}\}$};
    \node[fill=none, draw=none] (c4_) at (2.25, -6) {$\{\vec{A_4}, \vec{0}\}$};

    \draw[draw=black, line width=1pt, ->] (x1) -- (x2);
    \draw[draw=black, line width=1pt, ->] (x2) -- (x3);
    \draw[draw=black, line width=1pt, ->] (x3) -- (x4);

    \draw[->] (x1) -- ++(-90:0.8cm);
    \draw[->] (x2) -- ++(-90:0.8cm);
    \draw[->] (x3) -- ++(-90:0.8cm);
    \draw[->] (x4) -- ++(-90:0.8cm);

    \draw[<-, dashed] (x1) -- node[midway, above, pos=0.5, draw=none, fill=none] {$\vec{b}$} ++(-150:1.2cm);
    \draw[->, dashed] (x4) -- node[midway, above, pos=0.5, draw=none, fill=none] {$\vec{0}$} ++(-30:1.2cm);
\end{tikzpicture}.
\end{center}
The charge conservation law \eqref{eq:charge_conservation} is then equivalent to constraints in \eqref{eq:problem_statement}, because it says, that $P(x) \neq 0$ only if $\vec{b}-x_1\vec{A_1}-x_2\vec{A_2}-x_3\vec{A_3}-x_4\vec{A_4}=\vec{0}$.

Building on an exact symmetric TT for the problem \eqref{eq:problem_statement} may be exponentially hard, because it requires an exponentially large number of link charges (thus, exponentially large rank of tensor train). In \cite{LopezPiqueres2023}, heuristic strategies are suggested, which limit the minimum required rank to construct a symmetric tensor train, but also build only the approximation of the required tensor. The two constructions are described below. We chose \textit{Method (ii)} (the most expressive) for this paper as the baseline:

\textit{Method (i)}: Given the feasible dataset $\mathcal{T}$, for each $x\in \mathcal{T}$, we calculate all the link charges $c_{\alpha_k}$ by propagating: $c_{1} = \vec{b}$, $c_{\alpha_{i}} = c_{\alpha_{i-1}}-x_i\vec{A_i}$. It requires $|b||\mathcal{T}| (N-1)$ operations, where $|\mathcal{T}|$ is the number of (unique) feasible samples, and $|b|$ is the number of equalities in the problem. Then, we initialize each TT-core $G^{(k)}$ with non-zero blocks for subspace with charges $(c_{\alpha_{k-1}}, x_k \vec{A_k}, c_{\alpha_k})$.

    \textit{Method (ii)}: Given the feasible dataset $\mathcal{T}$, for each $x\in \mathcal{T}$, we calculate all the link charges $c_{\alpha_k}$ by propagating: $c_{1} = \vec{b}$, $c_{\alpha_{i}} = c_{\alpha_{i-1}}-x_i\vec{A_i}$. It requires $|b||\mathcal{T}| (N-1)$ operations, where $|\mathcal{T}|$ is the number of (unique) feasible samples, and $|b|$ is the number of equalities in the problem. Then, for each $G^{(k)}$, we search for consistent (with the left link charges $c_{\alpha_{k-1}}$ and right link charges $c_{\alpha_k}$) site charges $c_{x_k} \in \{\vec{0}, \vec{A_k}\}$. This step requires checking for all $c_{\alpha_{k-1}}$ and $c_{x_k}$, if $(c_{\alpha_{k-1}}-c_{x_k}) \in \bigcup_{\alpha_k}\{c_{\alpha_k}\}$. We store $\bigcup_{\alpha_k}\{c_{\alpha_k}\}$ as a hash-table based type, thus all the checks take $O(N|b||\mathcal{T}|)$ operations on average, and $O(N|b||\mathcal{T}|^2)$ in the worst case (many hash collisions). TT-rank grows with the dataset size as $O(|\mathcal{T}|)$.
    
The process of building a nonzero subspace of the tensor train may be interpreted as building a directed acyclic \textit{charge graph}, with nodes corresponding to link charges, edges corresponding to values of $x_i$, and paths in this graph from the initial node $\vec{b}$ to the final node $\vec{0}$ corresponding to feasible solutions to $Ax=b$. For instance, look at the \cref{fig:example_graph} (a). In these terms, the algorithm of encoding constraints in $U(1)$-symmetric TT is expressive if it generates many feasible paths in the charge graph.

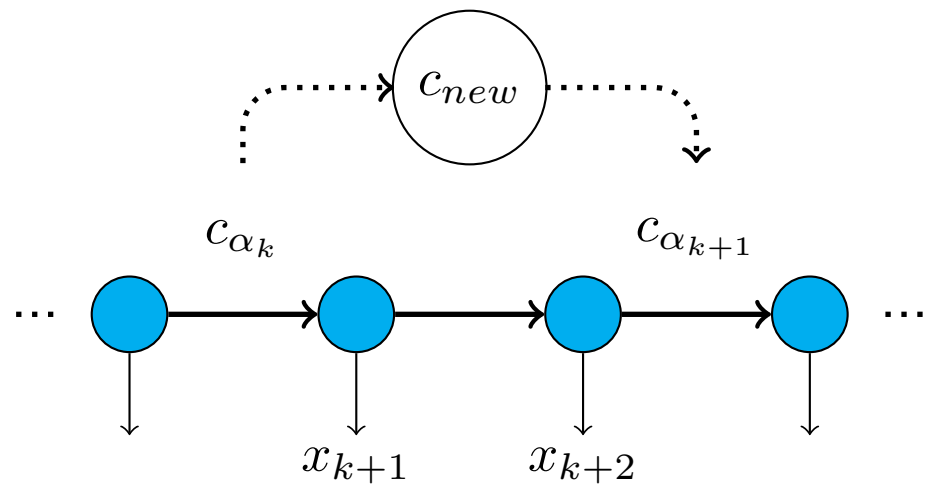
\begin{figure}[!tbp]
\begin{center}
\begin{tikzpicture}[scale=1.0, every node/.style={circle, draw, minimum size=5mm}]
    \node[fill=cyan] (x1) at (-2.25, -5) {};
    \node[fill=cyan] (x2) at (-0.75, -5){};
    \node[fill=cyan] (x3) at (0.75, -5) {};
    \node[fill=cyan] (x4) at (2.25, -5) {};

    \node[fill=none, draw=none] (x1_) at (-2.25, -4.5) {};
    \node[fill=none, draw=none] (x2_) at (-0.75, -4.5) {};
    \node[fill=none, draw=none] (x3_) at (0.75, -4.5) {};
    \node[fill=none, draw=none] (x4_) at (2.25, -4.5) {};

    \node[fill=none, draw=none] (c2_) at (-0.75, -6) {$x_{k+1}$};
    \node[fill=none, draw=none] (c3_) at (0.75, -6) {$x_{k+2}$};
    
    \draw[draw=black, line width=1pt, ->] (x1) -- (x2);
    \draw[draw=black, line width=1pt, ->] (x2) -- (x3);
    \draw[draw=black, line width=1pt, ->] (x3) -- (x4);

    \draw[->] (x1) -- ++(-90:0.8cm);
    \draw[->] (x2) -- ++(-90:0.8cm);
    \draw[->] (x3) -- ++(-90:0.8cm);
    \draw[->] (x4) -- ++(-90:0.8cm);

     \node[fill=none] (cnew) at (0, -3.5) {$c_{new}$};
     \draw[->, line width=1pt, rounded corners=8pt, dotted]
    (-1.5,-4) -- (-1.5,-3.5) -- (-0.5, -3.5);
    \draw[->, line width=1pt, rounded corners=8pt, dotted]
    (0.5, -3.5) -- (1.5, -3.5) -- (1.5, -4);
    \node[draw=none] (ck) at (-1.5, -4.5) {$c_{\alpha_k}$};
    \node[draw=none] (cknext) at (1.5, -4.5) {$c_{\alpha_{k+1}}$};
    \draw[-, dotted, line width=1pt] (-3, -5)--(-2.75, -5);
    \draw[-, dotted, line width=1pt] (2.75, -5)--(3, -5);
    
\end{tikzpicture}
\caption{One step of the EELS algorithm for $\gamma=1$.}
\label{fig:enhancement_scheme}
\end{center}
\end{figure}

\begin{figure*}[!tbp]
\centering
\begin{tikzpicture}[scale=1.0, every node/.style={transform shape, circle, draw, minimum size=8mm}]

    \node[fill=none, draw=none] (a) at (-4, 3) {$a)$};
    \node[fill=none, draw=none] (b) at (-4, -4) {$b)$};
    \node[fill=none, draw=none] (c) at (5.5, 3) {$c)$};
    \node[fill=none, draw=none] (d) at (5.5, -4) {$d)$};
    
    \node[line width=1pt] (t1) at (-4, 0) {$2$};
    \node[line width=1pt] (t2) at (-2, 2) {$1$};
    \node[line width=1pt] (t3) at (0, 2) {$0$};
    \node[line width=1pt] (t4) at (2, 2) {$0$};
    \node[line width=1pt] (t5) at (4, 0) {$0$};

    \node[draw=red, line width=2pt] (m) at (0, 0) {$1$};

    \node[line width=1pt] (b1) at (-2, -2) {$2$};
    \node[line width=1pt] (b2) at (0, -2) {$2$};
    \node[line width=1pt] (b3) at (2, -2) {$1$};

    \draw[draw=cyan, line width=1pt, ->] (t1) -- node[midway, above, pos=0.2, draw=none, fill=none] {$x_1=1$}  (t2);
    \draw[draw=cyan, line width=1pt, ->] (t2) -- node[midway, above, draw=none, fill=none] {$x_2=1$} (t3);
    \draw[draw=cyan, line width=1pt, ->] (t3) -- node[midway, above, draw=none, fill=none] {$x_3=0$}(t4);
    \draw[draw=cyan, line width=1pt, ->] (t4) -- node[midway, above, draw=none, fill=none, pos=0.8] {$x_4=0$}(t5);

    \draw[draw=cyan, line width=1pt, ->] (t1) -- node[midway, below, pos=0.2, draw=none, fill=none] {$x_1=0$}(b1);
    \draw[draw=cyan, line width=1pt, ->] (b1) -- node[midway, below, draw=none, fill=none] {$x_2=0$}(b2);
    \draw[draw=cyan, line width=1pt, ->] (b2) -- node[midway, below, draw=none, fill=none] {$x_3=1$}(b3);
    \draw[draw=cyan, line width=1pt, ->] (b3) -- node[midway, below, draw=none, fill=none, pos=0.8] {$x_4=1$}(t5);

    \draw[draw=ForestGreen!70!white, line width=1pt, ->] (t2) -- (m);
    \draw[draw=ForestGreen!70!white, line width=1pt, ->] (m) -- (b3);
    \draw[draw=Yellow!70!white, line width=1pt, ->] (b1) -- (m);
    \draw[draw=Yellow!70!white, line width=1pt, ->](m) -- (t4);

    \node[fill=cyan] (x1) at (-3, -5) {};
    \node[fill=cyan] (x2) at (-1, -5) {};
    \node[fill=cyan] (x3) at (1, -5) {};
    \node[fill=cyan] (x4) at (3, -5) {};

    \node[fill=none, draw=none] (x1_) at (-3, -6.5) {$x_1$};
    \node[fill=none, draw=none] (x2_) at (-1, -6.5) {$x_2$};
    \node[fill=none, draw=none] (x3_) at (1, -6.5) {$x_3$};
    \node[fill=none, draw=none] (x4_) at (3, -6.5) {$x_4$};

    \draw[draw=black, line width=1pt, ->] (x1) -- 
    node[midway, above, pos=0.5, draw=none, fill=none] {$\{1, 2\}$}  
    (x2);
    \draw[draw=black, line width=1pt, ->] (x2) -- 
    node[midway, above, pos=0.5, draw=none, fill=none] {$\{0, \textcolor{red}{1}, 2\}$}  
    (x3);
    \draw[draw=black, line width=1pt, ->] (x3) -- 
    node[midway, above, pos=0.5, draw=none, fill=none] {$\{0, 1\}$}  
    (x4);

    \draw[->] (x1) -- ++(-90:1.2cm);
    \draw[->] (x2) -- ++(-90:1.2cm);
    \draw[->] (x3) -- ++(-90:1.2cm);
    \draw[->] (x4) -- ++(-90:1.2cm);

    \draw[<-, dashed] (x1) -- node[midway, above, pos=0.5, draw=none, fill=none] {$\{2\}$} ++(-150:1.5cm);
    \draw[->, dashed] (x4) -- node[midway, above, pos=0.5, draw=none, fill=none] {$\{0\}$} ++(-30:1.5cm);

    \node[anchor=west, fill=none, draw=none] at (6,2) {$x_1 + x_2 + x_3 + x_4 = 2$};
    \node[anchor=west, fill=none, draw=none] at (6, 1.2) {$x_i \in \{0, 1\},\quad i = 1\dots4$};
    \node[anchor=west, fill=none, draw=none] at (5.5,-1) {\(
\mathcal{T} = \left[
\begin{array}{c|c}
\textcolor{cyan}{1} & \textcolor{cyan}{0} \\
\textcolor{cyan}{1} & \textcolor{cyan}{0} \\
\textcolor{cyan}{0} & \textcolor{cyan}{1} \\
\textcolor{cyan}{0} & \textcolor{cyan}{1}
\end{array}
\right]
\)};

    \node[anchor=west, fill=none, draw=none] at (5.5,-5) {\(
\mathcal{X} = \left[
\begin{array}{c|c|c|c}
\textcolor{cyan}{1} & \textcolor{cyan}{0} & \textcolor{ForestGreen!70!white}{1} & \textcolor{Yellow!70!white}{0} \\
\textcolor{cyan}{1} & \textcolor{cyan}{0} & \textcolor{ForestGreen!70!white}{0} & \textcolor{Yellow!70!white}{1} \\
\textcolor{cyan}{0} & \textcolor{cyan}{1} & \textcolor{ForestGreen!70!white}{0} & \textcolor{Yellow!70!white}{1}\\
\textcolor{cyan}{0} & \textcolor{cyan}{1} & \textcolor{ForestGreen!70!white}{1} & \textcolor{Yellow!70!white}{0}
\end{array}
\right]
\)};
\end{tikzpicture}
\caption{Example problem, where EELS reveals new samples, while standard approaches do not find any new feasible solutions. a) Charge graph of the problem. Cyan paths correspond to the training data, green and yellow lines show new feasible samples after local search. Red charge is a "newborn" charge. b) TT and the link charges nodes, obtained from the training data (black) and found by local search (red). c) Problem constraints and initial feasible dataset. d) Four illustrated feasible samples; the constraint has six feasible binary solutions in total. The sample color in $\mathcal{X}$ matches the color of the corresponding path in the graph.}
\label{fig:example_graph}
\end{figure*}
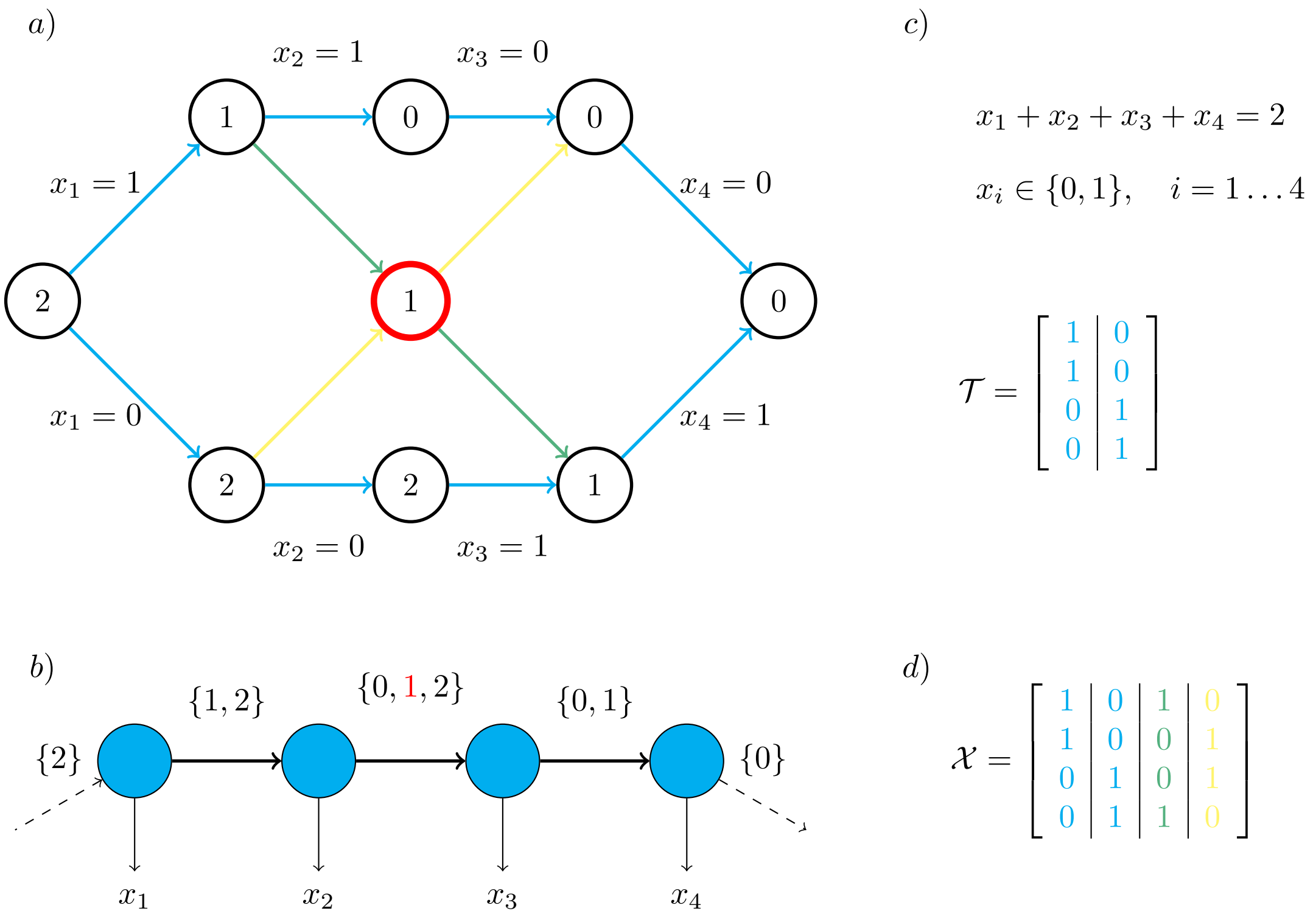

\subsection{Training the symmetric Born machine}
\label{subsec:training}

We use the Born parametrization $P_\theta(x)=\pi_\theta(x)^2/Z_\theta$
with real TT cores and $Z_\theta=\sum_x\pi_\theta(x)^2$.
Training follows the two-site, DMRG-like likelihood optimization of
\citet{LopezPiqueres2023}, Sec.~IV~B. For a fixed selected dataset
$\mathcal{T}=\{x^{(s)}\}_{s=1}^{M}$ (the selected $\mathcal{T}_{\mathrm{new}}$ in \cref{alg:tt_generative}, with retention determined by the experiment), define
\begin{equation}
 w_s=\frac{\exp[-(C(x^{(s)})-C_{\min})/\tau]}
 {\sum_{r=1}^{M}\exp[-(C(x^{(r)})-C_{\min})/\tau]},
 \qquad \tau>0,
 \label{eq:training_weights}
\end{equation}
where $C_{\min}=\min_{x\in\mathcal{T}}C(x)$ is a numerical shift
that cancels in the ratio. The training objective is
\begin{equation}
 \mathcal{L}=\log Z_\theta
       -2\sum_{s=1}^{M}w_s\log|\pi_\theta(x^{(s)})|.
 \label{eq:training_nll}
\end{equation}
The data and weights are held fixed throughout the inner training sweeps.
This objective fits the reweighted empirical distribution; it does not
require evaluating $C$ outside the selected data during a sweep.

\paragraph{Merge and differentiate.}
Put the TT in mixed canonical form with the orthogonality center at site
$k$ (or $k+1$ for a leftward update). Contract the shared bond of the
two adjacent cores:
\begin{equation}
 B^{(k)}_{aijb}=\sum_c G^{(k)}_{aic}G^{(k+1)}_{cjb}.
 \label{eq:training_merge}
\end{equation}
Here $i,j$ are physical values and $a,b,c$ are virtual indices, including
their charge and degeneracy labels. All cores outside this pair are
isometries. Writing $\|\cdot\|_F$ for the Frobenius norm, we have $Z_\theta=\|B^{(k)}\|_F^2$. Let
$L_a^{(s)}$ and $R_b^{(s)}$ be the contractions of the fixed sample with
the cores strictly to the left and right of the pair, respectively.
Removing the pair leaves the sample environment
\begin{equation}
 E^{(s)}_{aijb}=L_a^{(s)}R_b^{(s)}
   \delta_{i,x_k^{(s)}}\delta_{j,x_{k+1}^{(s)}}
   =\frac{\partial\pi_\theta(x^{(s)})}{\partial B^{(k)}_{aijb}}.
 \label{eq:training_environment}
\end{equation}
For real parameters the gradient and local update are
\begin{equation}
 \begin{aligned}
 \nabla_B\mathcal{L}
    &=\frac{2B}{Z_\theta}
       -2\sum_{s=1}^{M}w_s\frac{E^{(s)}}{\pi_\theta(x^{(s)})},\\
 \widetilde B&=B-\eta\nabla_B\mathcal{L},
 \end{aligned}
 \label{eq:training_gradient}
\end{equation}
where $\eta$ is the learning rate and the amplitudes and environments are
evaluated before the update. In particular, the norm term must be retained
even if the state was normalized before training.

\begin{figure}[!tbp]
\centering
\begin{tikzpicture}[x=1cm,y=1cm,font=\scriptsize,
 tn/.style={draw,rounded corners=1pt,minimum width=7mm,minimum height=6mm,fill=cyan!15},
 center/.style={tn,fill=orange!35},
 env/.style={tn,fill=gray!12}]
\node[anchor=west] at (0,1.05) {(a) Merge at the canonical center};
\node[center] (g1) at (1,0) {$G^{(k)}$};
\node[tn] (g2) at (2.4,0) {$G^{(k+1)}$};
\draw (0.35,0)--(g1); \draw (g1)--node[above] {$c$}(g2);
\draw (g2)--(3.05,0);
\draw (g1)--(1,0.55) node[above] {$i$};
\draw (g2)--(2.4,0.55) node[above] {$j$};
\draw[->] (3.35,0)--(4.1,0);
\node[center,minimum width=13mm] (b) at (5.25,0) {$B$};
\draw (4.25,0)--(b); \draw (b)--(6.25,0);
\draw (4.95,0.3)--(4.95,0.55) node[above] {$i$};
\draw (5.55,0.3)--(5.55,0.55) node[above] {$j$};
\node[anchor=west] at (0,-0.65) {(b) Remove the pair: sample environment $E^{(s)}$};
\node[env] (l) at (1,-1.85) {$L^{(s)}$};
\node[env] (r) at (5.5,-1.85) {$R^{(s)}$};
\draw (l)--(2.45,-1.85) node[above] {$a$};
\draw (r)--(4.05,-1.85) node[above] {$b$};
\draw[dashed,gray] (2.5,-2.15) rectangle (4,-1.55);
\node at (3.25,-1.85) {open slot};
\node[draw,inner sep=2pt] (di) at (2.85,-1.25) {$x_k^{(s)}$};
\node[draw,inner sep=2pt] (dj) at (3.65,-1.25) {$x_{k+1}^{(s)}$};
\draw (di)--(2.85,-1.55); \draw (dj)--(3.65,-1.55);
\node[anchor=west] at (0,-2.5) {(c) Update the allowed charge blocks};
\node[center] (old) at (1,-3.1) {$B$};
\draw (0.35,-3.1)--(old); \draw (old)--(1.65,-3.1);
\draw (0.85,-2.8)--(0.85,-2.65); \draw (1.15,-2.8)--(1.15,-2.65);
\draw[->] (2,-3.1)--node[above] {$-\eta\nabla_B\mathcal{L}$}(4.1,-3.1);
\node[center] (new) at (5.25,-3.1) {$\widetilde B$};
\draw (4.6,-3.1)--(new); \draw (new)--(5.9,-3.1);
\draw (5.1,-2.8)--(5.1,-2.65); \draw (5.4,-2.8)--(5.4,-2.65);
\node[anchor=west] at (0,-3.85) {(d) Block SVD and split; move the center right};
\node[center] (sv) at (1,-4.9) {$\widetilde B_q$};
\draw (0.3,-4.9)--(sv); \draw (sv)--(1.7,-4.9);
\draw (0.85,-4.6)--(0.85,-4.35); \draw (1.15,-4.6)--(1.15,-4.35);
\draw[->] (2,-4.9)--node[above] {SVD}(3,-4.9);
\node[tn] (u) at (3.7,-4.9) {$U_q$};
\node[center] (vs) at (5.45,-4.9) {$\Sigma_q V_q^T$};
\draw (3.1,-4.9)--(u); \draw (u)--node[above] {$q$}(vs);
\draw (vs)--(6.35,-4.9);
\draw (u)--(3.7,-4.35) node[above] {$i$};
\draw (vs)--(5.45,-4.35) node[above] {$j$};
\end{tikzpicture}
\caption{One two-site training update, following
\citet{LopezPiqueres2023}. Lines are tensor indices; connected lines are
contracted. Orange marks the canonical center. In (b), the physical
legs are fixed to the training sample by Kronecker deltas; the dashed
slot is not a tensor. In (d), factors are reshaped into TT cores.
The leftward update instead puts $\Sigma_q$ into the left core.}
\label{fig:training_sweep}
\end{figure}
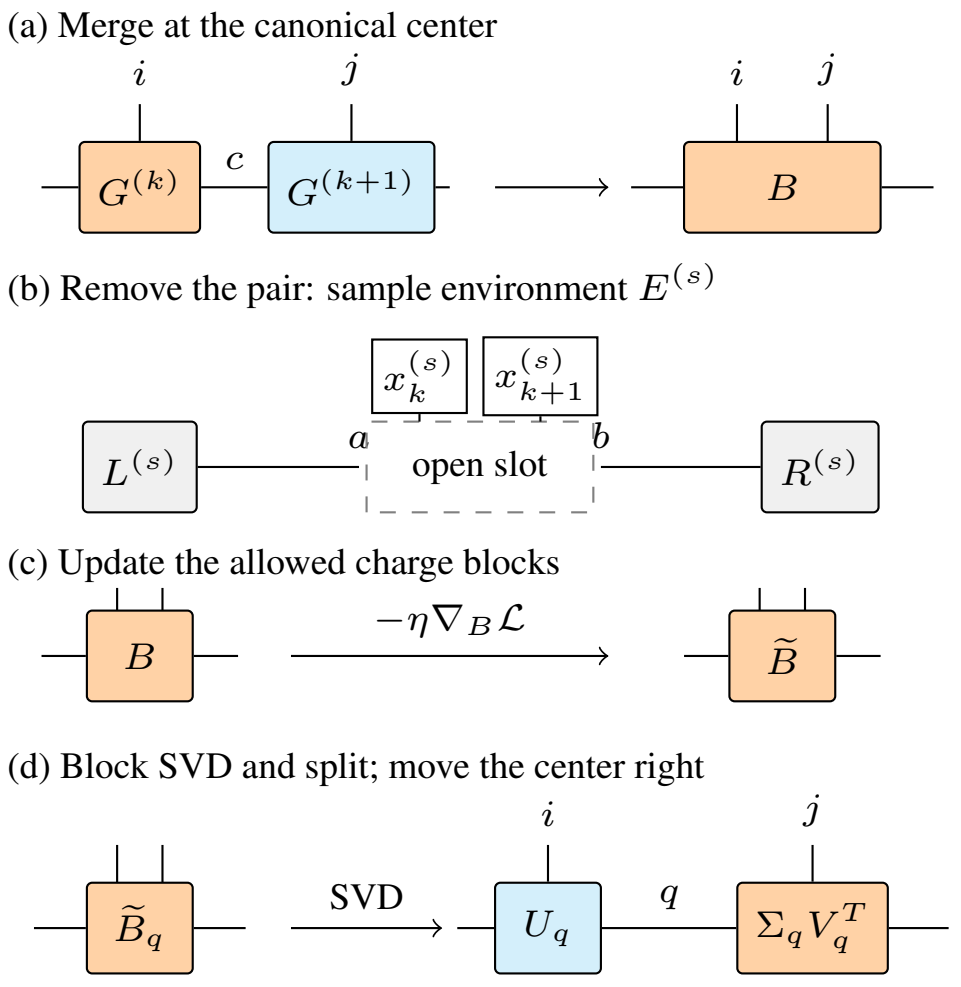

\paragraph{Split and sweep.}
Reshape $\widetilde B$ with row index $(a,i)$ and column index $(j,b)$.
Charge conservation separates this matrix into blocks labeled by the
intermediate residual charge
$q=c_a-i\vec A_k=c_b+j\vec A_{k+1}$.
For each retained charge compute
$\widetilde B_q=U_q\Sigma_q V_q^T$, truncate within the block, and reshape
back to cores. On a rightward pass use $(U_q,\Sigma_q V_q^T)$; on a
leftward pass use $(U_q\Sigma_q,V_q^T)$. Thus the orthogonality center
moves one site and forbidden charge sectors remain zero
(\cref{fig:training_sweep}). One full sweep visits
$k=1,\ldots,N-1$ and then $k=N-1,\ldots,1$.
Repeat for the prescribed number of sweeps, normalize the center core,
and sample sequentially from the trained model.

\paragraph{Specialization used here.}
During a training sweep, our current kernels keep the existing charge sets fixed and retain
only the largest singular value in each block, i.e., virtual charge
degeneracy one. This does not mean total TT rank one: the bond dimension
is the number of retained charges. Consequently, training changes weights
on the existing charge graph. 

\subsection{Expressibility Enhancement Local Search (EELS)}
\label{subsec:eels}
We noticed that all methods for building symmetric TT from data, as suggested in \cite{LopezPiqueres2023}, may fail for the same reason. Methods take link charges only from the data but do not apply any strategy to search for unknown charges. If the training data is too small or not diverse enough, then the required link charges, which encode near-optimal solutions, simply do not appear in TT. Moreover, this disadvantage may not be compensated during the training step, because the charge structure pre-defines all samples, which may be sampled from symmetric TT \textit{in principle}, and the learning procedure updates only the probabilities of principally possible samples.

\begin{figure}[!tbp]
    \centering

    \begin{subfigure}{0.90\linewidth}
        \centering
        \includegraphics[width=\linewidth]{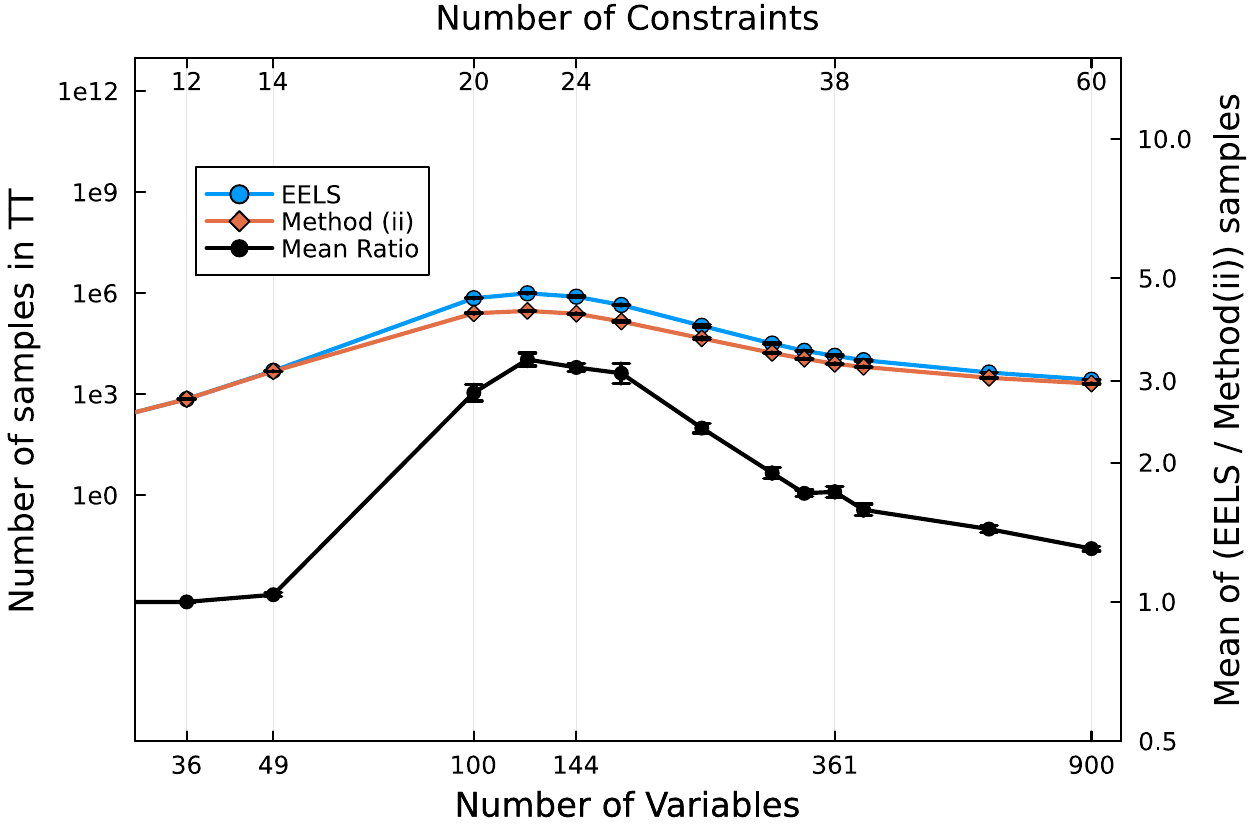}
        \caption{$\left|T\right|=400$.}
        \label{fig:assignemnt_400_abs}
    \end{subfigure}
    \par\vspace{2pt}
    
    \begin{subfigure}{0.90\linewidth}
        \centering
        \includegraphics[width=\linewidth]{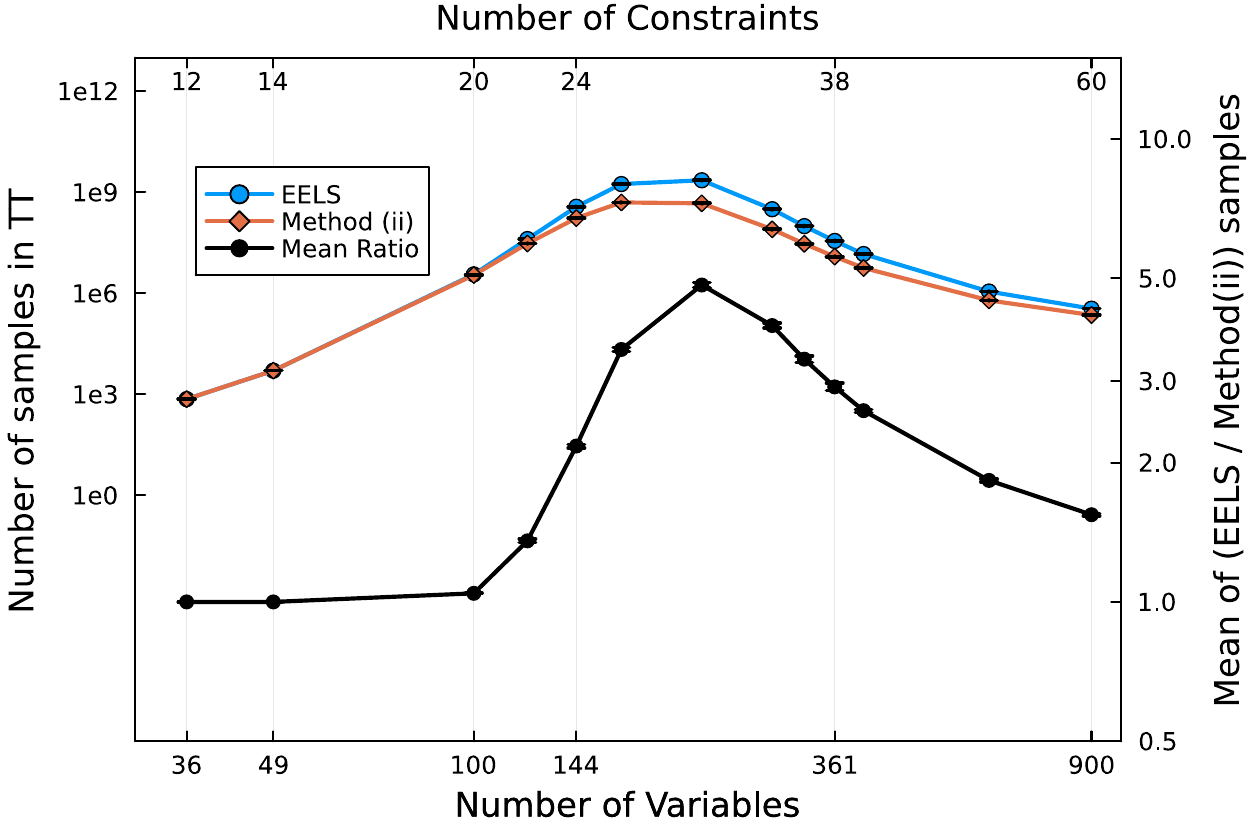}
        \caption{$\left|T\right|=4\mathrm{K}$.}
        \label{fig:assignment_4000_abs}
    \end{subfigure}
    \par\vspace{2pt}

    \begin{subfigure}{0.90\linewidth}
        \centering
        \includegraphics[width=\linewidth]{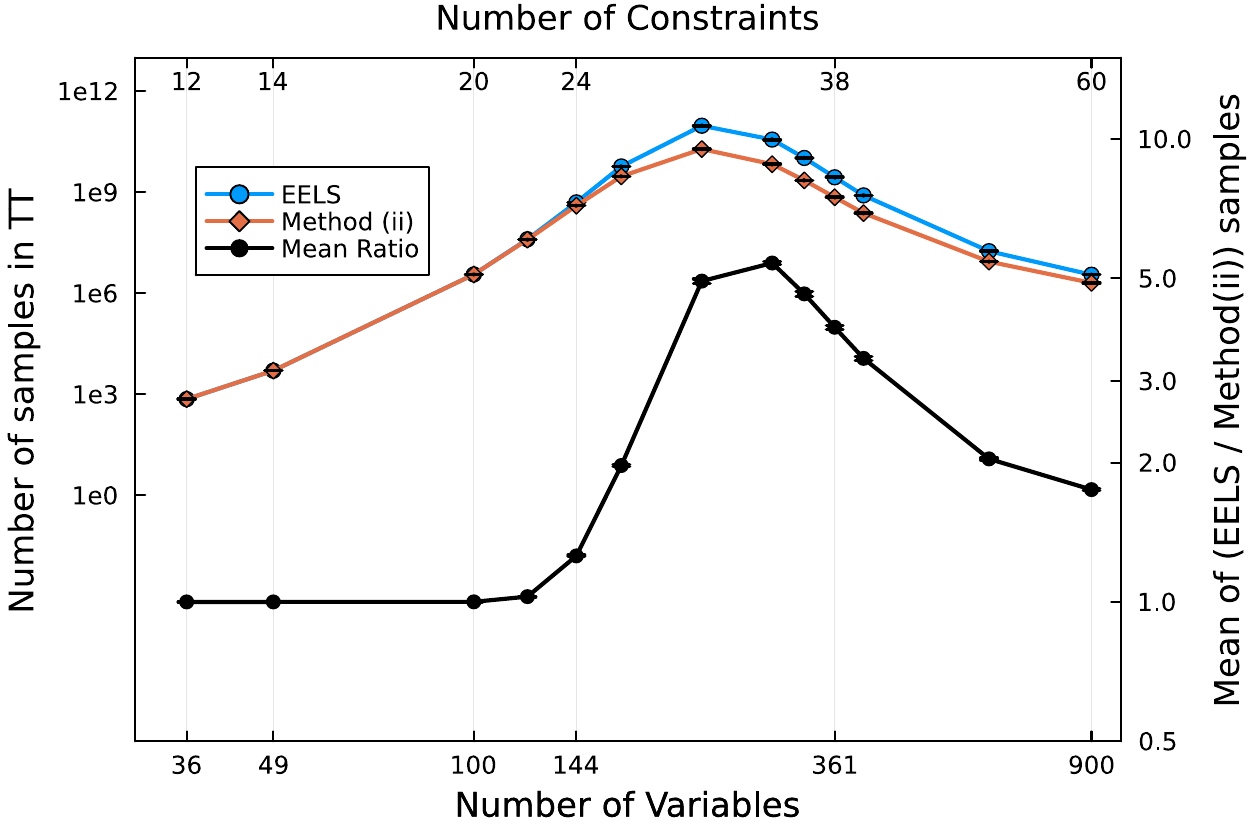}
        \caption{$\left|T\right|=10\mathrm{K}$. }
        \label{fig:assignment_10000_abs}
    \end{subfigure}
    \par\vspace{2pt}

    \caption{Number of unique samples in TT for different $U(1)$-symmetric TT building strategies on the assignment problem.}
    \label{fig:assignment_experiments}
\end{figure}

\begin{figure*}[!tbp]
    \centering

    \begin{subfigure}{0.44\textwidth}
        \centering
        \includegraphics[width=\linewidth]{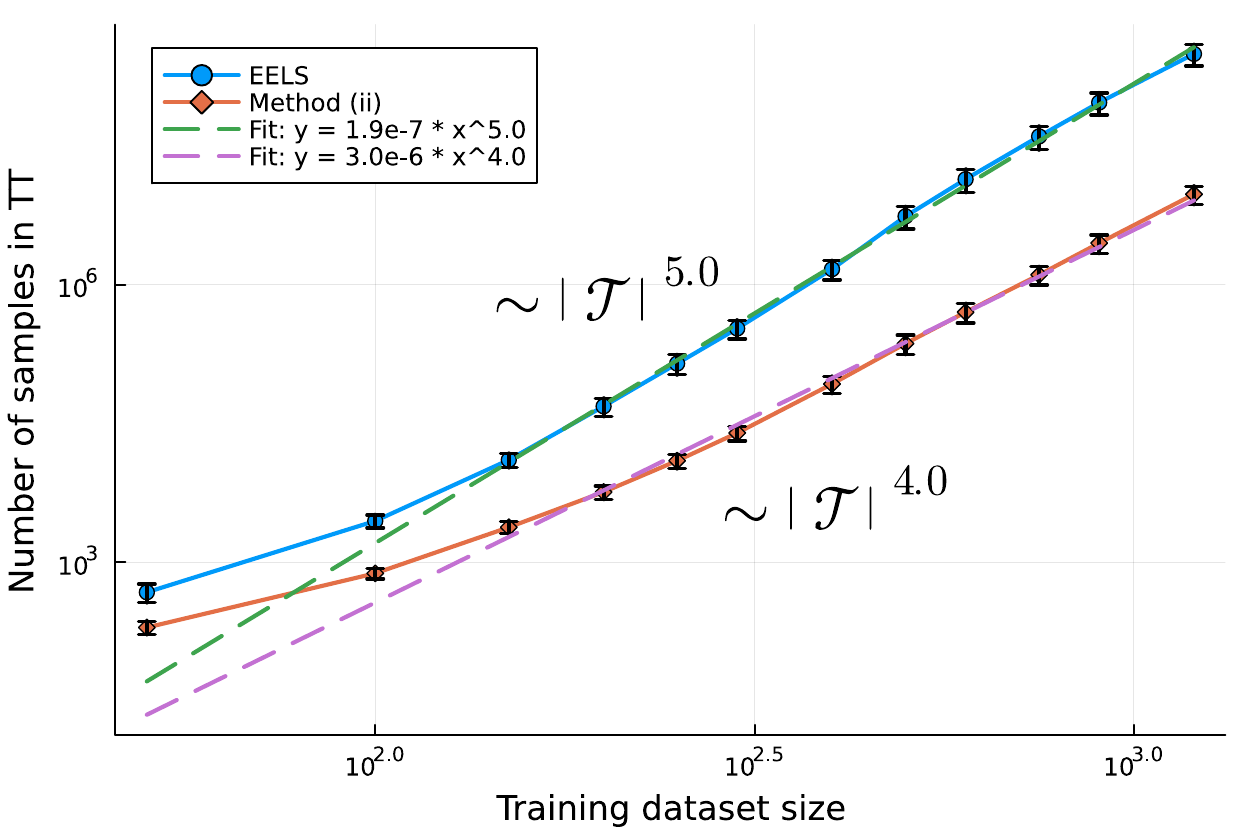}
        \caption{N=50, m=5, r=1.}
        \label{fig:T4_T5}
    \end{subfigure}
    \hspace{10pt}
    \begin{subfigure}{0.44\textwidth}
        \centering
        \includegraphics[width=\linewidth]{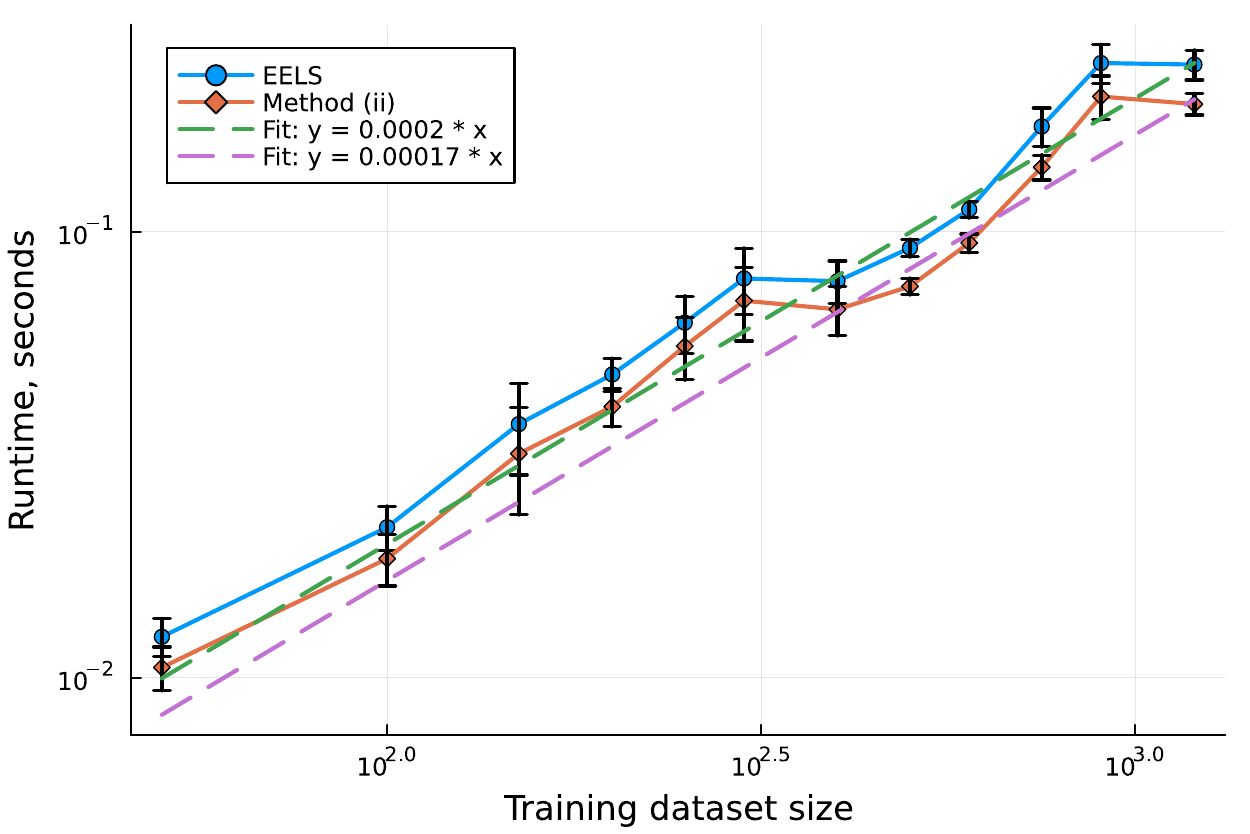}
        \caption{N=50, m=5, r=1.}
        \label{fig:num_samples_time}
    \end{subfigure}
    \hspace{10pt}
    \begin{subfigure}{0.44\textwidth}
        \centering
        \includegraphics[width=\linewidth]{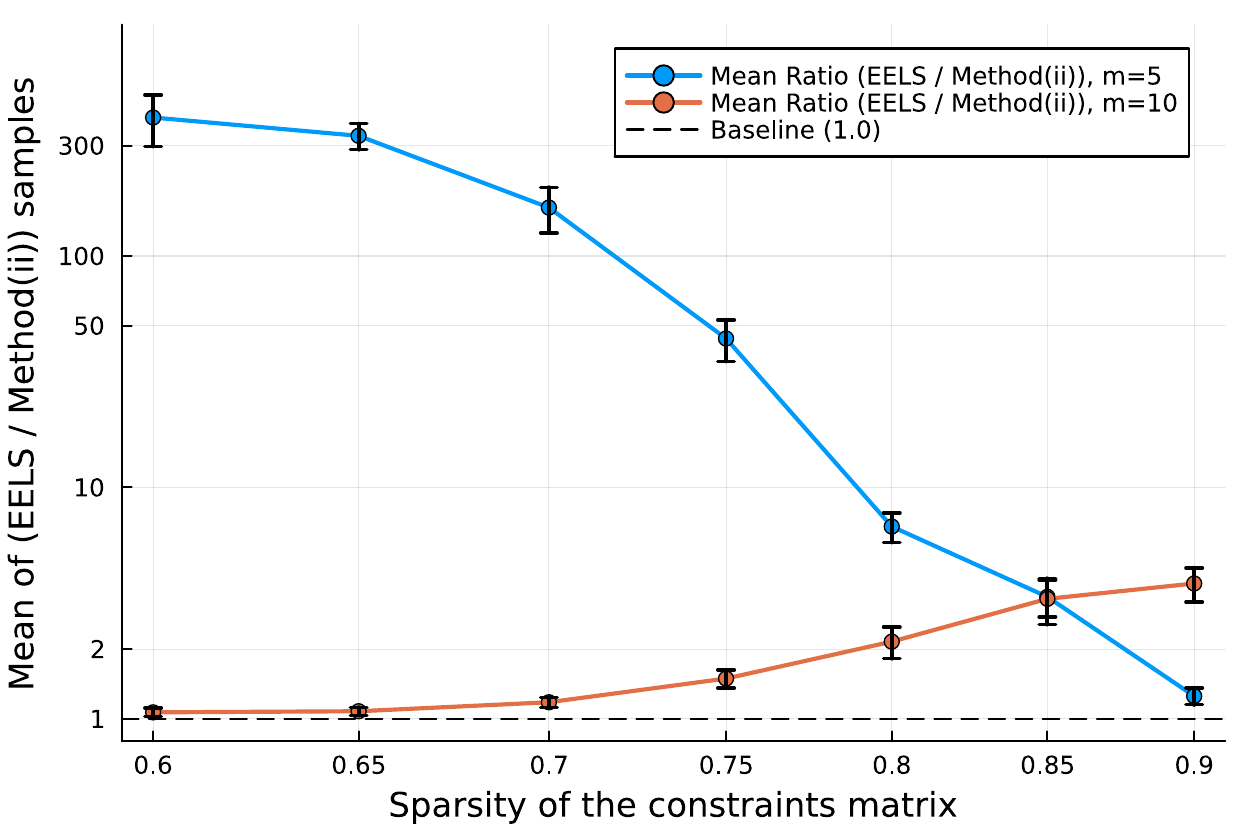}
        \caption{N=75, m=5,10; r=1.}
        \label{fig:random_r1_rel}
    \end{subfigure}
    \hspace{10pt}
    \begin{subfigure}{0.44\textwidth}
        \centering
        \includegraphics[width=\linewidth]{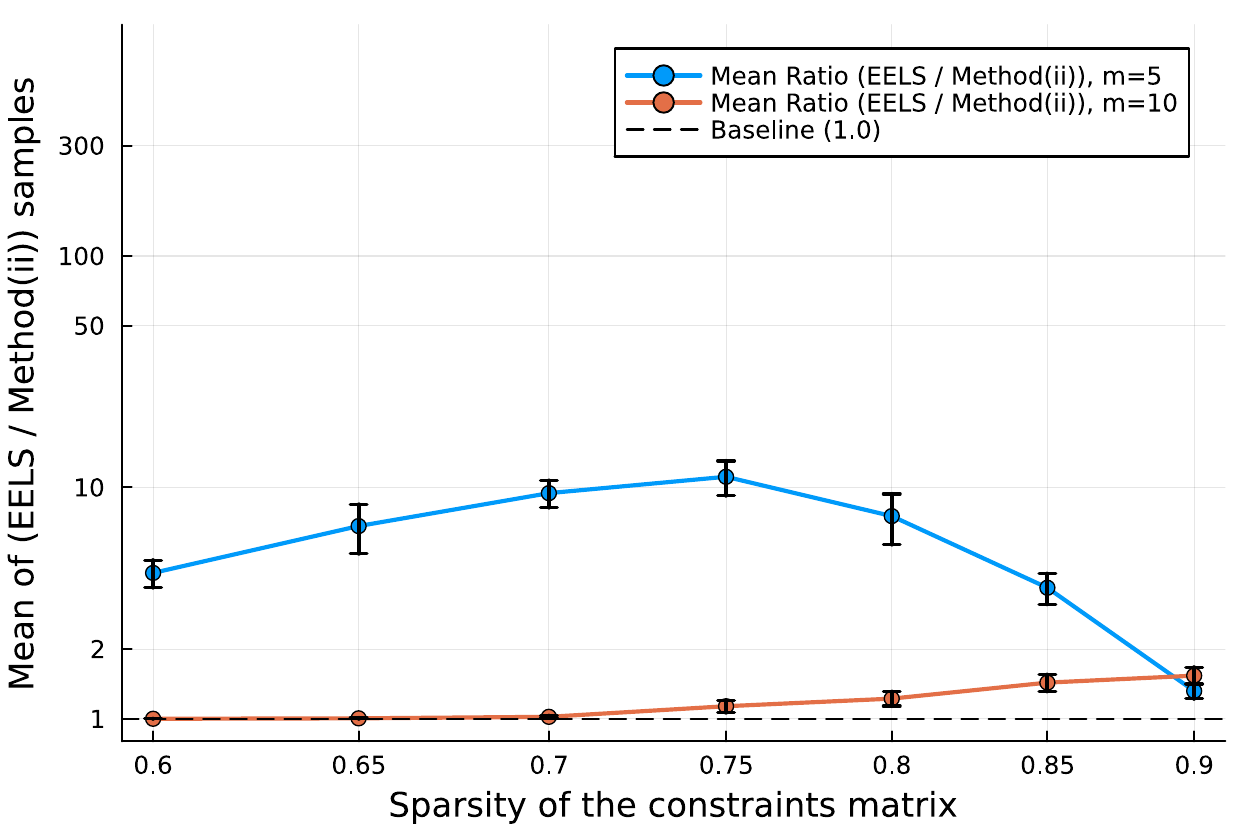}
        \caption{N=75, m=5,10; r=2.}
        \label{fig:random_r2_rel}
    \end{subfigure}

    \caption{Comparison of $U(1)$-symmetric TT building strategies on random problems with $N$ binary variables, $m$ linear equalities and constraints matrix sampled from uniform distribution in range $ [-r, r] $. (a) Number of samples dependency on training dataset size. (b) Wall time to build initial $U(1)$-symmetric TT. (c),(d) Sample-count ratio versus matrix sparsity.}
    \label{fig:num_samples_random}
\end{figure*}

For example, look at the combinatorial problem and dataset $\mathcal{T}$ on \cref{fig:example_graph} (c). Method (ii) fails to find any new samples that do not appear in $\mathcal{T}$. On \cref{fig:example_graph} (a), black circles are the link charges obtained from $\mathcal{T}$, and cyan arrows are all paths in the graph, which encode the dataset samples. Obviously, no new path exists through the black charges, which satisfy the charge conservation law $c_{\alpha_{k+1}} = c_{\alpha_k} - x_{k+1}$ on each edge; thus, no new samples were obtained in TT. 

We propose a new parameterized family of methods for initializing symmetric TT, called Expressibility Enhancement Local Search (EELS). First, we initialize symmetric TT with Method (ii). Then, we apply the expressibility enhancement phase. The \cref{fig:enhancement_scheme} illustrates the main idea of this phase: we fix the $k$-th link of TT, and then search for such charges $c_{new}$, which i) do not appear at the $k+1$-th link; ii) may connect the charges at the $k$-th link with the charges at the $(k+2)$-th link, satisfying charge conservation law for some $x_{k+1}$, $x_{k+2}$. Thus, we "jump" over the $k+1$-th link, and conduct a local search in a charge graph, finding such $c_{\alpha_{k}}$, $x_{k+1}$, $x_{k+2}$, which satisfy the $U(1)$-symmetry. Sometimes this procedure spawns new charges, which may be assigned to the $k+1$-th link. The length of the "jump" may be regulated with the parameter $\gamma$, which is equal to the number of links we jump over. At the \cref{fig:enhancement_scheme}, $\gamma=1$, and the \cref{alg:enhance_expressibility} also provides the detailed description for $\gamma=1$. Method (ii) may be associated with $\gamma=0$. For $\gamma>1$, we search for new charges at the links $k+1$, $k+2$, ..., $k+\gamma$ simultaneously, by considering all possible combinations of $c_{\alpha_k}$, $x_{k+1}$, ..., $x_{k+\gamma + 1}$, thus, the complexity of the algorithm is $O\left(2^{\gamma + 1}|\mathcal{T}| |b|(N-\gamma)\right)$ on average. \cref{fig:example_graph} shows that even for $\gamma=1$, our expressibility enhancement procedure may find some new charges, and thus, new paths in the charge graph, corresponding to feasible samples.

\begin{algorithm}
\caption{EELS algotithm, $\gamma=1$}
\label{alg:enhance_expressibility}
\begin{algorithmic}[1]

\REQUIRE $C_1,\dots, C_{N-1}$  - Sets of link charges, defined from data; $A, b$ - constraints for the problem; $G^{(1)}, \dots ,G^{(N)}$ - TT-cores;

\STATE $c_{\alpha_0}\gets \vec{b}$, $c_{\alpha_N} \gets \vec{0}$

\FOR{$k = 1$ \textbf{to} $N - 2$}
    \FOR{$c_{\alpha_k} \in C_k$}
        \FOR{$x_{k+1} \in \{0, 1\}$}
            \STATE $c_{new} \gets c_{\alpha_k} - x_{k+1}\vec{A_{k+1}}$
            \FOR{$x_{k+2} \in \{0, 1\}$}
                \IF {$c_{new} - x_{k+2}\vec{A_{k+2}} \in C_{k+2}$}
                    \STATE $c_{\alpha_{k+2}} = c_{new} - x_{k+2}\vec{A_{k+2}}$
                    \STATE $C_{k+1}$.add($c_{new}$)
                    \STATE Insert new nonzero blocks:
                    \STATE $G^{(k+1)}\gets (c_{\alpha_{k}}, x_{k+1}\vec{A_{k+1}}, c_{new})$
                    \STATE $G^{(k+2)}\gets (c_{new}, x_{k+2}\vec{A_{k+2}},c_{\alpha_{k+2}})$
                \ENDIF
            \ENDFOR
        \ENDFOR
    \ENDFOR
\ENDFOR
\end{algorithmic}
\end{algorithm}

We emphasize that EELS is not a regular local search, which may be applied to the data samples to obtain new solution candidates. It acts more like a crossover operator, entangling all the samples from $\mathcal{T}$. On \cref{fig:example_graph}(a), EELS obtains the solution, which takes the first bit of the solution $1100$, and the last bit from the solution $0011$, and the bits $00$ in the middle are found using the charge conservation law. Similarly, EELS also finds a feasible solution, taking the first bit from $0011$ and the last bit from $1100$, filling the middle with a guess of $11$. For large datasets, EELS may generate new feasible candidates by mixing all samples from $\mathcal{T}$, producing unpredictable bit concatenations with some "guessed" fillers. 

\subsection{Amplification of Data Diversity (ADD)}
\label{subsec:add}

EELS from \cref{subsec:eels} is useful when we lack training data or want to better exploit existing data. In this section, we consider the special case of constraints $Ax=b$ when a feasible dataset sampler is available, or it is easy to construct a large dataset $\mathcal{T}$ using extra knowledge. Such situations include, for instance, cardinality constraints $\sum_{i=1}^N x_i=\kappa$. In this case, because tensor train optimization trains on the fly, we probably should not use the entire available data but rather pick a good, diverse subset. Otherwise, the complexity of the training step may be (too) large.

We formulate a problem as follows. Let $\mathcal{T}$ be a data, available at the $t$-th iteration of the \cref{alg:tt_generative}. This data includes samples from TT and other feasible solutions generated by the feasible data sampler. How to choose a subset $\mathcal{T}'\subset \mathcal{T}$, such that $|\mathcal{T}'| = \mu < |\mathcal{T}|$, to achieve exploration-exploitation trade-off, and keep the complexity of the training step manageable? We need to define the measure of the sample's novelty relative to what the generative scheme has already explored. To do this, we give a formal definition of the charge graph for the problem \eqref{eq:problem_statement}, which we discuss informally in \cref{subsec:eels}.
Let $|b|=m$ be the number of equalities in problem \eqref{eq:problem_statement} and $N$ the number of binary variables. Let consider a DAG $D=(V, E, \omega, \xi)$, with $\omega: V\to \mathbb{R}^m$ and $\xi: E \to \{\{0\}, \{1\}, \{0,1 \}\}$. Here $\omega(v)$ is a charge of the vertex $v$, and the set of all vertices consists of several layers: $V = L_0 \sqcup L_1 \sqcup \dots \sqcup L_N.$
We call such a graph a charge graph of problem \eqref{eq:problem_statement} if all of the following stand:

\begin{enumerate}
    \item $|L_0| = |L_N| = 1$, $\omega(v_0\in L_0) = \vec{b}$, $\omega(v_N\in L_N)=\vec{0}$
    \item $(u, v) \in E \Rightarrow u\in L_i, \: v \in L_{i+1}$ for some $i$ 
    \item If $u_k\in L_k$, and $v_{k+1}\in L_{k+1}$, then $\omega(u_k)-\omega(v_{k+1})=\vec{0}$ or $\omega(u_k)-\omega(v_{k+1})=\vec{A_{k+1}}$, where $\vec{A_{k+1}}$ is a $k$-th column of the matrix $A$
    \item $\xi(e(u_k, v_{k+1})) = \{0, 1\}$, if $\vec{A_{k+1}} = \vec{0}$, else $\xi(e(u_k, v_{k+1})) = \{0\}$ if $\omega(u_k)=\omega(v_{k+1})$, or $\xi(e(u_k, v_{k+1})) = \{1\}$ if $\omega(u_k)-\omega(v_{k+1})=\vec{A_{k+1}}$.
\end{enumerate}

Charge graph on \cref{fig:example_graph}(a) satisfies all conditions from the definition. On each path $p$ with edges $(e_1, \dots e_N)$ from $v_0 \in L_0$ to $v_N \in L_N$, we define the operation of feasible subset evaluation: $F(p) = \xi(e_1) \times \dots \times \xi(e_N)$, where "$\times$" is the Cartesian product. We call a charge graph \textit{full} if the set $\mathcal{F}(D)=\bigsqcup_{p\in \mathcal{P}(\mathcal{D})}F(p)$ coincides with all the solutions of $Ax=b$, where $\mathcal{P}(\mathcal{D})$ is a set of all the paths from $v_0$ to $v_N$. We define the width of the charge graph as $|D|=\max_{k=0}^{N}\left(|L_k|\right)$. 

The search for the full graph charge may be exponentially hard, and the same is true for storing such a graph. That is why Method (ii) from \cref{subsec:U1_generative} builds only a subgraph $H\subset D$ with a limited width $|H|\leq |\mathcal{T}|$. The goal of the $U(1)$-generative approach is to find such $H\subset D$, that $x^{*} \in \mathcal{F}(H)$ for the optimal solution $x^{*}$ of problem \eqref{eq:problem_statement}. For the ADD selection-and-rebuild variant of \cref{alg:tt_generative}, the graph $H$ is rebuilt from the selected feasible data after external candidates are introduced. This makes the choice of samples control the charge topology as well as the training weights. The articles \cite{Alcazar2024} and \cite{LopezPiqueres2023} suggest only the selection of samples with the best cost function $C(x)$, or picking up some rare samples from $S$ for diversity. We suggest using graph distance as a diversity measure. Let $\mathcal{S}=(V^{\mathcal{S}}, E^{\mathcal{S}}, \omega^{\mathcal{S}}, \xi^{\mathcal{S}})$ be a charge graph for some feasible solution $s=(x_1,\dots, x_N)$. Let $V^{\mathcal{S}}=\bigsqcup_k \{v_k\}$. On each step of the generative optimization, we store a charge graph of already seen samples $\mathcal{H}$, as the sets of charges for each vertex layer: $C_k^{\mathcal{H}} = \bigsqcup_{v\in L_k^{\mathcal{H}}} \{\omega_{\mathcal{H}}(v)\}$. Then the graph distance between $\mathcal{S}$ and $\mathcal{H}$ is:

\begin{equation}
    \rho(\mathcal{S}, \mathcal{H}) = \sum_{k=1}^{N-1}\left(1 - \chi(\omega^{\mathcal{S}} (v_k) \in C_k^{\mathcal{H}})\right),
\label{eq:graph_distance}
\end{equation}

where $\chi(\cdot)$ is an indicator function of the condition in brackets. Thus, we count the number of new charges that appear in the solution $s$ and were not previously seen by the algorithm. In \cref{sec:results} we discuss how to balance this novelty measure with exploitation of the near-optimal samples.

\begin{table*}[!tbp]
\centering
\caption{Persistent-TT EELS versus Method (ii), with Best Cost selection in
both arms: 250 paired instances per objective. W/T/L denotes EELS wins,
exact ties, and losses. Means and medians refer to the paired relative
improvement $R$ (\%). Holm correction covers the two one-sided tests;
percentages are rounded to one decimal place.}
\label{tab:eels_optimization}
\begin{tabular}{lrrrrr}
\toprule
Objective & W/T/L & Mean (\%) & Median (\%) & Raw $p$ & Holm $p$ \\
\midrule
Linear & 125/116/9 & 3.1 & 0.0 & $2\times10^{-19}$ & $3.9\times10^{-19}$ \\
Quadratic & 128/100/22 & 3.1 & 0.2 & $1.4\times10^{-15}$ & $1.4\times10^{-15}$ \\
\bottomrule
\end{tabular}
\end{table*}

\begin{figure*}[!tbp]
\centering
\begin{subfigure}{0.48\textwidth}
\includegraphics[width=\linewidth]{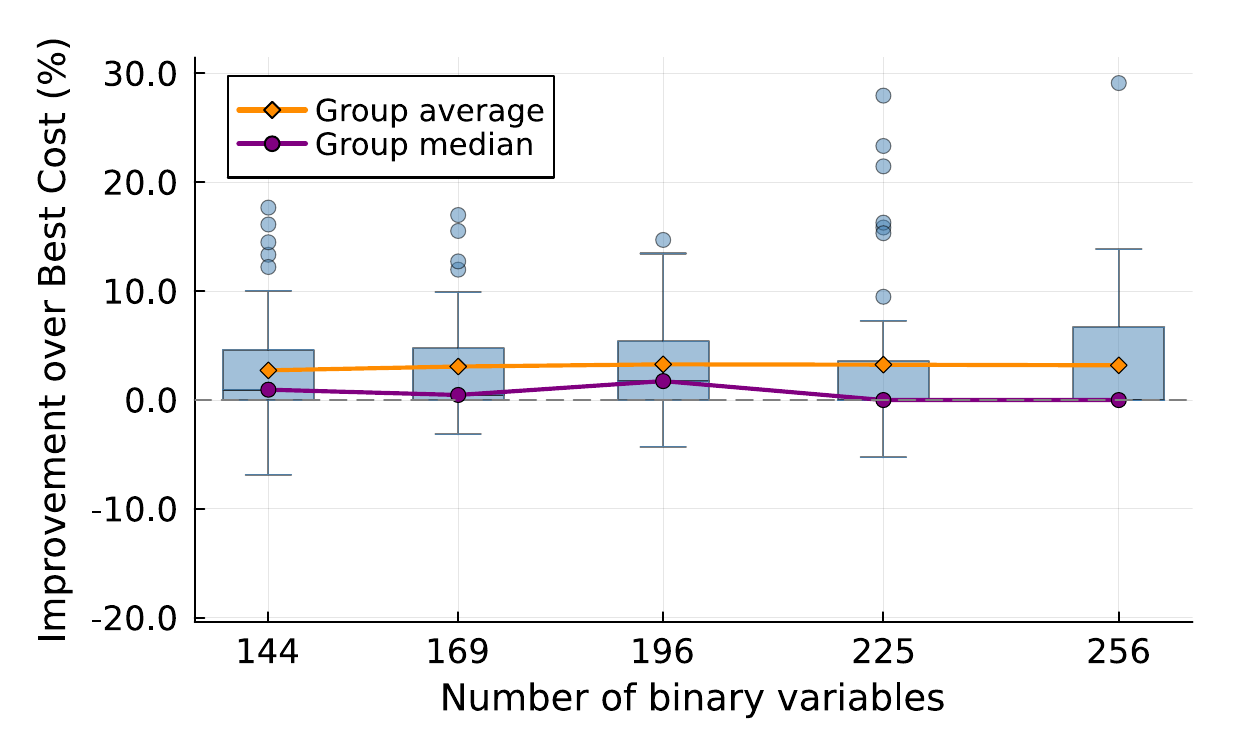}
\caption{Linear objective.}
\end{subfigure}\hfill
\begin{subfigure}{0.48\textwidth}
\includegraphics[width=\linewidth]{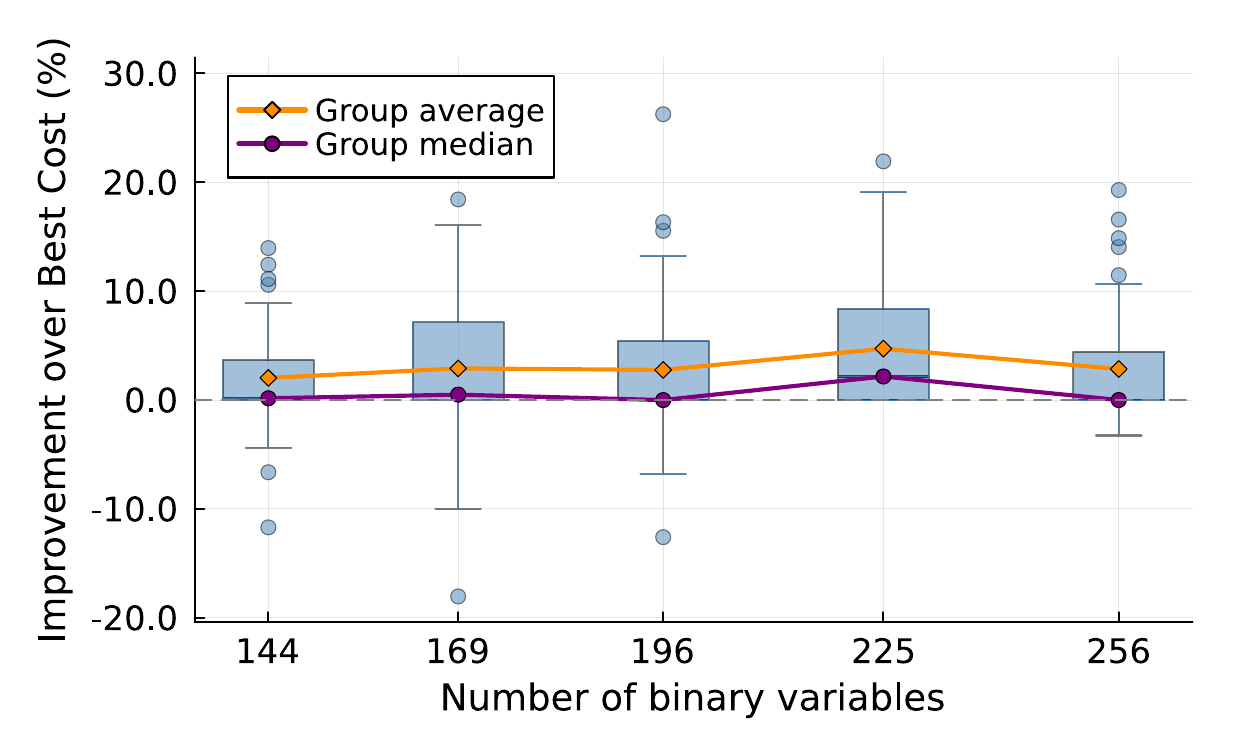}
\caption{Quadratic objective.}
\end{subfigure}
\caption{EELS ($\gamma=1$) versus Method (ii), both with Best Cost selection
and persistent learned amplitudes. Each box contains 50 paired instances.
Purple curves mark group medians and orange curves group arithmetic means.
Box edges are the 25th and 75th percentiles; whiskers extend to the most
extreme observations within 1.5 interquartile ranges of the box. Outliers
are shown individually. Both panels use identical axes. Positive values
favor EELS.}
\label{fig:eels_optimization}
\end{figure*}

\section{Results}
\label{sec:results}
We used DeepSeek-V4.1-Flash to help write the scripts that plot benchmark results from the saved experiment outputs.

\subsection{Comparison of EELS with other methods}
\label{subsec:eels_res}
We may compare the expressibility of the methods of building $U(1)$-symmetric tensor trains by fixing the dataset $\mathcal{T}$, defining all link charges degeneracy ranks in \eqref{eq:direct_sum} $\dim(V^{(c_k)}) = 1$ and then, calculating the number of nonzero elements. 
If we initialize $G^{(k)}(\alpha_{k-1}, x_k, \alpha_k)= 1$ for all positions, preserving $U(1)$-symmetry, then $\pi(x)=1$ for all $x$, conserving charge $b$. Thus, the number of unique samples in tensor train is just a Frobenius norm squared, which may be calculated efficiently in TT-format.

\cref{fig:num_samples_random} shows the group of experiments for the problem \eqref{eq:problem_statement} with the random matrices $A$ of size $m\times n$ and random vector $b$, both sampled from the uniform distribution on $[-r, r]$. \cref{fig:T4_T5} shows the dependency of the number of unique samples in $U(1)$-symmetric TT as the function of training dataset size. With EELS, the number of unique samples grows as $\sim|\mathcal{T}|^5$ for $r=1$, while applying only the baseline method (Method (ii) from \cref{subsec:U1_generative}) encodes only $\sim|\mathcal{T}|^4$ possible feasible solutions. EELS improves the baseline with a small increase in computational time to build an initial $U(1)$-symmetric TT, as shown in \cref{fig:num_samples_time}. Other experiments on \cref{fig:num_samples_random} demonstrate that the benefit of using EELS depends a lot on the sparsity of the matrix $A$, the number of constraints, and the value of parameter $r$. Generally, if $r$ is big, or the number of constraints is huge, then it is hard to satisfy the if-condition in line 7 of \cref{alg:enhance_expressibility}, and EELS spawns only a small fraction of new link charges $c_{new}$, if it spawns any at all. 
In contrast, if the matrix $A$ is sparse, EELS may be an efficient strategy for several dozen constraints. For example, we consider an assignment problem \cite{koopmans1957assignment}. This is NP-hard problem for a quadratic/black-box oracle  $C(x)$. The constraints matrix of size $2n \times n^2$ is given by the rule:

\begin{equation}
\left\{
\begin{aligned}
    & A_{j, (i-1)n+j} = 1, \quad \forall i,j = 1, \dots, n\\[4pt]
    & A_{n+i, (i-1)n+j} = 1, \quad \forall i,j = 1, \dots, n\\[4pt]
    & \text{else} \quad A_{k, m} = 0,
\end{aligned}
\right.
\label{eq:assignment_matrix}
\end{equation}
and $b_i=1 \quad \forall i = 1, \dots ,2n$. The density of matrix $A$ scales as $O(\frac{1}{n})$. \cref{fig:assignment_experiments} shows the group of experiments with the assignment problem. 
The total number of feasible samples grows monotonically with the problem size as $n!$; however, for large problems, building symmetric  TT methods struggle to find new feasible samples when the training data is insufficiently expressive. That is why the number of samples in TT saturates and then decreases. The peak of saturation shifts to the right as the initial dataset size grows. EELS provides approximately $3-6$ times more samples at the peak point than Method (ii) from \cref{subsec:U1_generative}.

\subsection{EELS impact on optimization results}
\label{subsec:eels_optimization}
To test whether enhanced expressibility translates into lower objective values,
we compare EELS with $\gamma=1$ against Method (ii), using Best Cost selection
in both arms. For each of the linear and quadratic objective families, we use
50 assignment instances at each $n\in\{12,13,14,15,16\}$: 250 paired comparisons
with $N=n^2\in\{144,169,196,225,256\}$ binary variables.

For each pair, 400 feasible assignments are drawn uniformly from random
permutations. Both arms receive the same deduplicated, cost-sorted initial
matrix.
Both initialize their TT once, with EELS applied only at initialization in
the enhanced arm. Each of 20 outer iterations first draws $10{,}000$ samples, then trains the existing TT on the best at most
400 distinct assignments from that new batch alone. Training uses one complete
right--left sweep, learning rate $0.05$, and virtual charge degeneracy one.
Boltzmann weights use the standard deviation of the current elite costs as
temperature, with uniform weights when this deviation is zero. Learned
amplitudes and structural charge sets persist across iterations. Initial and
historical points are not merged into subsequent training batches; the best
evaluated solution is retained separately. No external data, restarts, or
subsequent TT reinitializations are used. 

For best recorded costs $C_{\mathrm{ii}}$ and $C_{\mathrm{EELS}}$, define
\begin{equation}
 R=100\,\frac{C_{\mathrm{ii}}-C_{\mathrm{EELS}}}{|C_{\mathrm{ii}}|}.
 \label{eq:eels_relative_improvement}
\end{equation}
All baseline costs are nonzero; positive $R$ favors EELS.
\Cref{fig:eels_optimization} shows the distributions by size, including
group medians and arithmetic means. A one-sided
paired Wilcoxon signed-rank test is applied to unscaled final costs, pooled
over the five sizes within each objective. \Cref{tab:eels_optimization}
reports raw $p$-values and Holm adjustment across the two objective-level
comparisons.

The mean improvement is 3.1\% for linear and 3.1\%
for quadratic objectives, with 125 and 128 wins out
of 250 pairs, respectively. Both advantages remain significant after the
Holm correction. However, there are 116 and 100
exact ties, and the pooled medians are only 0.0\% and
0.2\% after rounding. Thus the improvement is heterogeneous,
not a uniform gain on every instance. These results
support the one-time expressibility enhancement under equal sampling and
iteration settings. Improvements also persist at a fixed 200\,s budget
(\cref{fig:eels_equal_time,subsec:equal_time}).

\begin{figure*}[!tbp]
\centering
\includegraphics[width=\textwidth]{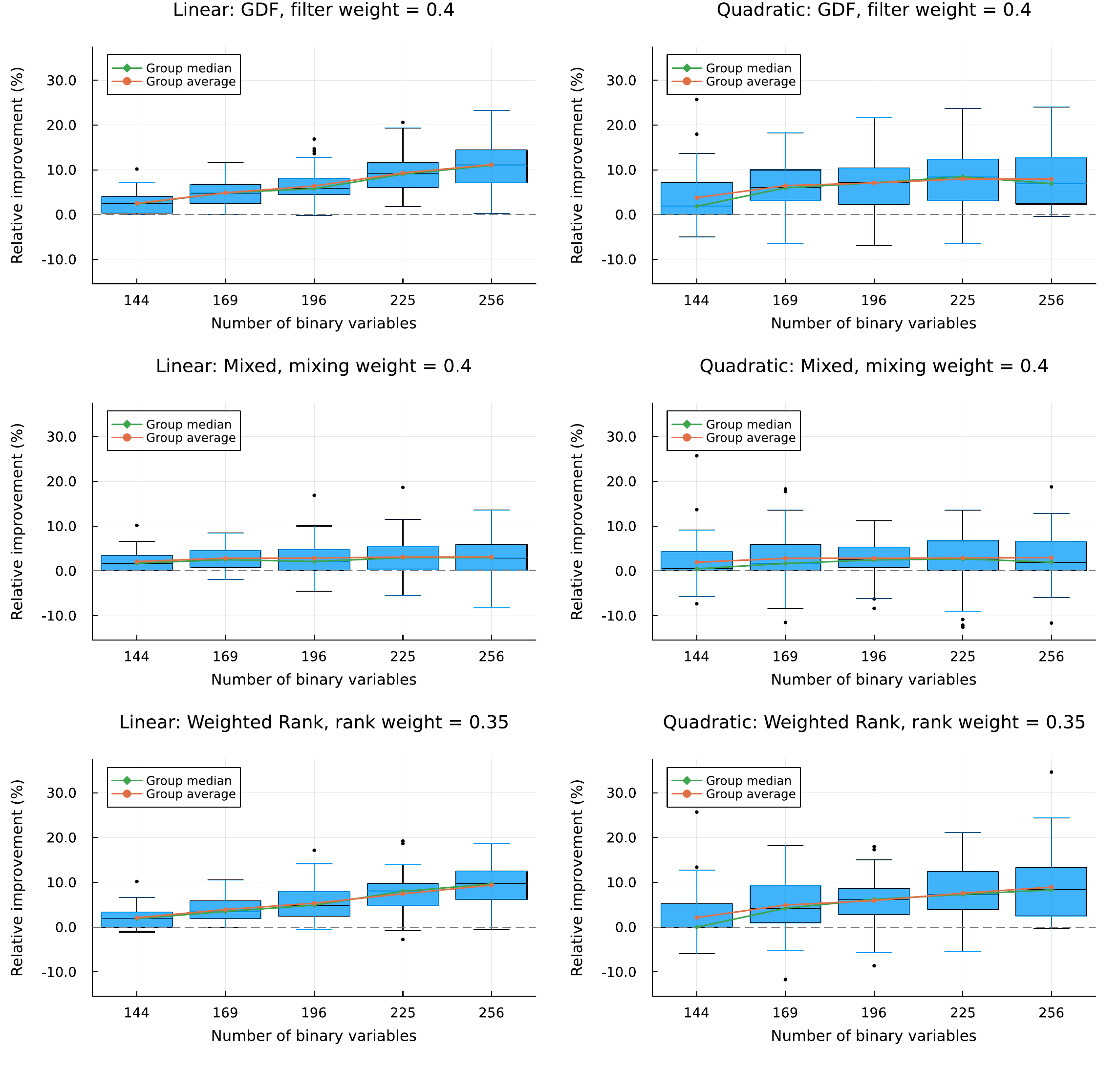}
\caption{ADD versus Best Cost on the parameter-screening instances.
Columns show linear (left) and quadratic (right) objectives; rows show
GDF ($w_{\mathrm{filter}}=0.4$), Mixed Selection ($w_{\mathrm{mix}}=0.4$), and Weighted Rank ($w_{\mathrm{rank}}=0.35$).
Each box contains 50 paired relative improvements at the indicated number
of binary variables. Positive values favor ADD. Box edges are the first
and third quartiles; whiskers extend to observations within 1.5 interquartile
ranges and remaining observations are shown as points. Green curves with
diamond markers and orange curves with circle markers show group medians
and arithmetic means, respectively.
All panels share axis limits. Each displayed weight maximizes mean relative
improvement within its strategy on these same instances.}
\label{fig:add_selected_strategies}
\end{figure*}

\begin{figure*}[!tbp]
\centering
\includegraphics[width=0.49\textwidth]{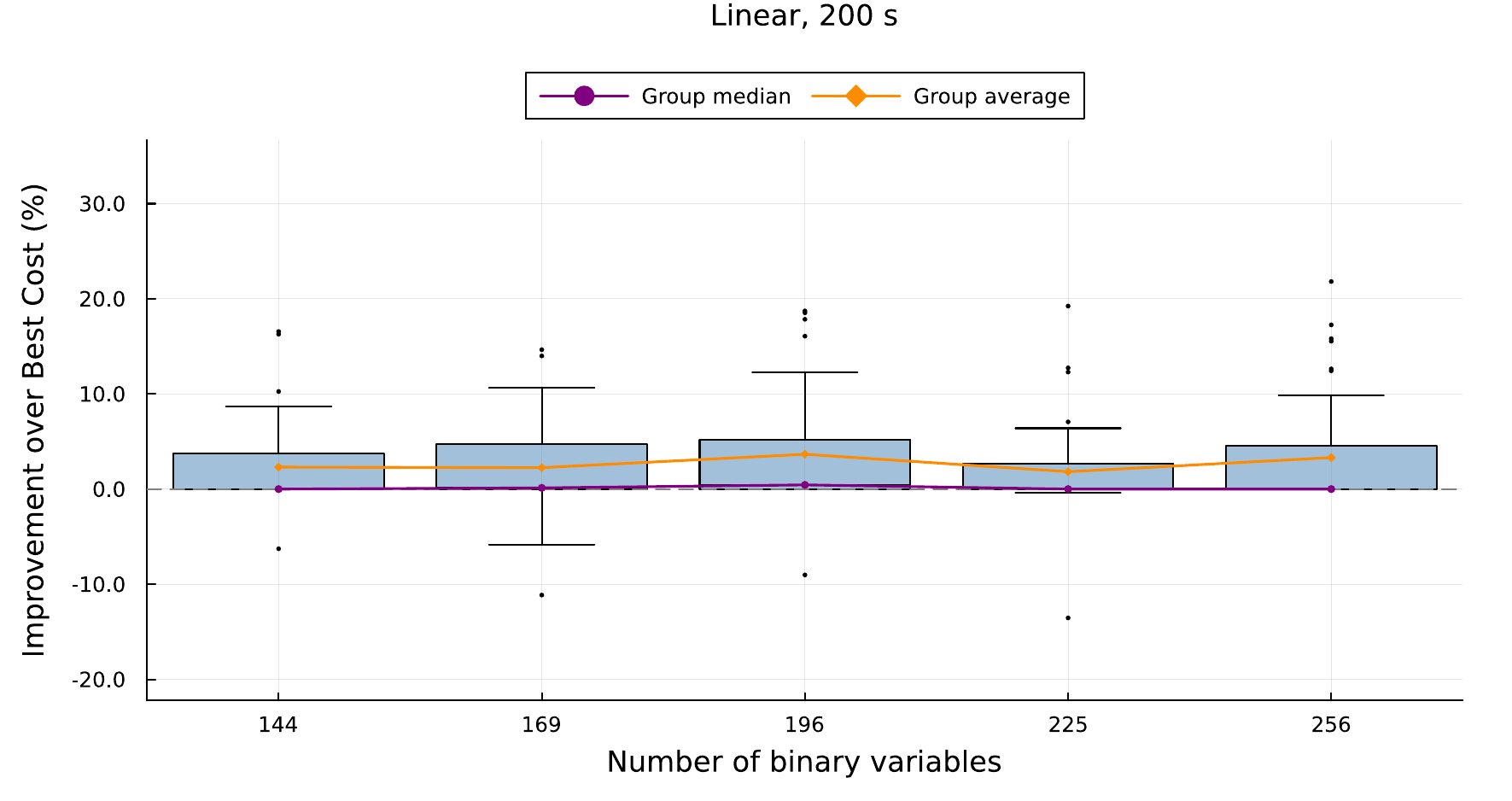}
\includegraphics[width=0.49\textwidth]{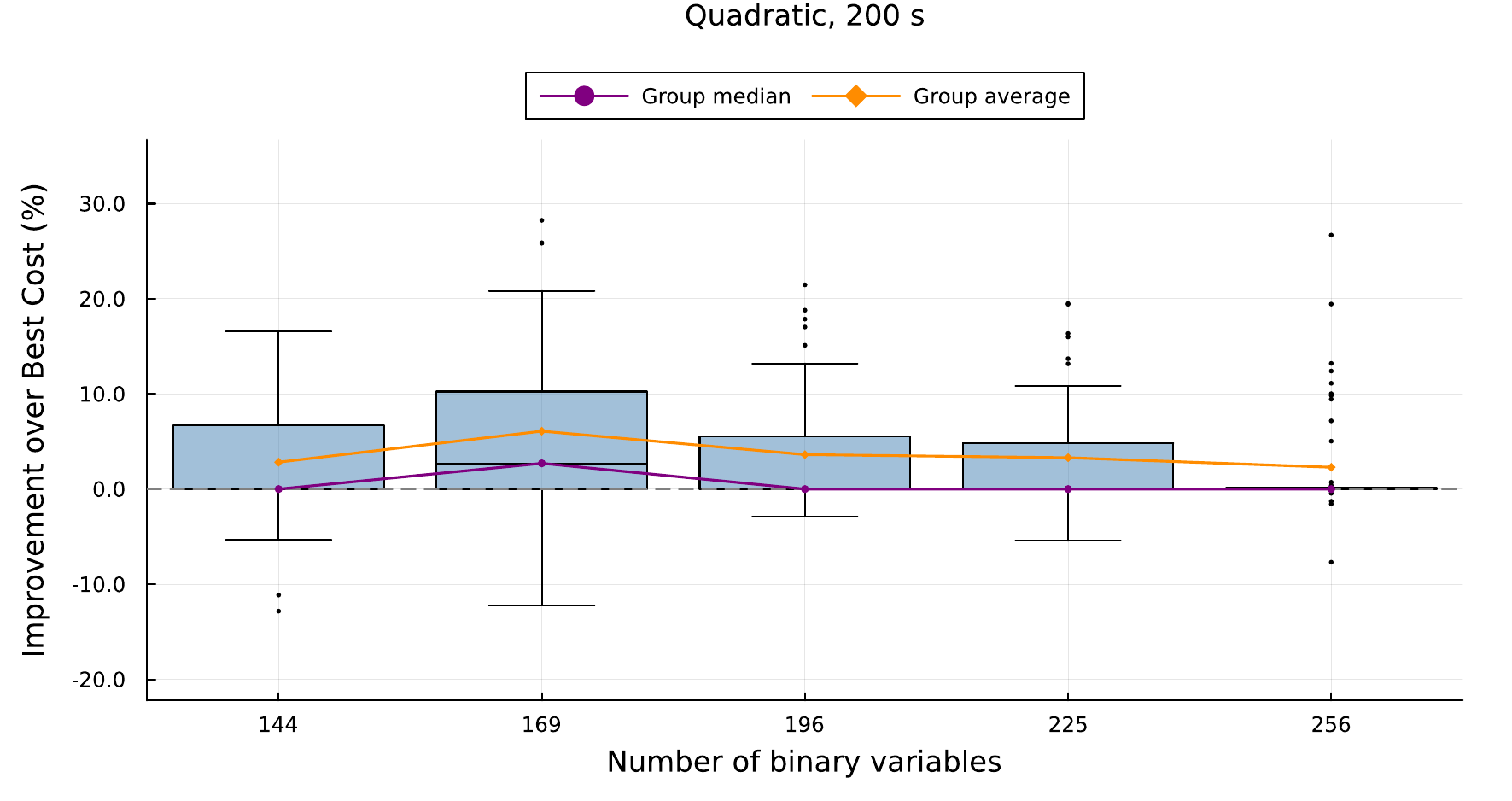}
\caption{EELS ($\gamma=1$) versus Best Cost with a 200\,s budget per arm: linear
(left) and quadratic (right) assignment, 50 paired instances per size.
Boxes show quartiles and medians; whiskers extend to observations within
1.5 interquartile ranges, with other observations plotted individually.
The labelled curves show group medians and means. Positive relative
improvement favors EELS. Both panels share axis limits. The training
protocol and paired tests are given in \cref{subsec:equal_time,tab:equal_time}.}
\label{fig:eels_equal_time}
\end{figure*}

\subsection {ADD impact on optimization results}
\label{subsec:diversify_data_res}
We use graph distance (\cref{eq:graph_distance}) to balance cost-based
selection with exploration of the charge graph. Assignment admits a simple
uniform feasible sampler based on random permutations, so this experiment
examines which available data to learn from, rather than how to obtain
feasible data. In the ADD variant of \cref{alg:tt_generative},
we reconstruct the TT from the selected data at every outer iteration.
We compare the following selection rules on a candidate pool $S$, with
selection quota $\mu$ and charge history $\mathcal H$:
\begin{enumerate}
\item \textit{Best Cost}: select the $\mu$ candidates with lowest cost.
\item \textit{Best Cost + Graph Distance Filter (GDF)}: first select
$(1-w_{\mathrm{filter}})\mu$ candidates by cost, then select up to $w_{\mathrm{filter}}\mu$ lowest-cost
remaining candidates satisfying $\rho(\mathcal S,\mathcal H)>0$.
If no such candidates remain, use Best Cost as a fallback.
\item \textit{Mixed Selection}: select $w_{\mathrm{mix}}\mu$ candidates with largest
graph distance and fill the remaining quota by cost from the other candidates.
\item \textit{Weighted Rank}: select the $\mu$ candidates with smallest
$r_C+w_{\mathrm{rank}} r_\rho$, where rank one denotes lowest cost or largest graph
distance, respectively.
\end{enumerate}
The history stores charge sets at each TT link. Computing one graph distance
requires $O(|b|(N-1))$ operations on average with hashed charge sets;
history storage is an additional cost of the diversity-aware rules.

For each of the linear and quadratic objectives we use 50 instances at
each $n\in\{12,13,14,15,16\}$, with $N=n^2$ binary variables.
Every configuration receives the same 1K initial feasible assignments
for each instance, verified by stored hashes. Each run performs 20 outer
iterations, one training sweep per iteration, learning rate $0.05$ and
bond degeneracy one. After training a newly initialized TT, we draw
10K TT samples, deduplicate them together with 50 retained elite
assignments, and supplement the pool with fresh uniform feasible samples
to reach 20K candidates. Each rule selects up to $\mu=10^3$ candidates;
retaining the 50 elites can increase the next training set to at most
1.05K distinct assignments. No worst-sample retention, EELS enhancement
or barrier filter is used. All methods, including Best Cost, receive
external candidates.

We evaluate $w_{\mathrm{mix}},w_{\mathrm{filter}}\in\{0.1,0.2,0.3,0.4\}$ for Mixed Selection and GDF,
and $w_{\mathrm{rank}}\in\{0.1,0.35,0.5,0.7\}$ for Weighted Rank.
For each strategy and objective, we select the weight with the largest
mean per-instance relative improvement
$R=100(C_{\mathrm{BestCost}}-C_{\mathrm{ADD}})/|C_{\mathrm{BestCost}}|$.
\Cref{fig:add_selected_strategies} shows all six selected configurations;
the complete grid and paired statistics are in
Supplementary Section~S1.

All three strategies improve on Best Cost at their selected weights.
The mean gains for Mixed Selection, Weighted Rank and GDF are,
respectively, $2.8\%$, $5.6\%$ and $6.9\%$ for linear objectives,
and $2.7\%$, $5.9\%$ and $6.6\%$ for quadratic objectives.
GDF with $w_{\mathrm{filter}}=0.4$ gives the largest mean improvement in both families.
It wins, ties and loses on 235/14/1 linear instances and 209/26/15
quadratic instances; its pooled median improvements are $5.9\%$ and
$6.2\%$. One-sided paired Wilcoxon tests on raw final costs give
$p=9.5\times10^{-41}$ and $p=2\times10^{-34}$, respectively;
Holm adjustment over all 24 baseline comparisons gives
$2.2\times10^{-39}$ and $3.5\times10^{-33}$.
The mean linear improvement of GDF increases from $2.5\%$ at $n=12$
to $11.2\%$ at $n=16$; its quadratic mean improves from $3.8\%$ to
$7.9\%$. This pattern is consistent with a benefit from charge-space
exploration as problem size grows.

\begin{figure*}[!tbp]
\centering
\includegraphics[width=0.49\textwidth]{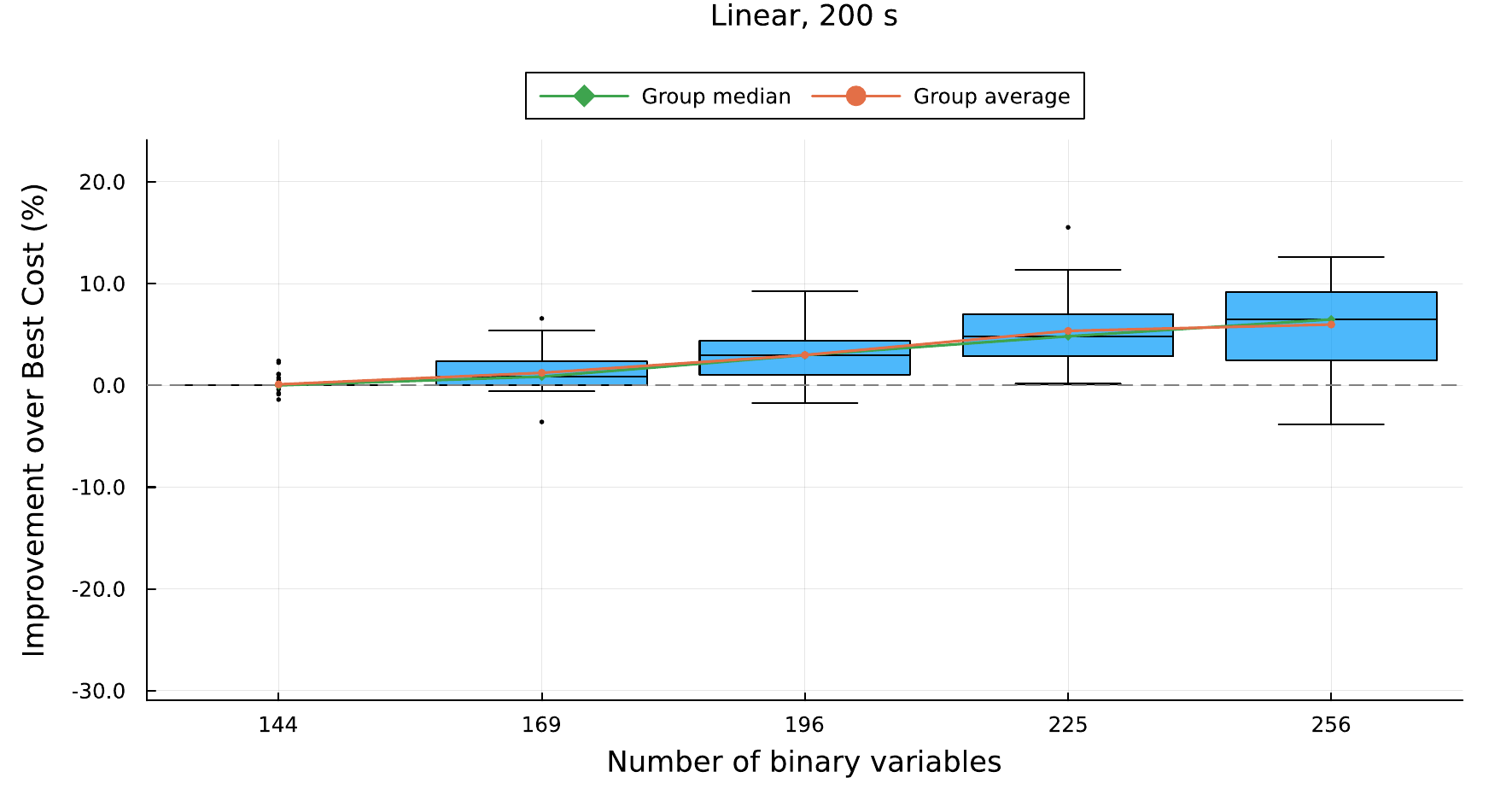}
\includegraphics[width=0.49\textwidth]{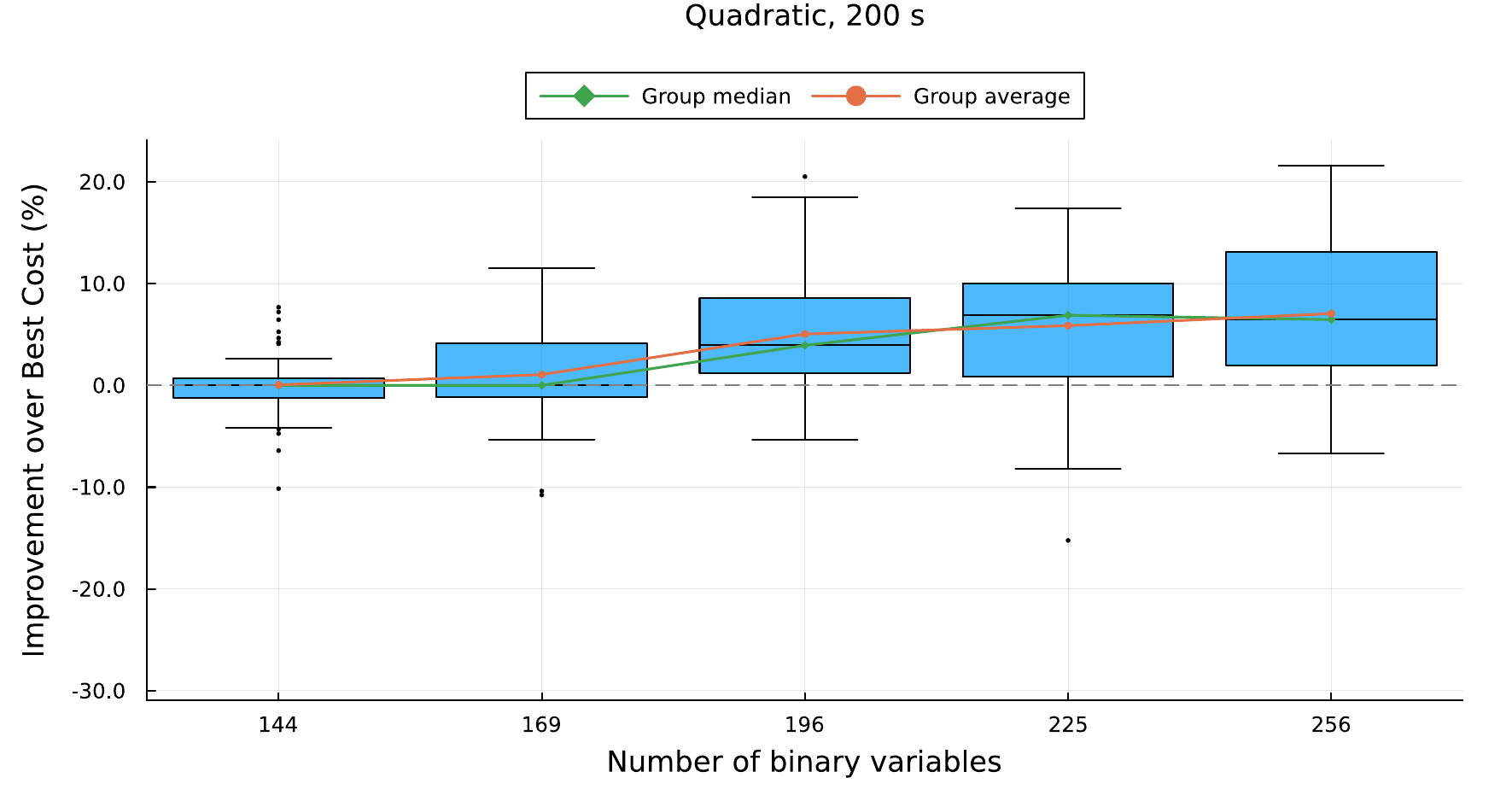}
\caption{ADD--GDF ($w_{\mathrm{filter}}=0.4$) versus Best Cost with a 200\,s budget per arm: linear
(left) and quadratic (right) assignment, 50 paired instances per size.
Boxes show quartiles and medians; whiskers extend to observations within
1.5 interquartile ranges, with other observations plotted individually.
The labelled curves show group medians and means. Positive relative
improvement favors ADD. Both panels share axis limits. The training
protocol and paired tests are given in \cref{subsec:equal_time,tab:equal_time}.}
\label{fig:add_equal_time}
\end{figure*}

A separate held-out comparison also supports the selected GDF rule at a
fixed 200\,s budget (\cref{fig:add_equal_time,subsec:equal_time}).
These are within-benchmark comparisons with weights selected on the same
instances. They support GDF with $w_{\mathrm{filter}}=0.4$ as a practical choice for this
setting.

\begin{figure*}[!tbp]
\centering
\includegraphics[width=0.64\textwidth]{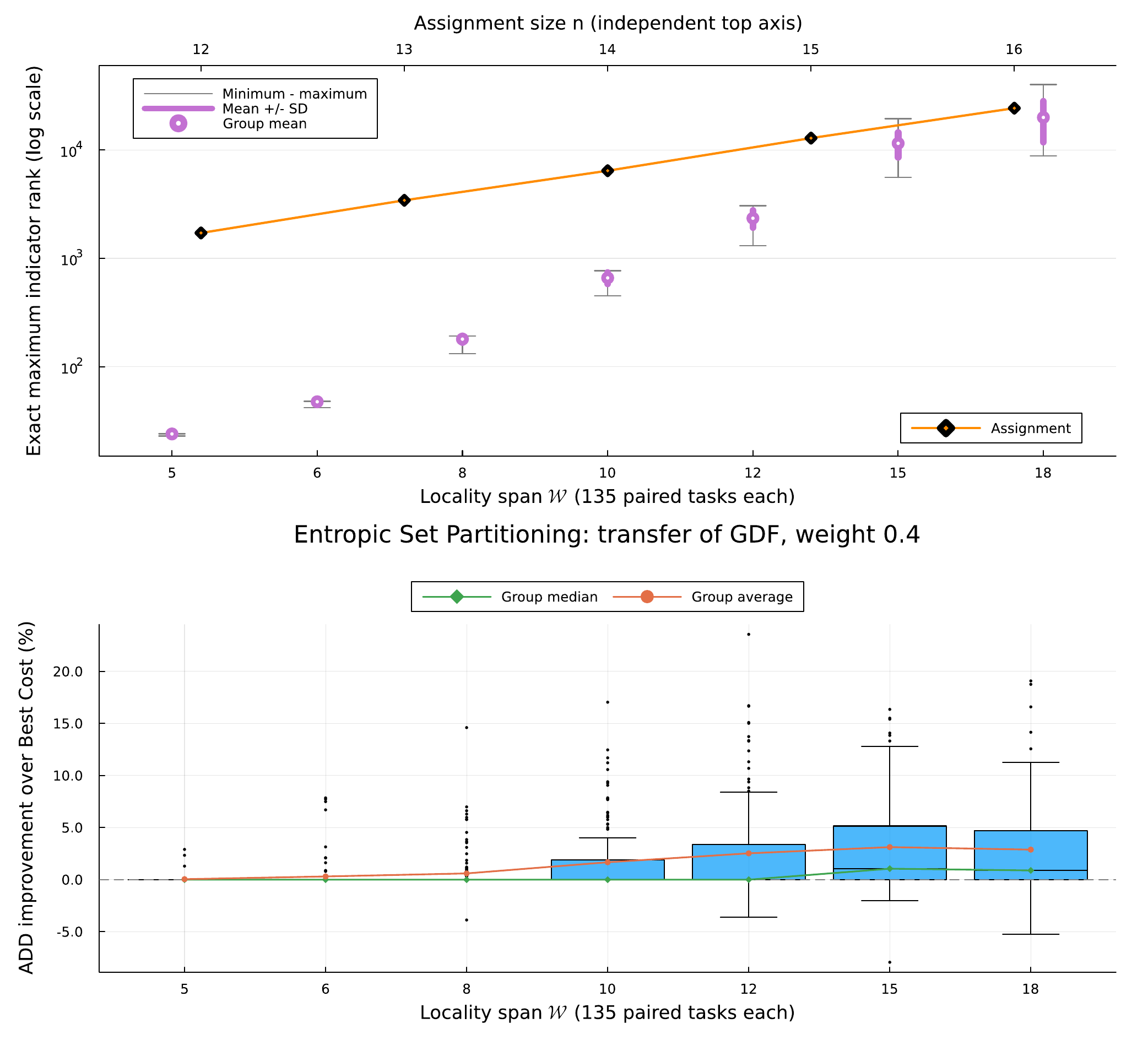}
\caption{Transfer of assignment-selected GDF ($w_{\mathrm{filter}}=0.4$)
to entropic Set Partitioning. Each locality group contains 135 paired tasks,
balanced across nine $(m,n)$ combinations. Top: maximum internal rank of
the full structural indicator, on a logarithmic axis. Thin gray intervals
show minima--maxima; thick purple intervals show the arithmetic mean
$\pm$ sample SD; hollow points mark means. These describe task dispersion,
not uncertainty in the mean. The orange line uses the independent top
axis and gives exact row-major assignment ranks for $n=12,\ldots,16$;
horizontal positions do not identify matched tasks across families. Bottom: relative final-cost improvement
$100(C_{\mathrm{BestCost}}-C_{\mathrm{ADD}})/|C_{\mathrm{BestCost}}|$.
Boxes span the first and third quartiles, with median lines; whiskers reach
the most extreme observations within 1.5 interquartile ranges and other
observations are points. Green diamonds and orange circles denote group
medians and means. Positive values favor ADD. Training uses data-built TTs,
not the full indicators whose ranks are shown above.}
\label{fig:add_partition_transfer}
\end{figure*}
\subsection{ADD transfer to nonpolynomial Set Partitioning}
\label{subsec:partition_benchmark}
To study a feasible space distinct from assignment, consider a collection
of tiles $S_j\subseteq\{1,\ldots,m\}$, one per binary variable $x_j$.
Set Partitioning selects disjoint tiles covering every element exactly once:
$A_{uj}=\mathbf1\{u\in S_j\}$ and $Ax=\mathbf1$.
Singleton tiles guarantee feasibility. We use $N=n^2$ variables; here $n$
is a size label, not the number of equalities or binary variables.
An exact counting sampler generates uniform feasible covers, without solving
the objective or constructing the TT; its recurrence is given in
Supplementary Section~S2.

For transfer of the assignment-selected ADD rule, the fixed locality grid
has $m\in\{18,19,20\}$, $n\in\{13,14,15\}$ and
$\mathcal{W}\in\{5,6,8,10,12,15,18\}$. Tiles have one to four elements and span at most
$\mathcal{W}$ consecutive element positions. Each $(\mathcal{W},m,n)$ cell contains 15 instances:
945 in total, or 135 for each $\mathcal{W}$. A fixed spatial variable order processes
tiles by their smallest element, with the corresponding singleton last.
It allows completed element constraints to leave the active frontier.
Small $\mathcal{W}$ therefore limits the number of simultaneous charge states; increasing
$\mathcal{W}$ broadens the exact-indicator rank profiles. Mean maximum ranks range from
24.0 at $\mathcal{W}=5$ to 19,962.3 at $\mathcal{W}=18$ over the balanced parameter grid
(Supplementary Table~S3). These are full-indicator ranks, not the ranks
of a TT learned from 500 points. Locality also changes the feasible covers,
so rank is a descriptive covariate rather than an isolated causal treatment.

The nonpolynomial objective is a linear cost plus entropic scenario risk:
\begin{equation}
\begin{aligned}
C_{\mathrm{ent}}(x)=\frac{c^\top x}{s}
 +\tau\log\left[\frac1{32}\sum_{a=1}^{32}
 \exp\left(\frac{W_a^\top x}{\tau s}\right)\right],\\
 s=\sqrt m.
\end{aligned}
\label{eq:partition_entropic}
\end{equation}
with $\tau=0.5$, independent $c_j\sim U[-1,1]$ and
$W_{aj}\sim\mathcal N(0,1)$. The second term smoothly emphasizes adverse
scenario costs; the factor $1/32$ makes it an empirical exponential-utility
risk rather than a sum depending on scenario count. Log-sum-exp is evaluated
stably. The transfer protocol fixes GDF at $w_{\mathrm{filter}}=0.4$, selected
on assignment, and compares it with Best Cost using matched initial and
external feasible draws, without retuning on this grid
(Supplementary Section~S2). Each pair uses 500 initial feasible draws,
15 iterations, and 5K TT plus 5K fresh external draws per iteration.
Both methods rebuild the TT after selection and train with one full sweep,
learning rate 0.05 and charge degeneracy one, using the table sampler.
The budgets match iterations and requested samples.

\textbf{Transfer and rank dependence.}
Across all 945 pairs, ADD wins on 305 tasks, ties on 628 and loses on 12.
Mean relative improvement is 1.6\% (stratified bootstrap 95\% interval
1.4--1.8\%); the median is 0.0\%. Thus the assignment-selected rule
transfers to a nonpolynomial objective and different constraint matrices,
but its benefit is concentrated in a subset of tasks.
In \cref{fig:add_partition_transfer}, the mean gain rises from less than
0.1\% at $\mathcal{W}=5$ to 3.1\% at $\mathcal{W}=15$, then levels off at 2.9\% for $\mathcal{W}=18$.
The paired one-sided Wilcoxon tests support an advantage for every group
from $\mathcal{W}=6$ onward after Holm correction across seven groups, but not for
$\mathcal{W}=5$ (Supplementary Table~S4).

\begin{figure}[!tbp]
\centering
\includegraphics[width=\columnwidth]{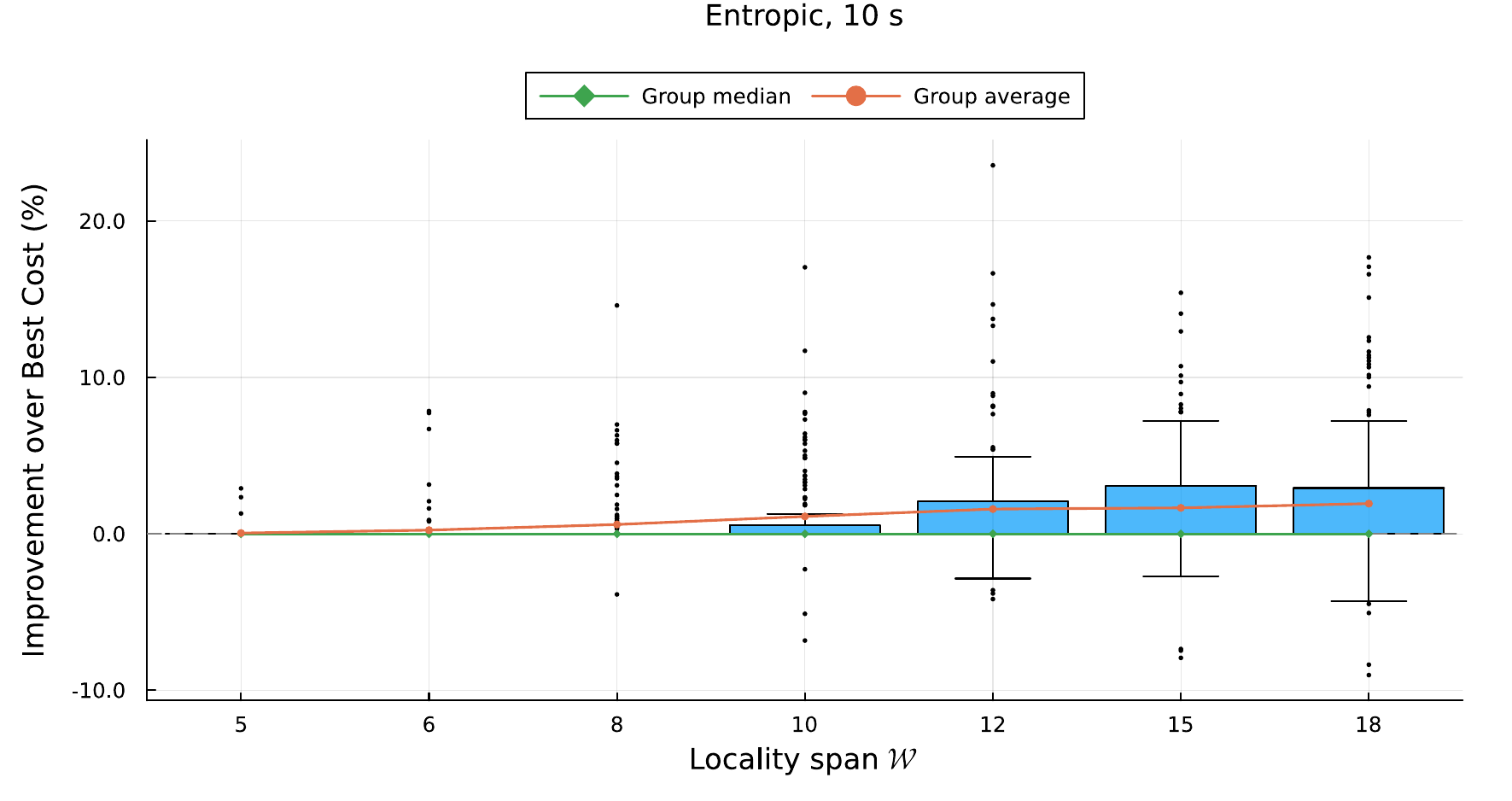}
\caption{ADD--GDF versus Best Cost on entropic Set Partitioning at 10\,s
per arm, with 135 pairs per locality group. Boxplot conventions and the
relative-improvement definition match \cref{fig:add_partition_transfer}.
Curves show group medians and means; positive values favor ADD.}
\label{fig:partition_equal_time}
\end{figure}

This pattern is consistent with the limitation of cost-only selection:
it can repeatedly favor already represented charge paths and neglect
structurally novel candidates. A larger charge graph is harder to cover
with a fixed training quota, plausibly aggravating this limitation.
Best Cost consequently falls further behind ADD on average in the
higher-rank families, where the diversity filter can help compensate for
incomplete exploration.
The positive rank--gain association persists within each fixed $(m,n)$
group. We therefore infer greater usefulness
of ADD in the higher-rank families studied. The improvement also persists
at a fixed 10\,s budget (\cref{fig:partition_equal_time,subsec:equal_time}).

\subsection{Improvement under fixed time budgets}
\label{subsec:equal_time}
The preceding comparisons fix the sampling and iteration settings. We
add equal-time comparisons to account for the cost of diversity filtering
and TT training. Each arm runs on one pinned CPU with one Julia and BLAS
thread. Paired arms share their initial feasible sample and execute
sequentially, with their order alternated across instances. Checkpoints at
10, 50 and 200\,s belong to the same continuous run; only improvements
completed before the corresponding deadline are retained. Loading,
warm-up and common input preparation are excluded, whereas initial
objective evaluation, TT construction and optimization are timed.
There is no fixed iteration cap.

For ADD assignment, the initial and retained training quotas are 1K,
with 10K TT draws and a candidate pool of 20K completed from the external
feasible sampler; 50 elites are retained. A new TT is built each iteration.
The selected rule is GDF with $w_{\mathrm{filter}}=0.4$; the comparison uses
250 held-out instances per objective, balanced over $n=12,\ldots,16$.
For EELS, the initial quota is 400, each iteration draws 10K samples and
trains on at most 400 best unique current samples. The TT is initialized
once, with or without one-time EELS ($\gamma=1$), and its learned amplitudes
persist; there are no external samples or restarts. Each objective again
has 250 pairs. Initial random samples differ from those of the earlier
fixed-iteration experiment, so the two experiments do not isolate time
budget as their only difference.

The Set Partitioning comparison reuses the 945-instance locality grid
and its fixed initial samples: quota 500, 5K TT draws plus 5K external
feasible samples, a 5\% elite fraction and TT reconstruction each iteration.
All three comparisons use one training sweep per iteration with learning
rate 0.05. ADD uses the table sampling backend; EELS uses legacy.
Objective caching is enabled for both arms of EELS and Set Partitioning,
and disabled for both assignment ADD arms. Analysis-only graph diagnostics
are excluded from the timed Set Partitioning runs; selection is unchanged.

\begin{table*}[!tbp]
\centering
\small
\begin{tabular}{llrrrrrr}
\hline
Comparison & Objective & Seconds & W/T/L & Mean (\%) & Median (\%) & Raw $p$ & Holm $p$\\
\hline
ADD: assignment & Linear & 200 & 179/58/13 & 3.1 & 2.4 & $8.9\times10^{-31}$ & $1.2\times10^{-29}$ \\
ADD: assignment & Quadratic & 200 & 158/41/51 & 3.8 & 2.7 & $1.6\times10^{-18}$ & $1.7\times10^{-17}$ \\
EELS: assignment & Linear & 200 & 120/122/8 & 2.7 & 0.0 & $7.5\times10^{-18}$ & $7.5\times10^{-17}$ \\
EELS: assignment & Quadratic & 200 & 113/112/25 & 3.6 & 0.0 & $9.6\times10^{-16}$ & $7.7\times10^{-15}$ \\
ADD: Set Partitioning & Entropic & 10 & 244/656/45 & 1.0 & 0.0 & $6.0\times10^{-27}$ & $7.7\times10^{-26}$ \\
\hline
\end{tabular}
\caption{Paired equal-time comparisons with Best Cost. W/T/L denotes
wins, exact ties and losses. One-sided paired Wilcoxon signed-rank tests
use objective values; descriptive gains are relative percentages.
Holm correction covers all 15 pooled tests (five comparisons at three
budgets), including those not selected for this table. Assignment rows
contain 250 pairs each; Set Partitioning contains 945.}
\label{tab:equal_time}
\end{table*}

\textbf{Quality and budget.}
At 200\,s, ADD improves mean assignment costs by 3.1\% and 3.8\% for
linear and quadratic objectives; EELS improves them by 2.7\% and 3.6\%
(\cref{fig:add_equal_time,fig:eels_equal_time,tab:equal_time}). At 10\,s,
ADD improves mean entropic Set Partitioning costs by 1.0\%.
Zero medians coexist with positive mean gains because many pairs tie.

\begin{figure*}[!tbp]
\centering
\includegraphics[width=\textwidth]{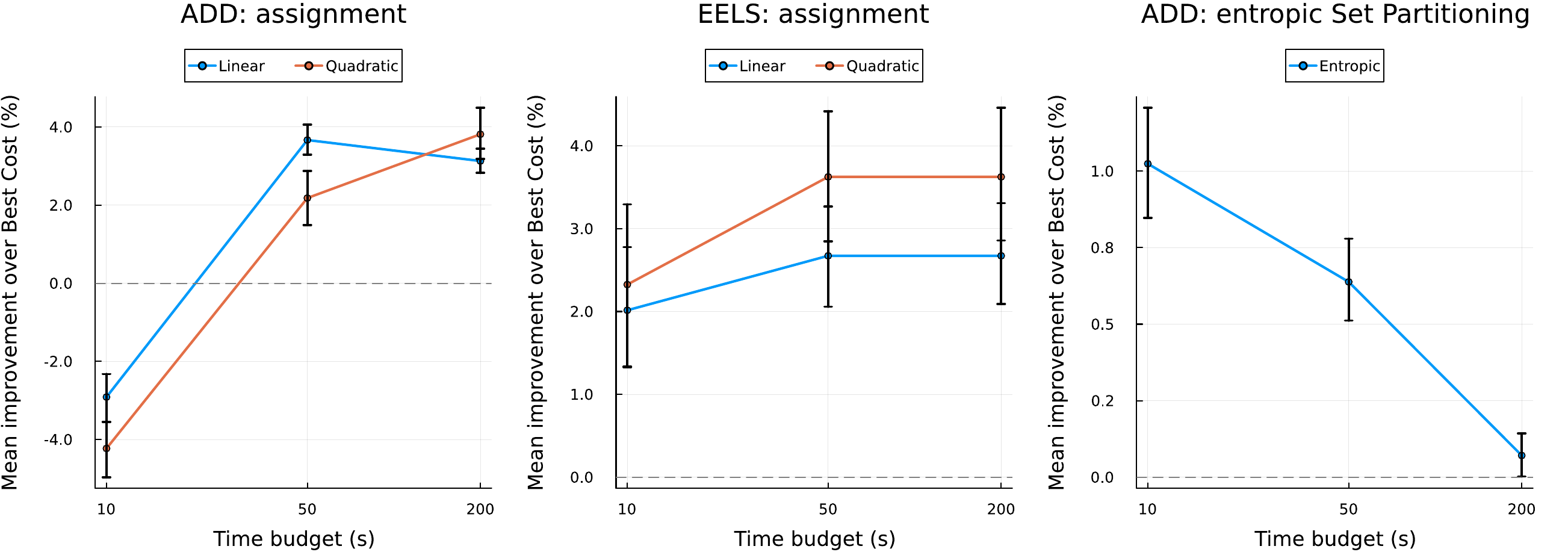}
\caption{Mean relative improvement over Best Cost versus elapsed-time
budget: ADD assignment, EELS assignment, and ADD entropic Set Partitioning.
Bars are asymmetric pointwise 95\% paired-bootstrap intervals, stratified
by problem size and, for Set Partitioning, by $(m,n,\mathcal{W})$.
Checkpoints share runs and are not independent repetitions.}
\label{fig:budget_overview}
\end{figure*}

\Cref{fig:budget_overview} reveals contrasting ADD dynamics. On assignment,
the filtering overhead is not repaid at 10\,s, but an advantage emerges
at 50--200\,s. On entropic Set Partitioning, the advantage diminishes as
Best Cost catches up. EELS reaches the same recorded incumbents at 50
and 200\,s in this experiment. 

\begin{figure}[!tbp]
\centering
\includegraphics[width=\columnwidth]{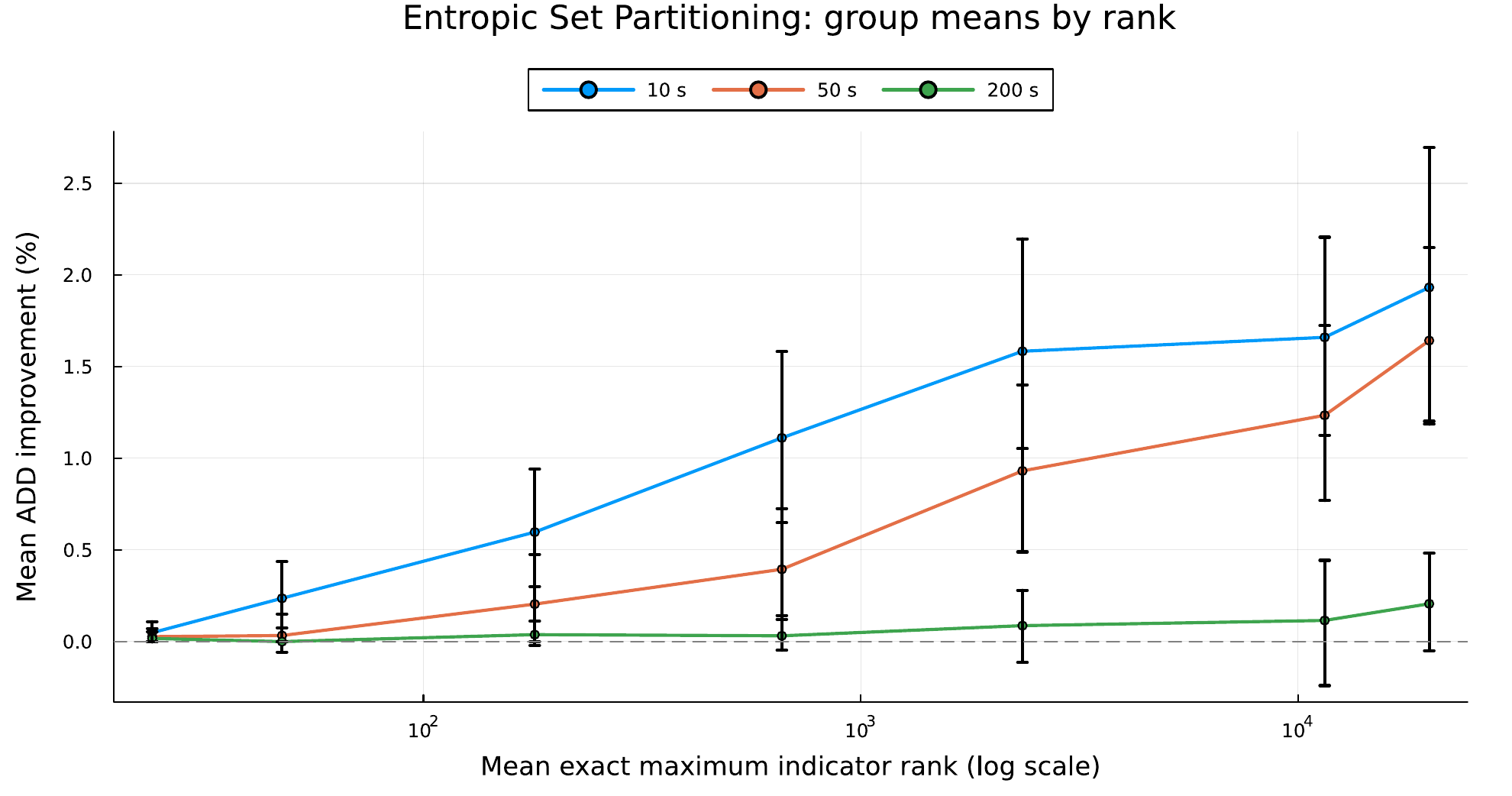}
\caption{Entropic Set Partitioning: mean ADD gain against the mean
maximum rank of the full indicator, one point per locality group and
135 pairs per point. Colors denote elapsed-time budgets; bars are
pointwise stratified paired-bootstrap 95\% intervals. Lines connect
group summaries and do not represent a fitted causal relationship.}
\label{fig:budget_rank_means}
\end{figure}

\textbf{Rank--budget trade-off.}
At 10 and 50\,s the larger-rank Set Partitioning groups exhibit larger
mean gains (\cref{fig:budget_rank_means}). Even at 200\,s a weak upward
tendency remains, reaching about 0.2\% at the largest locality, but it
is not strictly monotone: the pooled gain is only 0.1\%.
A plausible interpretation is that cost-only selection covers the
relevant paths more readily in compact feasible spaces when external
samples and sufficient time are available. Diversity filtering can then
offer little extra benefit. Harder-to-cover spaces may benefit more from
ADD once the budget absorbs filtering and training overhead. The rank
comparison in \cref{fig:add_partition_transfer} gives context: assignment
maximum ranks are approximately 1.7K--24K for $n=12$--16, while the Set Partitioning
group means span roughly 24--20K. The results suggest
a rank--budget--selection trade-off.

\subsection{Comparison with simulated annealing}
\label{subsec:sa_comparison}
For the completed SA experiment we use a separate low-rank Set Partitioning
family: connected tiles on a $3\times n$ strip, with $n=10,12,14$ and
$m=3n$. These instances are not the locality grid. Their fixed spatial order
and short tiles yield mean maximum exact ranks 27.7, 36.1 and 45.4,
respectively. Both quadratic $x^\top Qx$ and the entropic objective are used.
The two constructions and their fixed input datasets are distinguished in
Supplementary Section~S2 and Supplementary Section~S3.

We compare two symmetric-TT variants with two simulated annealing (SA)
baselines on assignment and the low-rank strip Set Partitioning family.
Each family and size $n\in\{10,12,14\}$ contains 40 paired instances with
the same 400 initial feasible draws. Separate runs use 10, 50 and 200-second
budgets on one CPU core per solver; both quadratic and entropic objectives
are evaluated as black boxes. We show $n=14$, 200-second trajectories in
\cref{fig:sa_time,fig:sa_calls}; the corresponding $n=10,12$ curves and
complete protocol appear in Supplementary Section~S3.

\textbf{Information and model variants.} Data-initialized EELS builds a
symmetric TT from the current selected data and expands its graph with
$\gamma=1$. Full-indicator TT starts from the exact indicator of all feasible
solutions. In this experiment EELS rebuilds its model, and full-indicator TT
resets its amplitudes, at every outer iteration; this differs from the
persistent-amplitude EELS--Best Cost comparison in
\cref{subsec:eels_optimization}. Informed SA receives explicit
feasibility-preserving moves. Penalty SA receives no such moves: it flips a
random bit and minimizes $C(x)+\alpha\|Ax-b\|_2^2$ with a fixed scale-based
choice of $\alpha$. EELS and penalty SA both use the constraints $A,b$ and
the initial data, but exploit them differently: EELS enforces feasibility
structurally, whereas penalty SA explores infeasible points. 

\textbf{Finite-time quality.} EELS outperforms the tested penalty-SA
configuration in every 200-second group in mean paired final quality
(Supplementary Table~S6). Penalty SA does improve its initial feasible
best in 630 of 1,440 runs, including improvements in all twelve 200-second
groups, but is substantially less effective within these budgets.

The EELS curves typically gain quickly and then plateau, whereas penalty SA
continues to make slower gains on the displayed horizon. This is consistent
with EELS acting as a rapid intensification step on data-defined support.
A flat incumbent curve does not prove a local optimum, and an improving
mean does not establish that individual SA chains never stagnate. A possible
hybrid would supply feasible solutions from SA to a subsequent EELS stage.
Such complementarity is motivated by the SA-seeded TN-GEO booster studied
by \citet{Alcazar2024}, but the proposed SA--EELS hybrid is not tested here.

\textbf{Rank and informed SA.} On assignment, informed SA has better final
cost on average in every group, with no individual TT wins (some ties).
Its transpositions encode
the permutation structure directly, while EELS receives no feasible-move
oracle. Full-indicator TT approaches informed SA at $n=10$, but loses more
at larger sizes. For row-major assignment, the full indicator has maximum TT ranks 461, 1,715 and 6,434 at $n=10,12,14$, respectively;
construction and training overhead are consistent with reduced search
throughput.
On low-rank Set Partitioning, the TT endpoints are much closer to informed
SA, particularly for the entropic objective, with many paired ties and
some TT wins.
The observed rank contrast explains a plausible computational mechanism.

\textbf{Objective-query efficiency.} The evaluation-axis curves tell a
different story from wall time. On Set Partitioning the TT methods commonly
reach comparable mean gains with far fewer objective calls; for example,
in the $n=14$ entropic 200-second runs, informed SA averages 58.3 million
calls versus 0.66 million for EELS and 0.54 million for full-indicator TT.
Both TT variants also have higher mean gains than informed SA at the
100,000-call slice on every Set Partitioning group of the 200-second runs.
Full-indicator TT has this advantage on assignment at $n=10,12$, while
EELS does not. These are descriptive slices of time-limited trajectories:
SA cooling depends on elapsed time, and every curve is truncated to a query
range reached by all 40 tasks. All calls, including repeats and infeasible penalty-SA queries,
are counted.

\begin{figure*}[!tbp]
\centering
\includegraphics[width=0.96\textwidth]{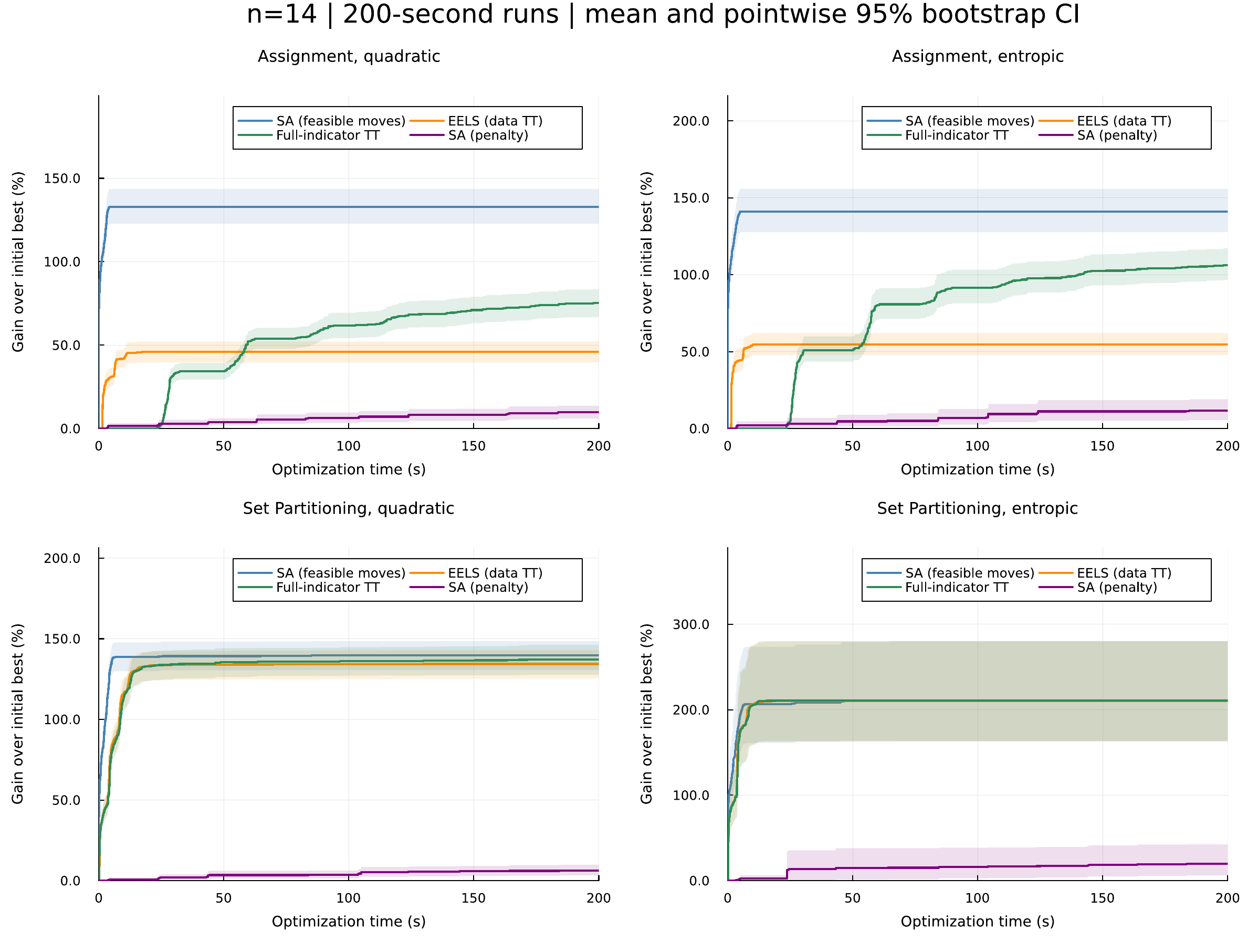}
\caption{Mean feasible-objective gain versus optimization time at $n=14$
($N=196$), with 40 paired tasks per panel and a 200-second nominal budget.
Top: assignment; bottom: spatial Set Partitioning. Left: quadratic cost;
right: entropic risk. Gain is $100(C_0-C_{\mathrm{best}})/|C_0|$ per task,
then averaged. Shading gives pointwise asymmetric 95\% bootstrap intervals
for the mean, not task standard deviations. The same initial data are used
by all four methods; penalized or infeasible costs never enter the ordinate.
EELS and full-indicator TT use the reset/rebuild protocol described in
Supplementary Section~S3, not persistent-amplitude training.}
\label{fig:sa_time}
\end{figure*}

\begin{figure*}[!tbp]
\centering
\includegraphics[width=0.96\textwidth]{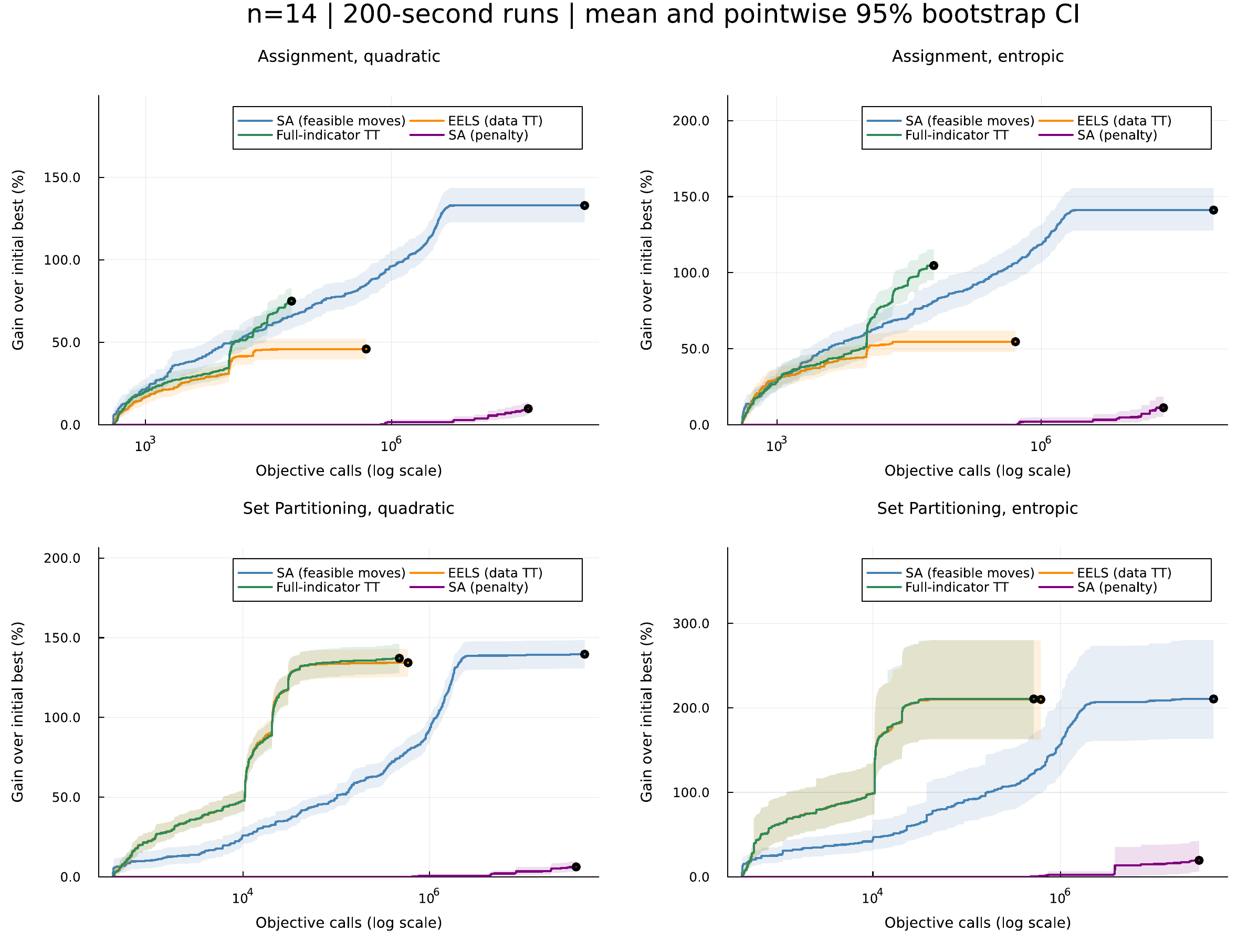}
\caption{The same $n=14$, 200-second runs as \cref{fig:sa_time}, plotted
against cumulative objective calls (logarithmic axis). Calls include the
initial evaluations, duplicates, and infeasible penalty-SA points.
Each method stops at the smallest final call count among its 40 tasks
(endpoint marker); no stopped trajectory is extrapolated. Comparisons
are meaningful only on overlapping query ranges. Confidence intervals and
gain normalization match \cref{fig:sa_time}. These curves describe
query efficiency under time-based cooling, rather than a fixed-query
benchmark.}
\label{fig:sa_calls}
\end{figure*}


\begin{figure}[!tbp]
    \centering
    \includegraphics[width=\columnwidth]{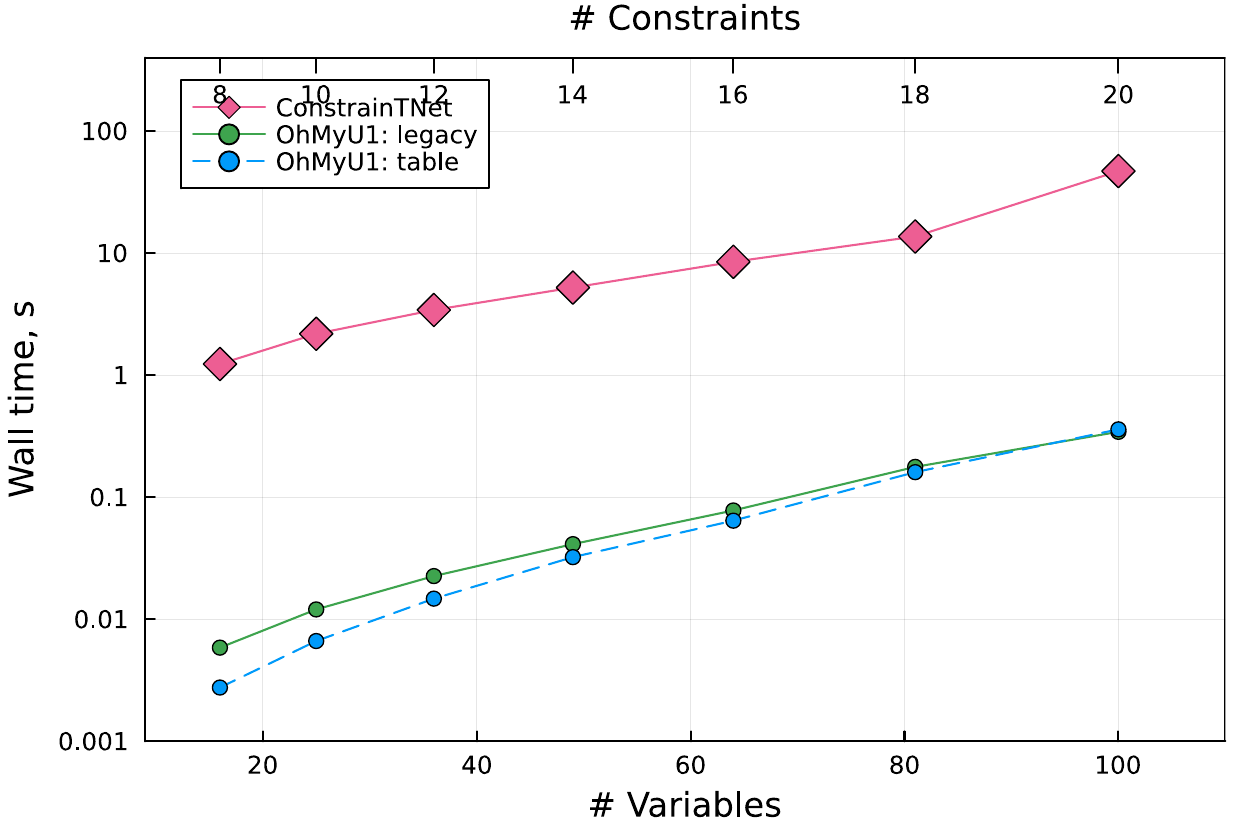}
    \caption{Median wall time for one right--left training sweep and 400
    sequential samples on the assignment problem. Both methods receive the
    same training matrix and start from the same exact uniform feasible state.
    Initial-state construction is excluded. Both OhMyU1 samplers are shown. The degeneracy of the block $\dim(V^{(c_k)}) = 1$ for each link index.}
    \label{fig:lib_comparison_fair}
\end{figure}
\begin{figure*}[!tbp]
\centering
\begin{subfigure}[t]{0.36\textwidth}
\centering\includegraphics[width=\linewidth]{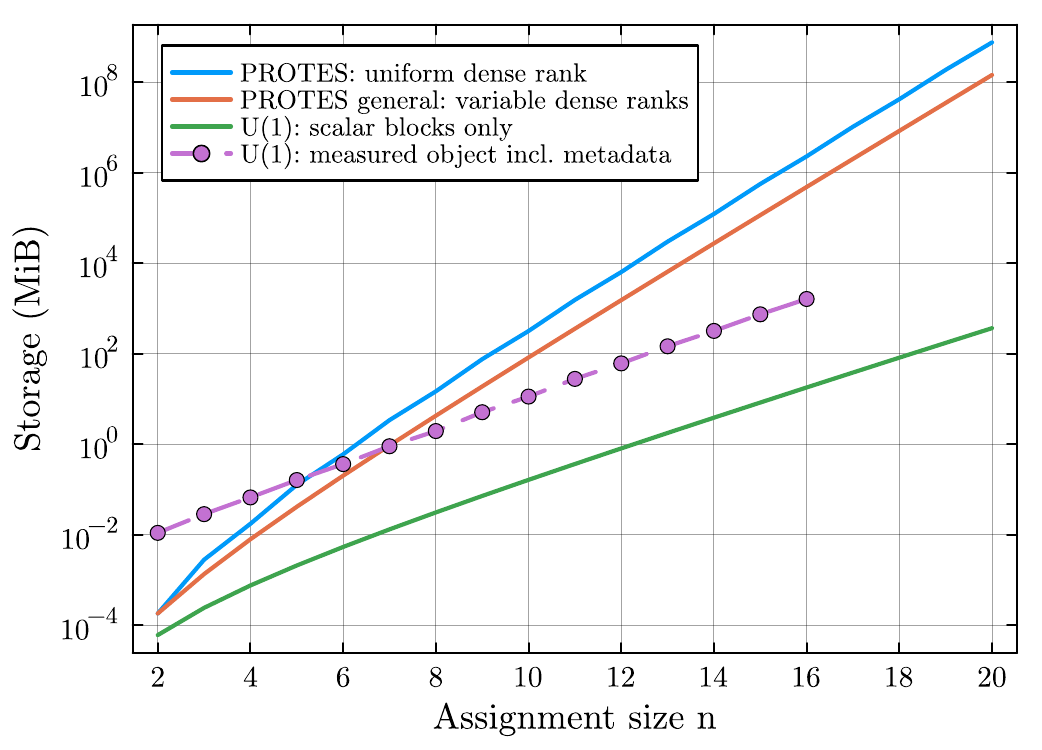}
\caption{Storage of the full assignment indicator.}
\end{subfigure}\hspace{0.02\textwidth}
\begin{subfigure}[t]{0.36\textwidth}
\centering\includegraphics[width=\linewidth]{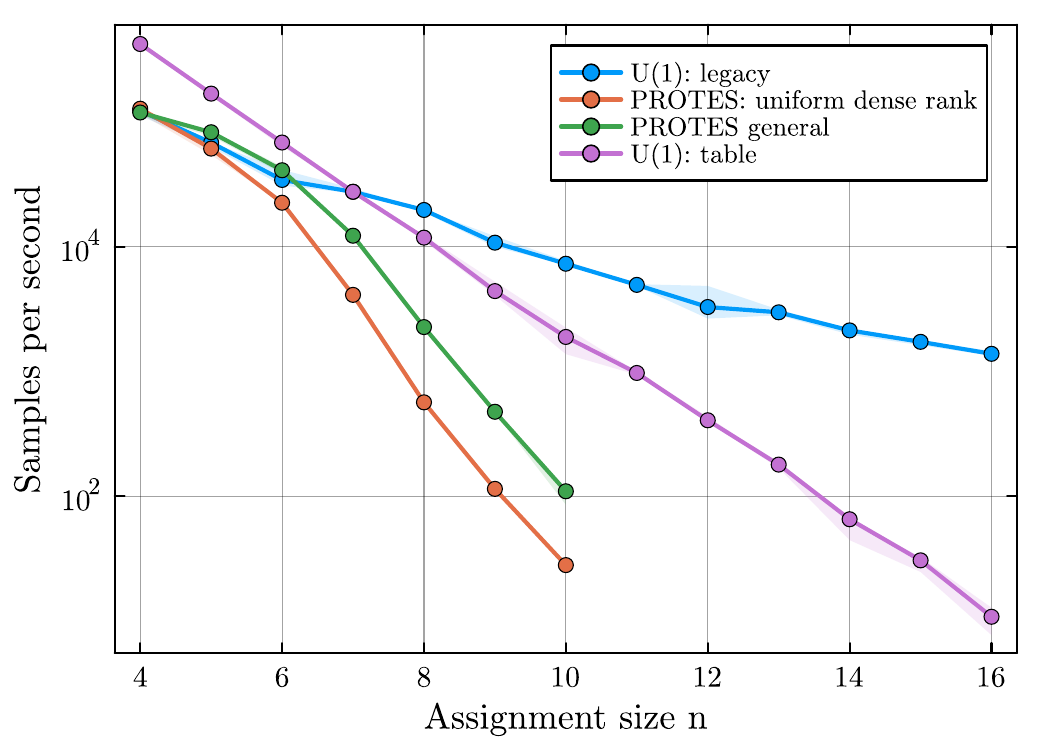}
\caption{Sampling from the full uniform feasible distribution.}
\end{subfigure}
\caption{Exact-indicator benchmarks, with no training or finite-data initialization.
(a) The two PROTES curves count float64 coefficients in dense TT cores with
uniform or minimal variable ranks; the scalar-block $U(1)$ curve counts only
allowed coefficients. Markers measure the complete Julia object, including
charges and dictionaries, through $n=16$. These are storage sizes, not peak
training memory. (b) Median sampling throughput over three seeds; shading
shows the observed minimum--maximum range. Each trial accumulates at least
one second of timed generation in batches of 100. Counts include repeated
draws. Construction, canonicalization, compilation, and diagnostic bookkeeping
are excluded. Complete dense models are actually constructed only for
$4\le n\le10$; complete $U(1)$ models are constructed through $n=16$.
The distributions and binary variable order are identical. Both OhMyU1 sampling backends are shown.}
\label{fig:protes_exact}
\par\vspace{6pt}
\centering\includegraphics[width=0.76\textwidth]{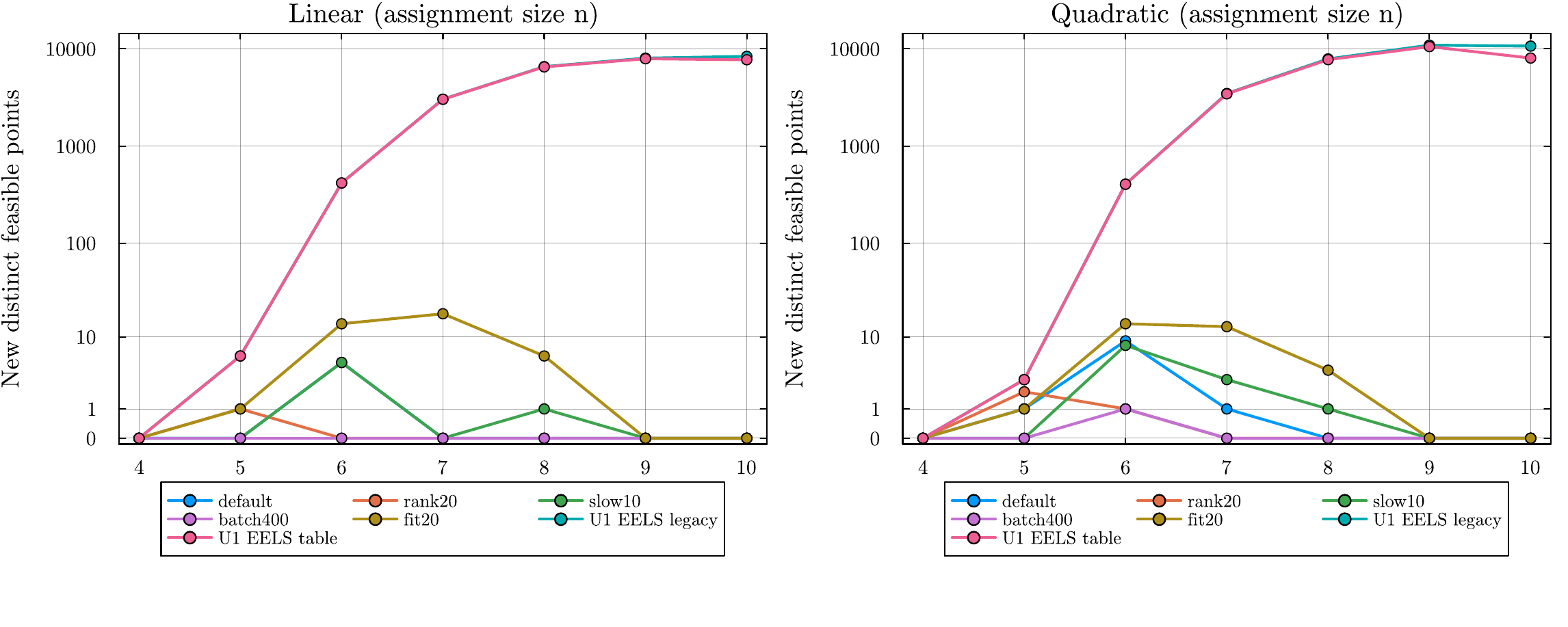}
\caption{New distinct feasible assignments within a 5-second budget, for
linear (left) and quadratic (right) objectives, versus assignment size $n$.
Each point is one matched diagnostic instance/seed; all methods receive the
same 400 initial observations. PROTES settings are in
Supplementary Table~S8; $U(1)$-EELS uses the fixed settings above, with both sampling backends shown.
Initialization, fitting, sampling, and scoring count toward the budget.
At $n=4$ the initial data already cover all 24 solutions, so zero novelty is
unavoidable. Both panels share the ordinate transformation $\log_{10}(1+z)$,
with ticks labelled by the original count $z$, allowing zeros to be shown.
No uncertainty estimate can be obtained from this single-seed screen.}
\label{fig:protes_data_short}
\end{figure*}
\subsection{Benchmarking library performance}
\label{subsection: bench_perf_res}


The implementation of \cite{LopezPiqueres2023} is not open-source. However,
the more recent article by the same authors has an accompanying GitHub
repository, \texttt{ConstrainTNet} \cite{ConstrainTNet2024}. It implements an
exact construction of a $U(1)$-symmetric tensor train: no data are needed for
initialization, although the construction itself may become intractable for
difficult constraint systems. ConstrainTNet also depends on \texttt{ITensors}
\cite{ITensor}, whereas OhMyU1 uses a specialized Julia back-end for
non-degenerate charge sectors.

We compared one right--left training sweep and 400 samples on assignment
instances with $4\leq n\leq10$, using identical full-indicator initial states
and training data. Both OhMyU1 sampling backends were measured on one CPU;
construction and compilation were excluded. The two libraries use two-site
and one-site updates, respectively. The complete input, backend and timing
protocol is provided in Supplementary Section~S4.

The results are shown in \cref{fig:lib_comparison_fair}, and detailed comparison is in Supplementary Section~S4. OhMyU1 is between
$77\times$ and $211\times$ faster with legacy, and between $85\times$ and $448\times$ with table, on the measured grid. The speedup is not
a universal comparison of the feature sets of the packages: ConstrainTNet
supports inequalities and degenerate charge subspaces, while the present
OhMyU1 kernels are specialized for non-degenerate equality-constrained
problems. The result instead quantifies the advantage of the specialized
back-end for this non-degenerate assignment workload. The two libraries use
two-site and one-site updates, respectively, so equal sweep counts do not
imply equal learning progress. The benchmark does not isolate the effect
of individual implementation choices. All benchmark code and raw measurements
are provided with OhMyU1 \cite{OhMyU1_2026}.

\paragraph{Two PROTES access models.}
We additionally compare with PROTES \cite{batsheva2023protes}, separating
access to an exact feasible-set indicator from initialization using only
feasible observations. These experiments concern binary $n\times n$
assignment, with $N=n^2$ variables in row-major order. They are distinct from
the preceding trained-kernel comparison with ConstrainTNet.
All runs use float64 and one logical CPU on the same local machine;
Supplementary Section~S5 specifies the construction, timing, training,
software versions, and recorded sample and iteration counts.

\paragraph{Exact full indicator: storage and sampling.}
When all feasible charge states are known, both representations describe
the complete set of $n!$ assignments. Their minimal total bond dimensions
coincide; the advantage of the symmetric representation is its block sparsity,
not a smaller total exact TT rank. In \cref{fig:protes_exact}, dense storage
grows substantially faster than storage of the permitted scalar blocks.
Actual Julia object sizes include substantial metadata: at $n=16$, the full
$U(1)$ object occupies $1.59$ GiB, compared with $18.0$ MiB of numerical blocks.
The corresponding dense variable-rank cores would require approximately
$472$ GiB; this last number is a storage calculation, not an allocation.

For sampling we construct the complete indicator in both libraries and
use the same uniform distribution over assignments, with no training.
All observed samples are feasible. At $n=8$ and $n=10$, OhMyU1 with the legacy sampler generates
approximately $2.0\cdot10^{4}$ and $7.3\cdot10^{3}$ samples per second,
respectively, versus $2.3\cdot10^{3}$ and $1.1\cdot10^{2}$ for the more
memory-efficient \texttt{protes\_general}. The resulting speedups are
approximately $8.7\times$ and $67\times$.
Small instances do not exhibit universal superiority.
The table sampler yields approximately 12K and 1.9K samples per second at these sizes; its per-call preparation is included. Both variants are shown.
 These measurements compare particular dense and
block-sparse implementations, rather than establishing a lower bound on
all possible PROTES implementations.

\paragraph{Initialization from feasible observations.}
We next provide the same 400 feasible observations and objective instance
to PROTES and $U(1)$-EELS, without providing an indicator or permutation
encoding to PROTES. OhMyU1 retains access to $A,b$, an intentional difference
in structural information. EELS uses $\gamma=1$, sector degeneracy one,
one training sweep at learning rate $0.05$, and up to 10K generated points
per outer iteration. All methods run until a common time budget, including
initial fitting; the number of iterations is not fixed.
Five prespecified PROTES configurations vary rank, learning rate, batch size,
Adam steps, and initial-data passes (Supplementary Table~S8).

\Cref{fig:protes_data_short} and Supplementary Figure~S5 count distinct feasible
points absent from the initial dataset, rather than repeated draws.
The 50-second, three-seed comparison in Supplementary Figure~S5
shows that longer initial fitting improves PROTES. Thus, inability to generate
any feasible point is not a robust conclusion; the observed limitation is rare
feasible generation under this data-based initialization and update rule.
New points need not improve the objective.


\section{Acknowledgements}
We thank Margarita Veshchezerova for providing the initial impetus for this research. We are grateful to Dmitry Gribanov, Elizaveta Iarovikova, Maxim Klimenko, and Ivan Samoylenko for insightful discussions. We appreciate the anonymous reviewers for their helpful feedback.

\FloatBarrier
\section*{Declaration of generative AI and AI-assisted technologies}
During the preparation of this work, the authors used DeepSeek-V4.1-Flash to assist with code review and debugging, test development, experiment automation, performance exploration, reproducible data visualization, preparation of selected manuscript passages, and language and formatting revisions. The mathematical formulation and core methodological ideas were developed by the human authors. The authors reviewed and edited the AI-assisted outputs and take responsibility for the content of this work.

\bibliography{refs}
\bibliographystyle{icml2026}

\clearpage

\clearpage
\setcounter{section}{0}
\setcounter{subsection}{0}
\setcounter{figure}{0}
\setcounter{table}{0}
\setcounter{equation}{0}
\setcounter{algorithm}{0}
\renewcommand{\theHsection}{supp.\arabic{section}}
\renewcommand{\theHfigure}{supp.\arabic{figure}}
\renewcommand{\theHtable}{supp.\arabic{table}}
\renewcommand{\theHequation}{supp.\arabic{equation}}
\renewcommand{\theHalgorithm}{supp.\arabic{algorithm}}
\icmltitlerunning{Supplementary Material}

\makeatletter
\global\icml@noticeprintedtrue
\makeatother
\onecolumn
\renewcommand{\thesection}{S\arabic{section}}
\renewcommand{\thefigure}{S\arabic{figure}}
\renewcommand{\thetable}{S\arabic{table}}
\renewcommand{\theequation}{S\arabic{equation}}
\renewcommand{\thealgorithm}{S\arabic{algorithm}}
\setcounter{topnumber}{4}
\setcounter{bottomnumber}{3}
\setcounter{totalnumber}{6}
\renewcommand{\topfraction}{0.92}
\renewcommand{\bottomfraction}{0.85}
\renewcommand{\textfraction}{0.08}
\renewcommand{\floatpagefraction}{0.75}
\raggedbottom
\begin{center}
{\Large\bfseries Supplementary Material}\par\medskip
{\large Balancing Exploration and Exploitation in Constrained Combinatorial Optimization with U(1)-Symmetric Tensor Trains}\par\medskip
Sergei Iudin, Mikhail Podobrii, Dmitry Zheltkov, Stanislav Moiseev
\end{center}
This document provides the complete strategy screening, benchmark construction,
experimental protocols and additional comparisons accompanying the main article.
Sections, figures, tables and equations prefixed with S belong to this supplement;
references to the main article are identified explicitly.

\FloatBarrier
\section{Complete ADD strategy and weight comparison}
\label{app:selection_strategy_screening}
\Cref{tab:add_linear_all,tab:add_quadratic_all} report every configuration
under the protocol of Section~4.3 of the main article, using the same
250 paired objective instances per table. Best Cost was computed once
per instance and objective (500 baseline runs); each of the 12 ADD
configurations was run on the corresponding 250 instances (6,000 ADD runs).
All 6,500 outcomes were checked for completeness, feasibility and consistency
of the stored final objective.

The linear coefficients follow $c_i\sim U[-50,50]$.
For the quadratic family a matrix $B$ is generated with independent
entries from $U[-7,7]$ and symmetrized as $Q=(B+B^T)/2$; hence
off-diagonal entries of $Q$ are averages of two uniform draws, rather
than uniform themselves. The cost is $x^TQx$, including both off-diagonal
directions and the diagonal. The source objective files are fixed across
strategies. For objective index $o=0$ (linear) or $o=1$ (quadratic),
instance index $i$ and size $n$, the initial/external source uses a
Mersenne Twister seeded by $20261012+100000o+100n+i$; the model sampling
stream is seeded by this value plus $10^7$. The source stream is independent
of MPS sampling. Subsequent source batch sizes can differ because the
number of unique MPS samples differs; identical starting data therefore
do not imply identical full adaptive trajectories.

We apply the paired Wilcoxon signed-rank test to raw differences
$C_{\mathrm{BestCost}}-C_{\mathrm{ADD}}$, with the right-tailed alternative
favoring ADD, as implemented by \texttt{SignedRankTest} in
\texttt{HypothesisTests.jl}. Zeros are handled by that implementation;
win/tie/loss counts use strict comparisons and exact equality of stored costs.
The tables give the positive-rank sum $W$, the unadjusted one-sided $p$,
Holm-adjusted $p$ over the 12 comparisons in that objective family, and
the more conservative adjustment over all 24 baseline comparisons.
Mean and median gains use the per-instance relative improvement $R$;
Rows are ordered by mean relative improvement,
matching the weight-selection criterion.
Of the 24 baseline comparisons, 23 remain significant at $0.05$ after
the joint Holm adjustment; Mixed Selection with $w_{\mathrm{mix}}=0.1$ on quadratic
objectives is the exception ($p=0.086$).

\begin{table}[!tbp]
\centering
\small
\setlength{\tabcolsep}{3pt}
\caption{Linear: every strategy and weight. All rows contain 250 pairs against Best Cost. W/T/L denotes wins/ties/losses for ADD; win rate, mean and median gains are percentages. Holm (12) and Holm (24) adjust the one-sided baseline comparisons within this objective and across both objectives, respectively.}
\label{tab:add_linear_all}
\begin{tabular}{lrcrrrrrrr}
\toprule
Strategy & Weight & W/T/L & Win (\%) & Mean (\%) & Median (\%) & $W$ & $p$ & Holm (12) & Holm (24) \\
\midrule
GDF & $w_{\mathrm{filter}}=0.4$ & 235/14/1 & 94.0 & 6.9 & 5.9 & 27963 & $9.5\!\times\!10^{-41}$ & $1\!\times\!10^{-39}$ & $2.2\!\times\!10^{-39}$ \\
GDF & $w_{\mathrm{filter}}=0.3$ & 239/10/1 & 95.6 & 6.5 & 6.0 & 28908 & $2.4\!\times\!10^{-41}$ & $2.8\!\times\!10^{-40}$ & $5.6\!\times\!10^{-40}$ \\
GDF & $w_{\mathrm{filter}}=0.2$ & 232/18/0 & 92.8 & 5.8 & 5.1 & 27028 & $4.1\!\times\!10^{-40}$ & $4.1\!\times\!10^{-39}$ & $9.1\!\times\!10^{-39}$ \\
Weighted Rank & $w_{\mathrm{rank}}=0.35$ & 222/21/7 & 88.8 & 5.6 & 5.0 & 26207 & $6.8\!\times\!10^{-39}$ & $5.5\!\times\!10^{-38}$ & $1.4\!\times\!10^{-37}$ \\
Weighted Rank & $w_{\mathrm{rank}}=0.5$ & 227/15/8 & 90.8 & 5.5 & 5.0 & 27531 & $1.7\!\times\!10^{-39}$ & $1.5\!\times\!10^{-38}$ & $3.5\!\times\!10^{-38}$ \\
Weighted Rank & $w_{\mathrm{rank}}=0.7$ & 225/18/7 & 90.0 & 4.9 & 4.2 & 26746 & $1.5\!\times\!10^{-38}$ & $1.1\!\times\!10^{-37}$ & $2.9\!\times\!10^{-37}$ \\
GDF & $w_{\mathrm{filter}}=0.1$ & 214/33/3 & 85.6 & 4.3 & 3.7 & 23630 & $1.6\!\times\!10^{-37}$ & $9.8\!\times\!10^{-37}$ & $2.9\!\times\!10^{-36}$ \\
Mixed & $w_{\mathrm{mix}}=0.4$ & 189/38/23 & 75.6 & 2.8 & 2.4 & 20727 & $2.4\!\times\!10^{-26}$ & $1.2\!\times\!10^{-25}$ & $3.4\!\times\!10^{-25}$ \\
Mixed & $w_{\mathrm{mix}}=0.3$ & 178/36/36 & 71.2 & 2.1 & 1.7 & 19522 & $4.7\!\times\!10^{-19}$ & $1.4\!\times\!10^{-18}$ & $4.7\!\times\!10^{-18}$ \\
Weighted Rank & $w_{\mathrm{rank}}=0.1$ & 157/74/19 & 62.8 & 2.0 & 1.2 & 14457 & $3.4\!\times\!10^{-23}$ & $1.4\!\times\!10^{-22}$ & $4\!\times\!10^{-22}$ \\
Mixed & $w_{\mathrm{mix}}=0.2$ & 139/62/49 & 55.6 & 1.3 & 0.6 & 13656 & $8.4\!\times\!10^{-11}$ & $1.7\!\times\!10^{-10}$ & $4.2\!\times\!10^{-10}$ \\
Mixed & $w_{\mathrm{mix}}=0.1$ & 94/109/47 & 37.6 & 0.6 & 0.0 & 6813 & $1\!\times\!10^{-4}$ & $1\!\times\!10^{-4}$ & $3\!\times\!10^{-4}$ \\
\bottomrule
\end{tabular}
\vspace{1.5em}
\centering
\small
\setlength{\tabcolsep}{3pt}
\caption{Quadratic: every strategy and weight. All rows contain 250 pairs against Best Cost. W/T/L denotes wins/ties/losses for ADD; win rate, mean and median gains are percentages. Holm (12) and Holm (24) adjust the one-sided baseline comparisons within this objective and across both objectives, respectively.}
\label{tab:add_quadratic_all}
\begin{tabular}{lrcrrrrrrr}
\toprule
Strategy & Weight & W/T/L & Win (\%) & Mean (\%) & Median (\%) & $W$ & $p$ & Holm (12) & Holm (24) \\
\midrule
GDF & $w_{\mathrm{filter}}=0.4$ & 209/26/15 & 83.6 & 6.6 & 6.2 & 24425 & $2\!\times\!10^{-34}$ & $2.5\!\times\!10^{-33}$ & $3.5\!\times\!10^{-33}$ \\
GDF & $w_{\mathrm{filter}}=0.3$ & 197/33/20 & 78.8 & 5.9 & 5.4 & 22423 & $1.3\!\times\!10^{-30}$ & $1.4\!\times\!10^{-29}$ & $2\!\times\!10^{-29}$ \\
Weighted Rank & $w_{\mathrm{rank}}=0.35$ & 189/28/33 & 75.6 & 5.9 & 5.4 & 22988 & $8.3\!\times\!10^{-29}$ & $8.3\!\times\!10^{-28}$ & $1.2\!\times\!10^{-27}$ \\
Weighted Rank & $w_{\mathrm{rank}}=0.5$ & 184/23/43 & 73.6 & 5.2 & 4.2 & 22510 & $2.2\!\times\!10^{-22}$ & $1.8\!\times\!10^{-21}$ & $2.4\!\times\!10^{-21}$ \\
GDF & $w_{\mathrm{filter}}=0.2$ & 179/44/27 & 71.6 & 4.6 & 3.5 & 19655 & $4.3\!\times\!10^{-26}$ & $3.9\!\times\!10^{-25}$ & $5.7\!\times\!10^{-25}$ \\
Weighted Rank & $w_{\mathrm{rank}}=0.7$ & 151/28/71 & 60.4 & 3.5 & 2.0 & 19289 & $2.7\!\times\!10^{-13}$ & $1.6\!\times\!10^{-12}$ & $2.2\!\times\!10^{-12}$ \\
Mixed & $w_{\mathrm{mix}}=0.4$ & 165/35/50 & 66.0 & 2.7 & 1.8 & 18031 & $1\!\times\!10^{-12}$ & $5.1\!\times\!10^{-12}$ & $7.2\!\times\!10^{-12}$ \\
GDF & $w_{\mathrm{filter}}=0.1$ & 145/73/32 & 58.0 & 2.5 & 1.1 & 13019 & $2.5\!\times\!10^{-14}$ & $1.8\!\times\!10^{-13}$ & $2.2\!\times\!10^{-13}$ \\
Mixed & $w_{\mathrm{mix}}=0.3$ & 140/45/65 & 56.0 & 2.1 & 1.3 & 15044 & $6.6\!\times\!10^{-8}$ & $2\!\times\!10^{-7}$ & $2.6\!\times\!10^{-7}$ \\
Weighted Rank & $w_{\mathrm{rank}}=0.1$ & 137/71/42 & 54.8 & 1.9 & 0.5 & 12589 & $3.3\!\times\!10^{-11}$ & $1.3\!\times\!10^{-10}$ & $2\!\times\!10^{-10}$ \\
Mixed & $w_{\mathrm{mix}}=0.2$ & 122/57/71 & 48.8 & 1.4 & 0.0 & 12125 & $1.9\!\times\!10^{-4}$ & $3.8\!\times\!10^{-4}$ & $3.8\!\times\!10^{-4}$ \\
Mixed & $w_{\mathrm{mix}}=0.1$ & 98/73/79 & 39.2 & 0.5 & 0.0 & 8810 & $0.086$ & $0.086$ & $0.086$ \\
\bottomrule
\end{tabular}
\end{table}


\FloatBarrier
\section{Set Partitioning construction and ADD transfer protocol}
\label{app:partition_protocol}
\paragraph{Elements, tiles and binary variables.}
Start with the ordered set $V=\{1,\ldots,m\}$. A candidate tile $S_j$ is
a nonempty subset of these elements, not necessarily a contiguous interval.
The binary variable $x_j$ indicates whether this tile is selected. Defining
$A_{uj}=\mathbf1\{u\in S_j\}$ and $b=\mathbf1_m$, the constraint $Ax=b$
requires the selected tiles to be disjoint and to cover every element exactly
once. There are $m$ equalities and $N$ binary variables. We set $N=n^2$ only
to retain the size labels of the assignment experiments.

\paragraph{Locality and the tile catalogue.}
Each tile contains at most four elements. Its span is
$\max S-\min S+1$, where $\min S$ and $\max S$ are its smallest and largest
element labels. For example, $S=\{3,5,6\}$ has size three and span four:
it lies inside the interval $\{3,4,5,6\}$. We admit only tiles with span at
most $\mathcal W$. Increasing $\mathcal W$ allows more distant elements to
occur in the same tile, without increasing the maximum tile size.

The tile-size proportions come from the reference case $m=18$, $n=15$,
$\mathcal W=5$. For a tile of size $s\ge2$ and span exactly $\ell$, its
containing interval has $m-\ell+1$ possible positions. Both endpoints
belong to the tile, and $s-2$ of the $\ell-2$ interior elements are chosen.
Consequently the number of eligible size-$s$ tiles is
\[
K_s(m,\mathcal W)=\sum_{\ell=s}^{\mathcal W}
(m-\ell+1)\binom{\ell-2}{s-2}.
\]
At the reference setting, the catalogue contains 18 singletons, 62 pairs,
88 triples and 57 quadruples: exactly $225=15^2$ tiles.
\paragraph{Generating one instance.}
For each $(\mathcal W,m,n)$, include all $m$ singletons, guaranteeing that
at least one exact cover exists. Allocate the remaining $N-m$ columns to
sizes 2, 3 and 4 in proportions $62:88:57$. Take the floor of each quota
and assign the leftover columns to the largest fractional remainders.
Include all $m-1$ adjacent pairs $\{u,u+1\}$ within the size-two quota.
Fill the remaining quotas by uniform sampling without replacement from
the eligible, not-yet-selected tiles of each size. Thus columns are distinct;
the algorithm checks that each catalogue is large enough. When a quota
exhausts its catalogue, that part of the matrix is deterministic. In
particular, the reference matrix at $(5,18,15)$ is deterministic, although
its objective arrays and initial draws vary between instances.

Sort columns by the smallest element of their tile. Within each such group,
put nonsingletons in lexicographic order and the singleton last. For example,
all tiles starting at element 3 precede tiles starting at element 4, and
$\{3\}$ is processed last in its group. This spatial ordering is fixed
by construction and is not selected by optimization performance.
The grid uses $\mathcal W\in\{5,6,8,10,12,15,18\}$,
$m\in\{18,19,20\}$ and $n\in\{13,14,15\}$, with 15 instances per cell:
945 tasks, or 135 per locality group. Constraints are generated independently
of the objective values.

\paragraph{Exact uniform feasible sampling.}
Let $U\subseteq V$ be the elements already covered, and let $F(U)$ count
the ways to complete their cover using available tiles. Set $F(V)=1$.
Otherwise choose the smallest uncovered element $u$ and compute
\begin{equation}
F(U)=\sum_{j:\ u\in S_j,\,S_j\cap U=\varnothing}F(U\cup S_j).
\label{eq:partition_counting}
\end{equation}
Every completion must contain exactly one tile covering $u$. Its other
elements must be uncovered, which explains the two restrictions in the sum.
Memoization stores these counts as arbitrary-precision integers. To draw
a cover, choose each next tile with probability $F(U\cup S_j)/F(U)$.
Each exact cover has one path under this smallest-uncovered-element rule,
and its probability telescopes to $1/F(\varnothing)$. The resulting binary
vectors are therefore exactly uniform over feasible solutions, with
replacement. Integer rejection sampling avoids converting large counts to
floating point.
The saved inputs contain 500 initial draws, objective arrays and random seeds.
The solver verifies input hashes and $AX=b$; external batches are also
checked for feasibility and identity between the paired arms.

\paragraph{Left-to-right construction of the full indicator and its exact ranks.}
The task is to determine which virtual charges can occur in a complete
feasible solution, before assigning entries to the TT cores. We process
the binary decisions $x_1,\ldots,x_N$ from left to right in the fixed
spatial column order. After the first $k$ decisions, the residual charge is
\[
q_k=\mathbf1-\sum_{j=1}^{k}A_jx_j.
\]
An entry equal to one denotes an uncovered element; zero denotes an
element already covered. Let $Q_k$ be the set of distinct residual charges
that are reachable by a prefix $(x_1,\ldots,x_k)$ and admit at least one
feasible suffix $(x_{k+1},\ldots,x_N)$. Initially
$Q_0=\{\mathbf1\}$, and at the end $Q_N=\{\mathbf0\}$.
The central issue is how to construct these sets without enumerating
complete solutions or solving a separate completion problem for each charge.

\textit{State enumeration and the role of singletons.}
Variables are grouped by the smallest element of their tile. In the group
with minimum element $u$, all multielement tiles precede the singleton
$\{u\}$. Every singleton is present. From each $q\in Q_{k-1}$, consider
the two decisions for tile $S_k$:
\begin{itemize}
\item If $x_k=0$, the candidate residual is $q'=q$.
\item If $x_k=1$, the candidate is $q'=q-A_k$. Allow this transition only
when $q_v=1$ for every $v\in S_k$, so no element is covered twice.
\item At the singleton $S_k=\{u\}$, also require $q'_u=0$.
After this variable no later tile contains $u$, so leaving it uncovered
cannot lead to a feasible solution.
\end{itemize}
Merge equal surviving residuals to obtain $Q_k$, retaining the allowed
transitions between successive layers. Reachability is automatic because
each candidate comes from a retained state of the preceding layer.
Completion is also guaranteed: before the singleton $\{u\}$, every
uncovered element still has its singleton among the unprocessed variables;
after it, the same is true once $u$ has been covered. Skip all remaining
multielement tiles and select exactly the singletons of uncovered elements.
This gives an explicit feasible suffix for every retained state.
Conversely, every feasible solution obeys the transition rules, so none is
lost. Thus the forward pass constructs exactly $Q_k$, rather than a
superset requiring a later backward feasibility search. This argument
depends on the singleton-last order and the presence of all singletons;
it does not apply to arbitrary constraint matrices or variable orders.

\textit{From states and transitions to the indicator TT.}
The forward pass determines $Q_k$ and its incoming transitions at every
binary variable, including those before a singleton. No waiting until the
end of a group or backward pruning of earlier layers is required. To form
the corresponding core, assign one virtual index to each charge in $Q_k$
and set
\[
G^{[k]}_{q,x,q'}=
\begin{cases}
1, & q\in Q_{k-1},\ q'\in Q_k,\ q'=q-A_kx,\\
0, & \text{otherwise}.
\end{cases}
\]
The boundary indices are fixed to $q_0=\mathbf1$ and $q_N=\mathbf0$.
A fixed binary vector determines at most one state path; it has a complete
path precisely when it is an exact cover. Contracting the cores therefore
gives one for every feasible vector and zero otherwise. For the rank
measurements reported here, the implementation only enumerates states,
counts transitions and accumulates prefix-path counts; it does not allocate
TT cores. Its final path count provides an independent check of the
sampler's $F(\varnothing)$.

\textit{Why locality bounds the work.}
While processing the group with minimum element $u$, all earlier elements
have been covered, and elements beyond $u+\mathcal W-1$ have not yet been
touched. Only the coverage of at most $\mathcal W$ consecutive elements
can vary. Store it as a bit mask, with one bit for a covered element
(the complement of the residual convention above). After processing
$\{u\}$, drop its coordinate and advance this active window. At any fixed
cut, the mask and the window position determine the full residual uniquely.
Consequently,
\[
R_k:=|Q_k|\leq 2^{\mathcal W}.
\]
Each state has at most two outgoing transitions. With machine-word masks
and constant-time state lookup, the forward traversal takes
$O(\sum_{k=1}^{N}R_{k-1})\subseteq O(NR_{\max})$ operations, where
$R_{\max}=\max_{0\leq k\leq N}R_k$.
The current implementation uses two mask-indexed count arrays of length
$2^{\mathcal W_{\mathrm{actual}}}$, where $\mathcal W_{\mathrm{actual}}\leq\mathcal W$ is the actual maximum tile span, as well as
lists of occupied masks. Including initialization of these arrays, its
state-enumeration cost is bounded by
\[
O\!\left(2^{\mathcal W}+\sum_{k=1}^{N}R_{k-1}\right),
\qquad\text{with }O(2^{\mathcal W})\text{ working memory}.
\]
These bounds concern the state pass after preparing the input and masks;
the code uses fixed-width masks and path counts for the benchmark sizes.
Retaining the entire transition graph or the sparse cores additionally
requires $O(\sum_{k=1}^{N}R_{k-1})$ storage. Dense cores would instead
store $2R_{k-1}R_k$ entries per site. The sparse traversal is therefore
linear in the number of visited states, not quadratic in adjacent ranks,
but it is not a polynomial-time guarantee in $\mathcal W$: the ranks
themselves may grow exponentially.

\textit{A complete three-element example.}
Take $V=\{1,2,3\}$ and the ordered tiles
\[
S_1=\{1,2\},\quad S_2=\{1\},\quad
S_3=\{2,3\},\quad S_4=\{2\},\quad S_5=\{3\}.
\]
Here the maximum span is two. Start with $Q_0=\{(1,1,1)\}$.
At the first variable, skipping $\{1,2\}$ leaves $(1,1,1)$, while
selecting it gives $(0,0,1)$. Thus
$Q_1=\{(1,1,1),(0,0,1)\}$. Both states have a singleton completion.
Next process $S_2=\{1\}$:
\begin{itemize}
\item From $(1,1,1)$, skipping leaves element 1 uncovered after its last
available tile, so discard that transition. Selecting gives $(0,1,1)$.
\item From $(0,0,1)$, selecting would cover element 1 twice, so disallow
it. Skipping gives $(0,0,1)$ and closes the group correctly.
\end{itemize}
Hence $Q_2=\{(0,1,1),(0,0,1)\}$. The transitions at the first two
variables correspond, for example, to the core entries
$G^{[1]}_{(1,1,1),1,(0,0,1)}=1$ and
$G^{[2]}_{(0,0,1),0,(0,0,1)}=1$.

At $S_3=\{2,3\}$, state $(0,1,1)$ allows both skipping and selecting,
giving $(0,1,1)$ and $(0,0,0)$; state $(0,0,1)$ allows only skipping.
Thus $Q_3=\{(0,1,1),(0,0,0),(0,0,1)\}$.
At singleton $S_4=\{2\}$, state $(0,1,1)$ must select it, whereas
$(0,0,0)$ and $(0,0,1)$ must skip it. Merge the identical resulting
residuals to obtain $Q_4=\{(0,0,0),(0,0,1)\}$.
Finally, singleton $S_5=\{3\}$ is skipped from
$(0,0,0)$ and selected from $(0,0,1)$, yielding $Q_5=\{(0,0,0)\}$.
The layer sizes are $(1,2,2,3,2,1)$, all at most $2^2=4$.
There are exactly three complete paths, with binary vectors
$(1,0,0,0,1)$, $(0,1,1,0,0)$ and $(0,1,0,1,1)$.
In particular, at cut 2 the charge $(0,0,1)$ has one feasible suffix
and $(0,1,1)$ has two: the bond has two states, although the indicator
contains three nonzero entries. Multiplying the layer sizes would not
count solutions, because only the recorded transitions connect layers.

\textit{Why these state counts are the exact TT ranks.}
At cut $k$, unfold the full indicator into a matrix with binary prefixes
as rows and binary suffixes as columns. A reachable prefix with residual
$q$ is compatible with a suffix exactly when
$\sum_{j=k+1}^{N}A_jx_j=q$. All prefixes with the same residual therefore
have identical compatible suffix sets, and different residuals have
disjoint suffix sets. Grouping rows and columns by residual produces one
all-ones rectangle per $q\in Q_k$, with zeros between these rectangles
and possibly additional zero rows or columns. Every rectangle is nonempty:
reachability supplies a prefix and the singleton completion supplies a suffix.
Each rectangle has matrix rank one, and their row and column supports are
disjoint. The unfolding rank is consequently exactly $|Q_k|=R_k$.
This is the minimal ordinary TT bond dimension for the specified variable
order, not merely an upper bound from a particular construction. The
symmetric representation has the same total bond dimension, with charge
degeneracy one. Maximum and mean ranks refer to the $N-1$ internal cuts.
Table~\ref{tab:partition_ranks} pools 135 instances per locality group;
its SD includes variation across the nine $(m,n)$ combinations, not only
within each cell. The bound $2^{\mathcal W}$ limits the number of charge
states at a cut, not the number of feasible solution paths.

\begin{table}[!tbp]
\centering\small
\caption{Exact-indicator ranks in the fixed locality grid: mean $\pm$ sample SD over all 135 instances at each span.}
\label{tab:partition_ranks}
\begin{tabular}{rrr}\toprule
$\mathcal{W}$ & Maximum internal rank & Mean internal rank\\\midrule
5 & $24.0 \pm 0.2$ & $17.2 \pm 0.8$\\
6 & $48 \pm 1$ & $31 \pm 2$\\
8 & $180 \pm 10$ & $100 \pm 10$\\
10 & $660 \pm 80$ & $300 \pm 40$\\
12 & $(2.3 \pm 0.4)\times10^{3}$ & $900 \pm 200$\\
15 & $(12 \pm 3)\times10^{3}$ & $(3.0 \pm 0.8)\times10^{3}$\\
18 & $(20 \pm 8)\times10^{3}$ & $(5 \pm 2)\times10^{3}$\\
\bottomrule\end{tabular}
\end{table}

\paragraph{Transfer configuration.}
We transfer Best Cost + GDF with $w_{\mathrm{filter}}=0.4$, the shared
assignment-selected choice, to the entropic objective in
Equation~(13) of the main article. Both arms start from the same 500 feasible draws
and use 15 outer iterations, one training sweep at learning rate 0.05,
charge degeneracy one, 5K TT samples and 5K independent exact-uniform
external draws per iteration. Five percent of the initial distinct training
size is retained as an elite, with a minimum of one point; no worst-point
retention is used. The model is rebuilt from selected data each iteration,
without EELS. Boltzmann temperature is the standard deviation of current
training costs, with uniform weights when it is zero. The table sampling
backend is explicit. External RNG streams are matched by instance and iteration
and separate from model sampling; matching seeds do not imply identical
TT samples after the distributions diverge. Graph-history
diagnostics record the fraction of externally drawn solutions with zero
novelty.

There are 945 pairs, or 1,890 solver runs. Four Slurm allocations each execute
32 isolated one-core processes, with one BLAS thread and private per-instance
checkpoints.
The completed comparison uses equal iteration and requested-sample budgets,
not equal runtime or objective-call budgets.

\paragraph{Statistical analysis.}
We apply one-sided paired Wilcoxon signed-rank tests to raw final-cost
differences, excluding exact ties, with Holm adjustment over the seven
locality groups.
a location-shift interpretation assumes symmetric paired differences.
The overall test gives $p=2\times10^{-48}$; a paired sign test also
supports the advantage ($p=6.8\times10^{-75}$). Wins/ties/losses are
unchanged by an absolute tolerance of $10^{-10}$. The smallest absolute
baseline cost is 0.842, so the relative gains do not arise from nearly zero
denominators. Mean-gain intervals use 5K paired bootstrap resamples
stratified by the 63 design cells; these are pointwise intervals.

\begin{table}[!tbp]
\centering\small
\caption{ADD transfer: 135 pairs per locality. Relative gains and pointwise bootstrap intervals are percentages. W/T/L denotes wins, exact ties and losses; Holm adjustment covers seven one-sided Wilcoxon tests.}
\label{tab:partition_transfer}
\begin{tabular}{rrrrrr}\toprule
$\mathcal{W}$ & W/T/L & Mean & 95\% interval & Median & Holm $p$\\\midrule
5 & 3/132/0 & 0.0 & [0.0, 0.1] & 0.0 & 0.125\\
6 & 14/121/0 & 0.3 & [0.1, 0.5] & 0.0 & $1.2\times10^{-4}$\\
8 & 27/107/1 & 0.6 & [0.3, 0.9] & 0.0 & $4.1\times10^{-6}$\\
10 & 44/91/0 & 1.7 & [1.2, 2.2] & 0.0 & $3.4\times10^{-13}$\\
12 & 63/69/3 & 2.5 & [1.9, 3.2] & 0.0 & $1.1\times10^{-10}$\\
15 & 75/56/4 & 3.1 & [2.5, 3.8] & 1.1 & $1.5\times10^{-12}$\\
18 & 79/52/4 & 2.9 & [2.2, 3.6] & 0.9 & $2.5\times10^{-13}$\\
\bottomrule\end{tabular}
\end{table}

The individual-task Spearman correlation between exact maximum rank and
relative gain is 0.435. It is positive within all nine $(m,n)$ groups,
with Holm-adjusted $p<0.006$ in each. A regression controlling for these
nine groups estimates 0.3 percentage points of additional gain per rank
doubling (95\% HC3 interval 0.3--0.4 percentage points, rounded).
The accompanying analysis notebook reports
all 63 cell-level A/B summaries, including losses and nonsignificant cells.

\paragraph{Separate strip family used for SA.}
The SA comparison uses a different construction: $m=3n$ cells on a
$3\times n$ rectangular strip, and $N=n^2$ binary tile-selection variables.
Eligible tiles are connected sets of one to four cells, with connectivity
through shared cell edges, spanning at most three strip columns.
All singletons are mandatory. We also include the dominoes corresponding
to every vertical edge and to the horizontal edges along the top row;
together these edges form a spanning tree of the strip.

Select the remaining tiles without replacement using size weights
$a(S)=\exp[-2(|S|-2)]$, where $|S|$ is the number of cells in a tile.
Among optional tiles, pairs thus have weight 1, triples $e^{-2}$ and
quadruples $e^{-4}$. The coefficient $-2$ is a fixed generator choice
favouring small tiles (the code's \texttt{light} setting), not a Boltzmann
temperature, a learned parameter, or a consequence of the objective.
The implementation draws independent $U_S\sim\operatorname{Uniform}(0,1)$,
forms exponential keys $-\log U_S/a(S)$,
and selects the required number of tiles with the smallest keys.

The benchmark inputs are generated and saved before optimization, then
reused unchanged for all methods. Spatial cell labels are restored and
tiles use the singleton-last order described above; the corresponding
variable permutation is also applied to the objective arrays. There is
no search over variable order. This family is a deliberately low-rank
contrast to row-major assignment, not a sample from the adjustable-locality
grid. Full-indicator states and ranks are obtained from feasible prefix
charges as above, and feasible covers are sampled using
\cref{eq:partition_counting}.

\FloatBarrier
\section{Simulated annealing comparison: full protocol and additional sizes}
\label{app:sa_protocol}
\paragraph{Instances, objectives and initialization.}
For each of assignment and strip Set Partitioning, we use 40 previously generated and saved instances
at each $n=10,12,14$, with $N=n^2$ binary variables. Assignment uses its
row-major permutation encoding ($2n$ equality rows); the strip construction
and its $3n$ equalities are specified in \cref{app:partition_protocol}.
Each instance is evaluated separately with $C(x)=x^\top Qx$ and the 32-scenario
entropic objective. The symmetric quadratic arrays are formed by symmetrizing
matrices with independent uniform entries in $[-7,7]$. Within each $(n,i)$,
both families use the same saved objective coefficient arrays and entropic
scale $s=\sqrt{3n}$, although their feasible spaces differ. Initial data are
400 independent uniform feasible draws, with replacement: random permutations
for assignment and exact counting for Set Partitioning. All four methods use
the same initial draw matrix within a problem, objective and budget.
The experiment comprises 5,760 runs: two families, two objectives, three sizes,
40 instances, three budgets and four methods.

\paragraph{Data-initialized EELS and full-indicator TT.}
Both start by deduplicating and evaluating the initial data. At every outer
iteration EELS rebuilds from the current training set and enhances the graph
with $\gamma=1$. The full TT constructs the exact complete indicator once
inside the timed run, then copies its uniformly normalized amplitudes at each
iteration. Both variants train for one sweep, using learning rate 0.05,
charge degeneracy one, and normalized Boltzmann weights with temperature
equal to the current cost standard deviation (uniform weights at zero spread).
After training they draw up to 10K samples, using the legacy backend in
batches of 100 for deadline checks. The next training set contains at most
400 distinct lowest-cost points from the previous training set and the new
batch; no external feasible samples are injected. The incumbent records all
evaluated feasible points, including a partially completed final batch.
The number of iterations is time-limited rather than fixed. 

\paragraph{Informed SA.}
Optim.jl's \texttt{SimulatedAnnealing} runs ten equally timed chains from the
ten best distinct initial points. For assignment a proposal swaps the assigned
columns of two distinct, uniformly chosen rows. For Set Partitioning a fixed
patch is chosen uniformly from the tile supports. If an occupied tile crosses
its boundary, the proposal is a self-loop; otherwise the current patch cover
is replaced by a uniformly chosen exact cover of that patch. Patch covers
are precomputed. These proposals are symmetric and preserve feasibility.
Each chain cools geometrically in wall time from the standard deviation of
initial costs to $10^{-3}$ times that value, with a machine-epsilon floor.
The moves and patch catalogue provide problem-specific information beyond
the explicit constraints and initial examples given to data-initialized EELS.

\paragraph{Penalty SA.}
The same Optim implementation and ten starting points minimize
$C(x)+\alpha\|Ax-b\|_2^2$ on the whole binary cube, proposing one uniformly
chosen bit flip. There is no repair or feasibility-preserving move catalogue.
With $R$ and $\sigma$ the range and sample SD of initial costs, set
\begin{align}
s_C&=\max\{R,\sigma,10^{-8}\max(1,\max_{x\in\mathcal T_0}|C(x)|)\},
&\alpha&=10s_C,\\
v&=\operatorname{median}_{j:\|A_j\|^2>0}\|A_j\|^2,
&T_0&=\max\{\sigma,\alpha v/\log2\},\\
&&T_{\mathrm{end}}&=\max\{10^{-3}\sigma,10^{-10}s_C\}.
\end{align}
Geometric cooling is restarted within each chain. Ignoring the objective
change, $T_0$ gives a typical constraint-violating one-bit proposal approximately
one-half acceptance probability, rather than freezing every move at a feasible
start. The multiplier is fixed across all instances and objectives. Only feasible raw costs
update the reported incumbent. Every objective evaluation on an infeasible
state nevertheless consumes a call and time. Residual computation, penalty
setup, feasible/infeasible counters and best penalized states are recorded.
No infeasible evaluated point had a penalized value below the initial feasible
best.

\paragraph{Timing, counters and reproducibility.}
Separate 10, 50 and 200-second runs reuse fixed per-instance seeds. The source
seed plus $10^6$ initializes SA's default RNG, plus $2\times10^6$ its proposal
RNG, and plus $3\times10^6$ the TT RNG. 
The cluster used Intel Xeon Gold 6226R CPUs at 2.90 GHz, Julia 1.12.5 and
Optim.jl 2.3.2. Each solver process has one Julia compute thread and one BLAS
thread; 32 independent processes share each exclusive allocation, with four
allocations at a time. Timed work includes initial objective
calls, TT construction, resets, training, sampling, penalty evaluation and
search. Loading, JIT warmup, pre-run validation and informed-SA patch-catalogue
construction are excluded; this last exclusion favors informed SA and limits
end-to-end runtime interpretation.

Termination is cooperative. Improvements after the requested deadline are
excluded, while actual elapsed times are retained. In particular, for quadratic
assignment at $n=14$ with a nominal 10-second budget, full TT averages 22.8
seconds (maximum 27.7) and completes no iteration on any task. At 200 seconds
all methods perform search. Exact objective counters include initial calls,
duplicates, repeated Optim evaluations, self-loops and penalty-SA infeasible
points. The same evaluator is used
by all methods: quadratic costs sum the active submatrix of $Q$, while the
entropic cost uses stable log-mean-exp.

\paragraph{Aggregation and inference.}
For each panel, average the per-task gain
$100(C_0-C_{\mathrm{best}})/|C_0|$ equally over its 40 instances. The curve
uses the exact union of improvement events with stepwise interpolation.
Pointwise asymmetric 95\% percentile-bootstrap intervals resample the 40
tasks 2,000 times (seed 20261014), using the same weights across time and
methods.
For the objective-call axis, each method stops at the minimum final count
over all 40 tasks; endpoints from slower tasks are not extrapolated. The
100,000-call diagnostic is a post hoc slice of the 200-second trajectories,
not a separately tuned fixed-query test. Full assignment at $n=14$ does not
reach that slice on all tasks.

Paired tests use one-sided Wilcoxon signed ranks of raw baseline-minus-TT
cost differences; differences below $10^{-9}+10^{-11}\max(|C_a|,|C_b|)$
are numerical ties. Tests versus penalty SA form a Holm family of 72 comparisons
(36 groups times two TT variants). The original informed-SA analysis, extended
to include penalty SA as a third candidate, uses 108 comparisons in a separate
Holm family. Both are exploratory analyses of this benchmark, with one seed
per task.

\begin{table}[!tbp]
\centering\small
\caption{200-second endpoints. Each cell is mean feasible-objective gain (\%) / mean objective calls (millions), over 40 tasks. Means include unsuccessful runs.}
\label{tab:sa_endpoints}
\setlength{\tabcolsep}{5pt}
\begin{tabular}{llrrrrr}\toprule
Family & Objective & $n$ & Informed SA & Penalty SA & EELS & Full TT\\\midrule
Assignment & Quadratic & 10 & 73.5 / 368.6 & 24.7 / 103.9 & 52.4 / 1.6 & 71.3 / 1.2\\
Assignment & Quadratic & 12 & 102.8 / 303.3 & 15.8 / 70.1 & 55.5 / 0.92 & 85.1 / 0.34\\
Assignment & Quadratic & 14 & 132.9 / 252.5 & 9.8 / 48.1 & 45.9 / 0.53 & 75.2 / 0.064\\
Assignment & Entropic & 10 & 93.3 / 125.2 & 27.7 / 62.7 & 76.4 / 1.7 & 92.7 / 1.2\\
Assignment & Entropic & 12 & 115.6 / 108.1 & 18.6 / 45.1 & 73.9 / 0.98 & 109.5 / 0.35\\
Assignment & Entropic & 14 & 141.1 / 95.6 & 11.7 / 32.4 & 54.6 / 0.57 & 106.2 / 0.067\\
Set Partitioning & Quadratic & 10 & 87.8 / 86.7 & 22.2 / 65.3 & 87.1 / 1.7 & 87.2 / 1.6\\
Set Partitioning & Quadratic & 12 & 117.2 / 79.7 & 17.3 / 59.8 & 115.9 / 1.0 & 116.7 / 0.88\\
Set Partitioning & Quadratic & 14 & 139.6 / 67.8 & 6.2 / 39.2 & 134.4 / 0.63 & 137.1 / 0.5\\
Set Partitioning & Entropic & 10 & 200.7 / 69.8 & 46.1 / 49.0 & 199.8 / 1.8 & 200.7 / 1.7\\
Set Partitioning & Entropic & 12 & 201.6 / 64.1 & 31.8 / 46.0 & 200.8 / 1.1 & 201.6 / 0.93\\
Set Partitioning & Entropic & 14 & 210.6 / 58.3 & 19.5 / 31.1 & 210.0 / 0.66 & 210.6 / 0.54\\
\bottomrule\end{tabular}
\end{table}

\begin{table}[!tbp]
\centering\footnotesize
\caption{Paired improvements over penalty SA at 200 seconds: mean $100(C_{\rm penalty}-C_{\rm TT})/|C_{\rm penalty}|$ (\%), wins/ties/losses, and one-sided Wilcoxon $p$ after Holm correction over all 72 TT--penalty comparisons (including 10 and 50 seconds).}
\label{tab:sa_penalty_tests}
\setlength{\tabcolsep}{4pt}
\begin{tabular}{llrrrrrrr}\toprule
Family & Objective & $n$ & EELS (\%) & W/T/L & Holm $p$ & Full TT (\%) & W/T/L & Holm $p$\\\midrule
Assignment & Quadratic & 10 & 23.3 & 40/0/0 & $6.5\times10^{-11}$ & 38.5 & 40/0/0 & $6.5\times10^{-11}$\\
Assignment & Quadratic & 12 & 35.0 & 39/0/1 & $6.5\times10^{-11}$ & 60.7 & 40/0/0 & $6.5\times10^{-11}$\\
Assignment & Quadratic & 14 & 33.5 & 40/0/0 & $6.5\times10^{-11}$ & 60.0 & 40/0/0 & $6.5\times10^{-11}$\\
Assignment & Entropic & 10 & 38.3 & 40/0/0 & $6.5\times10^{-11}$ & 51.3 & 40/0/0 & $6.5\times10^{-11}$\\
Assignment & Entropic & 12 & 48.4 & 38/1/1 & $6.5\times10^{-11}$ & 78.8 & 40/0/0 & $6.5\times10^{-11}$\\
Assignment & Entropic & 14 & 41.2 & 39/0/1 & $6.5\times10^{-11}$ & 87.9 & 40/0/0 & $6.5\times10^{-11}$\\
Set Partitioning & Quadratic & 10 & 53.7 & 40/0/0 & $6.5\times10^{-11}$ & 53.8 & 40/0/0 & $6.5\times10^{-11}$\\
Set Partitioning & Quadratic & 12 & 85.4 & 40/0/0 & $6.5\times10^{-11}$ & 86.0 & 40/0/0 & $6.5\times10^{-11}$\\
Set Partitioning & Quadratic & 14 & 121.8 & 40/0/0 & $6.5\times10^{-11}$ & 124.3 & 40/0/0 & $6.5\times10^{-11}$\\
Set Partitioning & Entropic & 10 & 109.3 & 40/0/0 & $6.5\times10^{-11}$ & 109.7 & 40/0/0 & $6.5\times10^{-11}$\\
Set Partitioning & Entropic & 12 & 125.5 & 40/0/0 & $6.5\times10^{-11}$ & 126.1 & 40/0/0 & $6.5\times10^{-11}$\\
Set Partitioning & Entropic & 14 & 159.1 & 40/0/0 & $6.5\times10^{-11}$ & 159.6 & 40/0/0 & $6.5\times10^{-11}$\\
\bottomrule\end{tabular}
\end{table}

\begin{figure}[!tbp]
\centering
\includegraphics[width=0.76\textwidth]{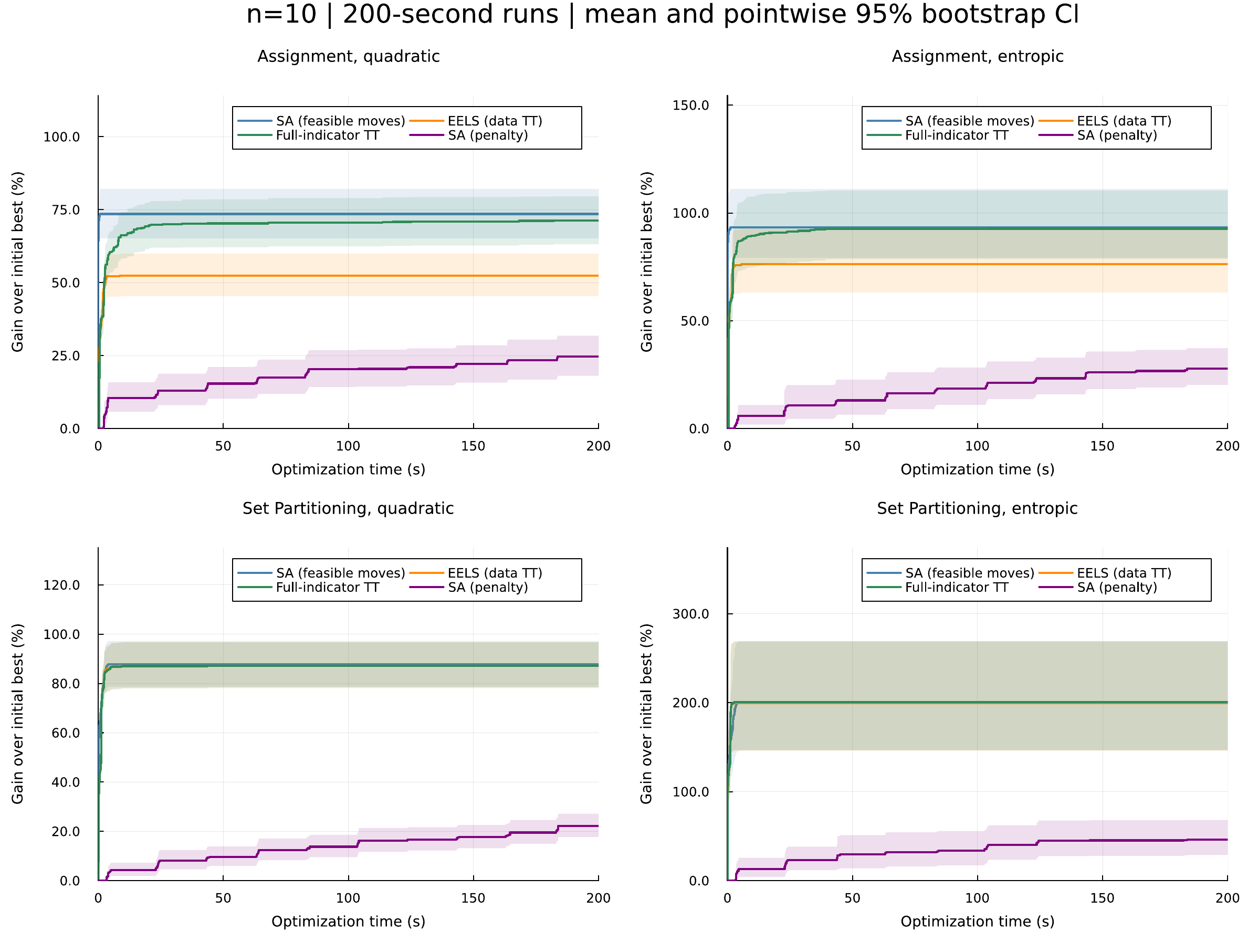}
\caption{Additional size $n=10$ ($N=100$): 200-second runs, 40 tasks per panel. Gain versus optimization time; conventions match Figure~14 of the main article. Top: assignment; bottom: strip Set Partitioning. Left: quadratic; right: entropic. All four methods use matched initial data; TT amplitudes are reset/rebuilt each iteration. Shading denotes pointwise 95\% bootstrap intervals for mean per-task gain.}
\label{fig:sa_n10_time}
\end{figure}

\begin{figure}[!tbp]
\centering
\includegraphics[width=0.76\textwidth]{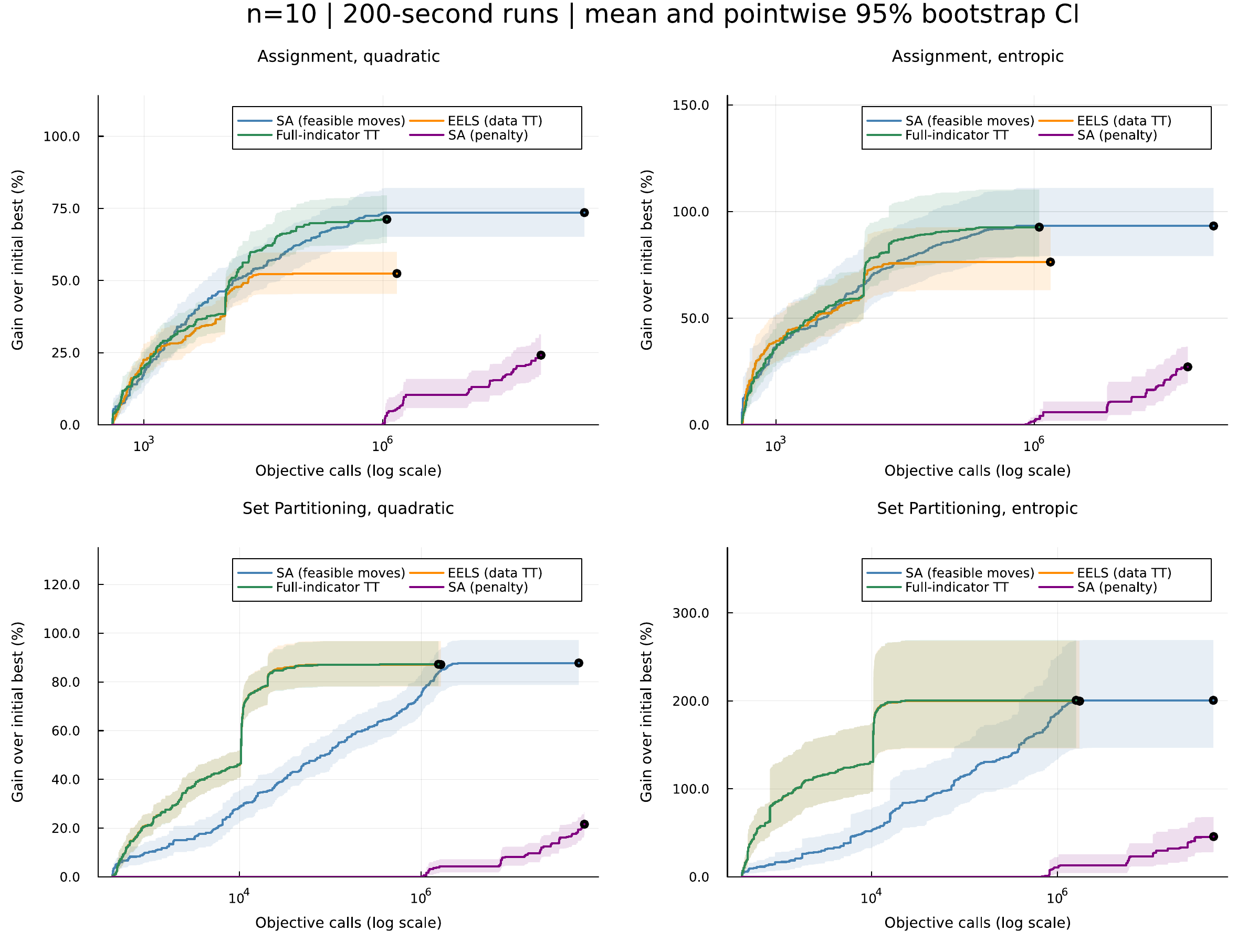}
\caption{Additional size $n=10$ ($N=100$): 200-second runs, 40 tasks per panel. Gain versus objective calls; conventions and common-support truncation match Figure~15 of the main article. Top: assignment; bottom: strip Set Partitioning. Left: quadratic; right: entropic. All four methods use matched initial data; TT amplitudes are reset/rebuilt each iteration. Shading denotes pointwise 95\% bootstrap intervals for mean per-task gain.}
\label{fig:sa_n10_calls}
\end{figure}

\begin{figure}[!tbp]
\centering
\includegraphics[width=0.76\textwidth]{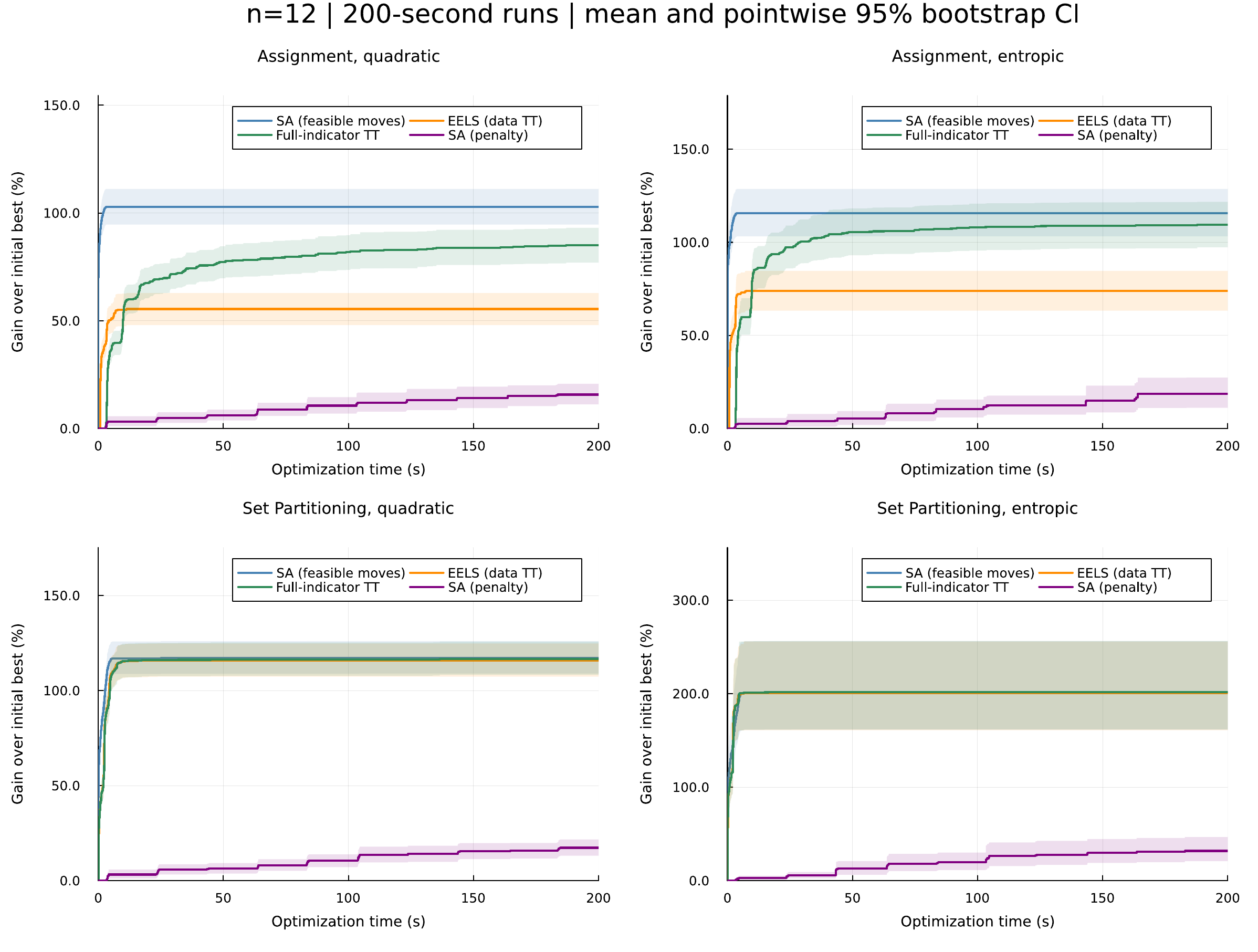}
\caption{Additional size $n=12$ ($N=144$): 200-second runs, 40 tasks per panel. Gain versus optimization time; conventions match Figure~14 of the main article. Top: assignment; bottom: strip Set Partitioning. Left: quadratic; right: entropic. All four methods use matched initial data; TT amplitudes are reset/rebuilt each iteration. Shading denotes pointwise 95\% bootstrap intervals for mean per-task gain.}
\label{fig:sa_n12_time}
\end{figure}

\begin{figure}[!tbp]
\centering
\includegraphics[width=0.76\textwidth]{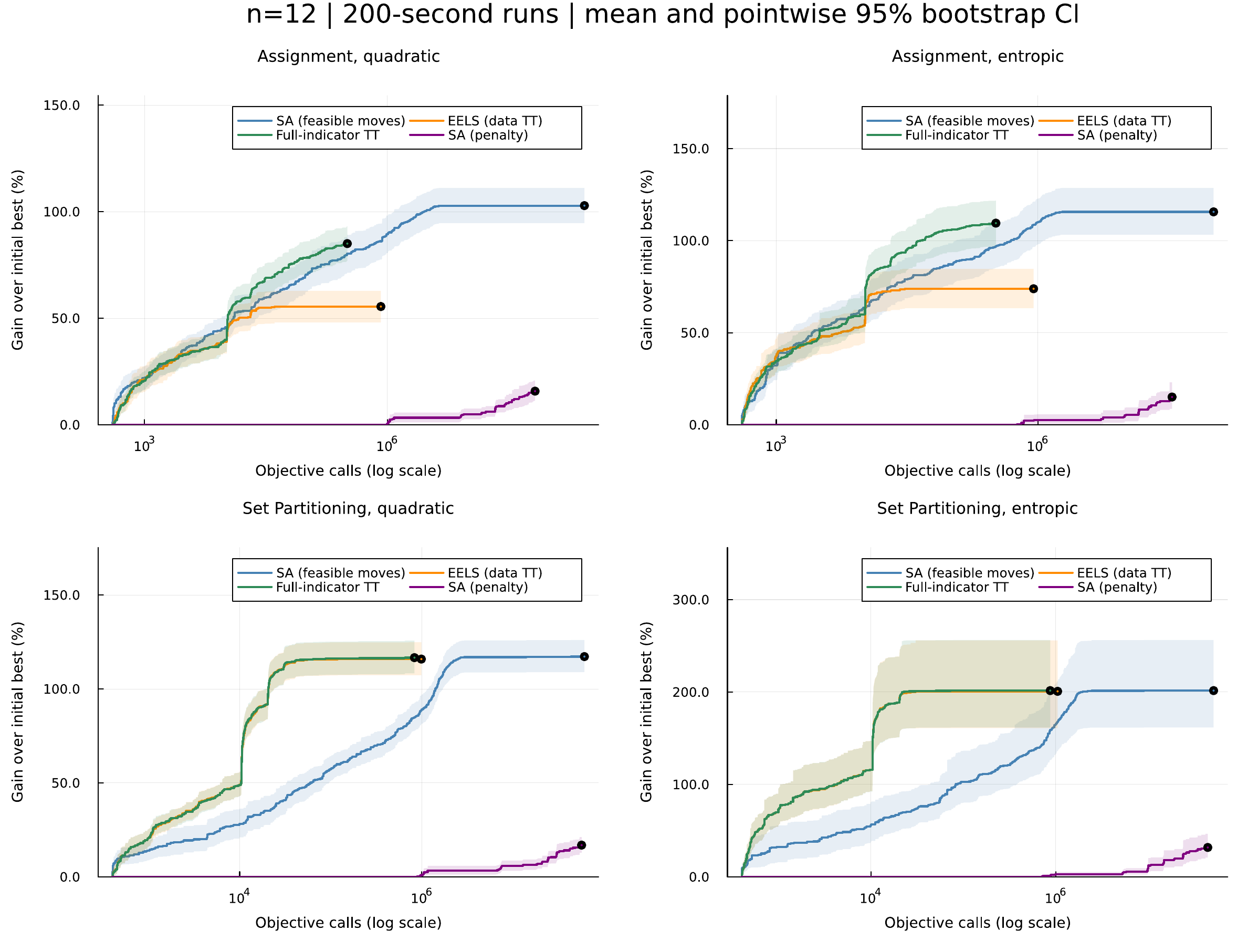}
\caption{Additional size $n=12$ ($N=144$): 200-second runs, 40 tasks per panel. Gain versus objective calls; conventions and common-support truncation match Figure~15 of the main article. Top: assignment; bottom: strip Set Partitioning. Left: quadratic; right: entropic. All four methods use matched initial data; TT amplitudes are reset/rebuilt each iteration. Shading denotes pointwise 95\% bootstrap intervals for mean per-task gain.}
\label{fig:sa_n12_calls}
\end{figure}


\FloatBarrier
\section{Detailed library performance benchmarking}
\label{app:libs_comparison}

\subsection{Matched inputs and benchmark protocol}
We compared the training and sampling kernels of the two libraries on the
binary assignment problem. For an $n\times n$ assignment instance there are
$n^2$ binary variables and $2n$ linear equality constraints. We used
$4 \leq n \leq 10$, one complete right--left training sweep, one batch of
$|T|=400$ samples after training, bond degeneracy one, and learning rate
$0.05$. A deterministic dense quadratic cost matrix with entries in
$\{-7,\dots,7\}$ was generated for every $n$.
Four hundred random feasible assignments were generated and resampled once
according to their Boltzmann weights. The resulting empirical training matrix,
including repeated rows, was passed unchanged to both implementations. Thus,
both libraries were trained on exactly the same data.
OhMyU1 is measured with two sampling backends. The \texttt{legacy}
backend follows sampled charge paths using the block dictionaries;
\texttt{table} first constructs transition tables for all charge states on
each sampling call, then amortizes that preparation over the batch.
Both implement the same conditional probabilities. Training is unchanged;
the recorded training-time differences between runs reflect timing variation,
not a different training algorithm.

The initial tensors represented the same mathematical state: the uniform
superposition of all feasible assignments. ConstrainTNet constructed its state
exactly from the constraints. For OhMyU1, the complete set of reachable charge
sectors at every cut was enumerated exactly by a dynamic construction specific
to the assignment variable order; no finite random initialization sample was
used. The maximum bond
supports were identical in both representations and ranged from 9 at $n=4$ to
461 at $n=10$.

Tensor construction, initial-state canonicalization and normalization, Julia compilation, and
copying of the prepared state were excluded from the reported time. Each
kernel was executed once as a warm-up, followed by three measured repetitions;
the median wall time is reported. The primary metric is the sum of the median
time for one training sweep and the median time for exactly 400 sequential
samples from the trained state. All runs used Julia 1.12.5, one default Julia thread, one BLAS thread,
and the same Intel Core i5-10210U CPU.

\subsection{Implementation details and validation}

The refreshed benchmark uses OhMyU1 revision \texttt{e23712d} and ConstrainTNet
revision \texttt{77d8c26} (latest upstream checked on 12 October 2026).
ConstrainTNet's algorithm sources are unmodified; vendored
ITensors~0.3.41/NDTensors~0.1.50 are retained.
The sole source adjustment increases ITensors' fixed charge capacity from
20 to 21: at 20 equalities, comparison with an empty charge attempts to insert
an inactive component and otherwise fails. The spare slot avoids this boundary
error.

The random seed for each problem is $20260923+10000n$. The Boltzmann
temperature is the sample standard deviation of the 400 candidate costs;
400 draws with replacement form the common empirical training multiset.

One default Julia thread, no interactive thread, and one BLAS thread are used
on Windows with an Intel Core i5-10210U processor. The two backends run in
separate sequential processes. Each timing starts from a fresh copy, with
garbage collection outside the timed region. Compilation, construction,
initial-state canonicalization/normalization, and copying are excluded.
Each kernel is warmed up and timed three times; the table adds the two kernel
medians rather than timing a high-level optimizer wrapper. This avoids the
extra pre-training sampling batch in \texttt{ConstrainTNet.optimizer}.

OhMyU1 performs projected two-site updates, whereas ConstrainTNet performs
one-site updates. Internal canonicalization, environment construction, and
loss evaluation performed by the public ConstrainTNet training routine remain
inside its timing. Thus, the measurements compare the cost of the available
training and sampling routines for the same inputs and sweep count.

\begin{table}[!tbp]\centering\small
\caption{Median times (seconds) for one training sweep and 400 serial samples.
L and T denote legacy and table. Total is the sum of the two kernel medians.
ConstrainTNet measurements are retained; both OhMyU1 variants use the same inputs.
Construction, warm-up, and copying are excluded.}
\label{tab:lib_comparison_fair}
\begin{tabular}{rrrrrrrr}\toprule
$n$ & Train L & Sample L & Total L & Train T & Sample T & Total T & ConstrainTNet\\\midrule
4 & 0.0025 & 0.0033 & 0.0058 & 0.0024 & 0.0004 & 0.0028 & 1.2\\
5 & 0.0061 & 0.0059 & 0.012 & 0.0058 & 0.0008 & 0.0066 & 2.2\\
6 & 0.013 & 0.0096 & 0.023 & 0.013 & 0.0018 & 0.015 & 3.4\\
7 & 0.028 & 0.013 & 0.041 & 0.028 & 0.0042 & 0.032 & 5.2\\
8 & 0.059 & 0.019 & 0.078 & 0.055 & 0.0091 & 0.064 & 8.5\\
9 & 0.15 & 0.029 & 0.18 & 0.14 & 0.021 & 0.16 & 14\\
10 & 0.3 & 0.044 & 0.34 & 0.31 & 0.045 & 0.36 & 47\\
\bottomrule\end{tabular}\end{table}


\FloatBarrier
\section{PROTES comparison: exact indicators and data-based initialization}
\label{app:protes_protocol}

OhMyU1 measurements were repeated on 12 October 2026 using revision
\texttt{e23712d}, with unchanged inputs, seeds, CPU affinity,
batch sizes, and time budgets. The PROTES measurements were retained.
The source hashes and separate memory-audit and timing runs are archived.
OhMyU1 uses the serial \texttt{sample\_nondeg!} routine, with both
\texttt{backend=:legacy} and \texttt{backend=:table} reported separately.
Full-indicator table timings are the archived repeat after the memory audit;
data-initialized table runs were measured separately with identical inputs.

\subsection{Exact assignment indicator, ranks, and storage}
The full structural tensor is $I(x)=\mathbf1\{Ax=b\}$. For assignment, reshape the $N=n^2$ binary variables as
$x_{ij}$ with exactly one selected entry in every row and column. Process
these bits row by row. At a cut after $k=rn+s$ bits, $r$ rows are complete
and $s$ bits of the next row have been processed, with $0\le s<n$.

A dynamic-programming state records the set $S$ of columns already used
by selected entries and the flag $h\in\{0,1\}$ indicating whether the
current, partially processed row already has its selected entry. The
completed rows contribute exactly $r$ selected entries, so $h=|S|-r$.
For example, for $n=3$ at $k=4$ we have $r=1,s=1$. State
$(S=\{2\},h=0)$ means that the first row selected column 2 and the
second row skipped column 1. State $(S=\{1,2\},h=1)$ means that the
second row has already selected column 1. These states have different
sets of feasible suffixes.

To process the next bit in column $s+1$, choosing zero leaves $S$ unchanged;
choosing one is allowed only if $h=0$ and $s+1\notin S$, and adds that
column to $S$. At the end of each row require exactly one selected entry
and reset the flag for the next row. Merge identical states. Keep only
states reachable from the empty prefix and extendable to a complete
assignment. For $h=0$, an unoccupied column must remain among
$s+1,\ldots,n$. For $h=1$, at least one occupied column must be among
$1,\ldots,s$, so the current row's selection can lie in the processed prefix.
At the endpoints, $S=\varnothing$ and $S=\{1,\ldots,n\}$.

Each retained state labels a virtual index. Set a core entry to one for
each permitted state transition and to zero otherwise. Every permutation
has exactly one binary path, hence unit amplitude; all other binary
vectors have zero amplitude. This constructs the full indicator without
enumerating $n!$ permutations. The symmetric and dense representations
encode this same graph; the dense one additionally stores forbidden
transitions as zero coefficients. Let $R_k$ denote the total virtual
dimension at cut $k$, equivalently the number of retained charge states,
as in the Set Partitioning construction.

With out-of-range binomial coefficients interpreted as zero, the number of
states at an internal cut is
\begin{equation}
 R_k=\binom nr-\binom{s}{r-n+s}
       +\binom n{r+1}-\binom{n-s}{r+1},\qquad R_0=R_N=1.
\end{equation}
Here $r$ is the number of completed rows and $s$ the number of processed bits in the current row. The first pair of terms counts $h=0$ states and the second counts $h=1$ states. After grouping prefixes and suffixes by residual constraint charge,
the indicator unfolding consists of disjoint nonempty all-ones blocks,
each of rank one. Consequently $R_k$ is the exact minimal ordinary TT rank
for this ordering. The maximum is $\Theta(2^n/\sqrt n)$ (e.g., 125 at $n=8$,
461 at $n=10$, and 24309 at $n=16$). The $U(1)$ representation has the same
total link dimension, but every sector has degeneracy one. This exact-rank
statement concerns the full indicator; an optimizer's distribution concentrated
on a smaller feasible subset need not have these ranks.

The two dense float64 storage counts in Figure~17 of the main article are
\begin{equation}
 B_{\rm uniform}=8[2(N-2)R_{\max}^2+4R_{\max}],\qquad
 B_{\rm variable}=16\sum_{k=1}^{N}R_{k-1}R_k.
\end{equation}
They describe the fast PROTES core layout and the variable-rank
\texttt{protes\_general} layout, respectively. For $U(1)$, the number of
allowed scalar transitions at position $(r,s)$ is
\begin{equation}
 E_{r,s}=\binom{n-1}{r}+\binom nr-\binom{s+1}{r-n+s+1}
                 +\binom n{r+1}-\binom{n-s}{r+1}.
\end{equation}
The first term counts transitions assigning one and the remaining terms
count transitions assigning zero. Thus numerical block storage is
$8\sum_{r=0}^{n-1}\sum_{s=0}^{n-1}E_{r,s}$ bytes. These counts were checked
against the constructed blocks and dense arrays. MiB and GiB denote $2^{20}$
and $2^{30}$ bytes. Full $U(1)$ objects were constructed and measured through
$n=16$ using \texttt{Base.summarysize}, including dictionaries, charge vectors,
the constraint matrix, and other reachable fields. The measurement excludes
process/runtime memory and temporary buffers; it depends on object sharing
in this constructor. Dense curves count coefficient arrays only. Adam states,
gradients, and training workspaces are excluded from all storage claims.

\subsection{Untrained full-indicator sampling protocol}
Both models sample uniformly over the same full feasible set. PROTES uses
the above zero/one transition cores, padded to a common internal rank for
its fast implementation and left at their exact ranks for
\texttt{protes\_general}. For OhMyU1, backward dynamic programming computes
the completion count $C_k(q)$ at every link: $k$ is the number of
processed binary variables and $q=b-\sum_{j=1}^{k}A_jx_j$ is the residual
constraint vector. Thus $C_k(q)$ counts suffixes $(x_{k+1},\ldots,x_N)$
whose contribution equals $q$. Starting with $C_N(0)=1$, compute backwards
\[
C_k(q)=\sum_{v:\,q\xrightarrow{v}q'}C_{k+1}(q'),
\qquad q'=q-vA_{k+1},
\]
where the sum includes only retained feasible transitions. For example, at
a row boundary of assignment, $C_k(q)=(n-r)!$ after $r$ completed rows:
the remaining rows can be matched bijectively to the unoccupied columns.
The scalar on an allowed edge $q\to q'$ is replaced by
$\sqrt{C_{k+1}(q')/C_k(q)}$. Products telescope to $1/\sqrt{n!}$, giving an
exact right-canonical normalized amplitude. Its squared amplitudes and the
PROTES probability indicator therefore induce the same uniform law.

We call the unchanged OhMyU1 serial \texttt{sample\_nondeg!} and the upstream
PROTES sampling kernels with their usual JAX vectorization. Batches contain 100 samples. Each of three seeds runs
until cumulative timed generator calls reach at least one second; a final
batch can exceed that duration. Throughput divides the actual generated
count by actual timed duration, and counts repetitions. Construction,
canonicalization, JIT warm-up, feasibility checks, and unique-point bookkeeping
are excluded. PROTES timing includes materializing the output and recomputing
right environments per batch, as in its normal sampling loop.
The feasible fraction is 100.0\% in all recorded batches. Raw per-trial
batch counts, times, warm-up/construction times, and distinct-point counts
are retained in the accompanying results snapshot. Full dense models are
allocated for $n=4,\ldots,10$; beyond this grid their storage curves are
analytical, not timing measurements. OhMyU1 is measured for $n=4,\ldots,16$.

The large $U(1)$ object traversal initially caused memory pressure at
$n=15,16$ on the local machine. Final sampling measurements were repeated
without repeating that deep traversal and with garbage collection before
warm-up; the initial memory-audit runs are retained separately. The shaded
ranges report run-to-run variation, not confidence intervals or an
asymptotic complexity result.

\subsection{Paired, data-initialized optimization protocol}

\begin{figure}[!tbp]
\centering\includegraphics[width=0.76\textwidth]{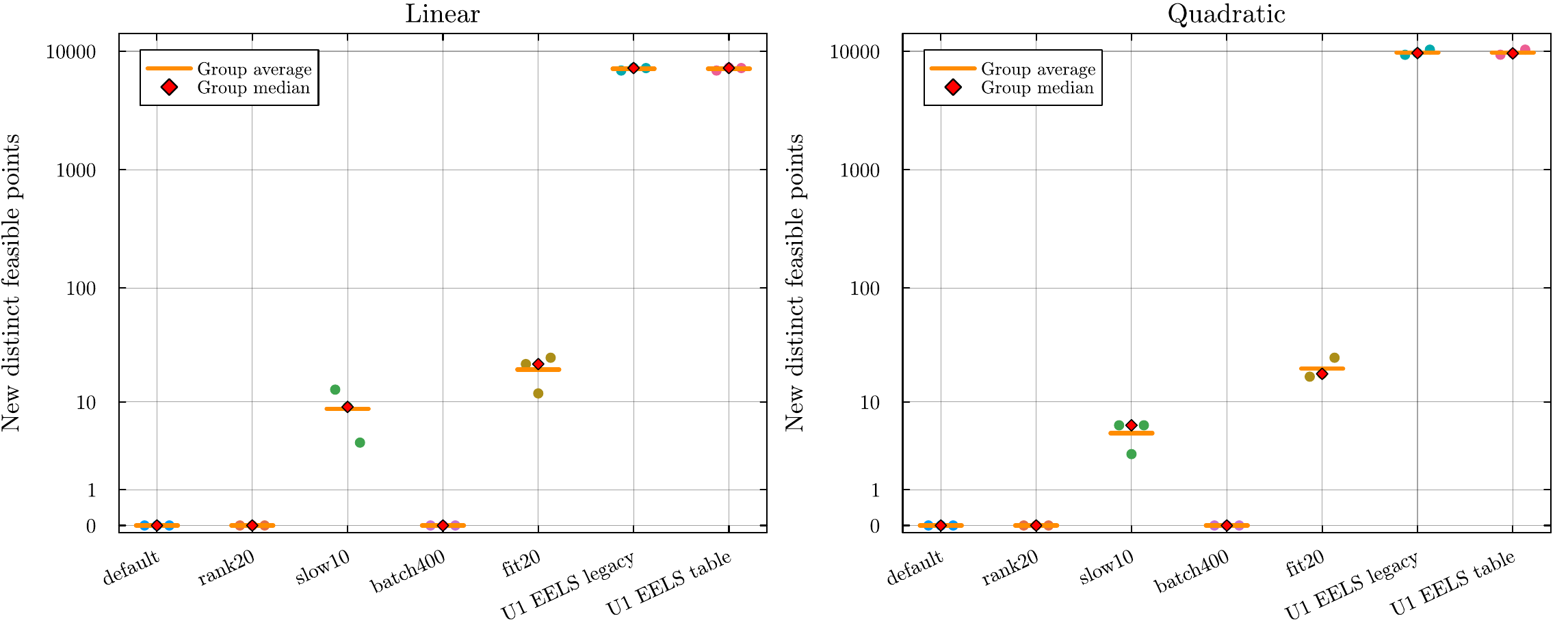}
\caption{The same data-initialized comparison at $n=8$ ($N=64$), with a
50-second budget and three matched instance/seed pairs for each objective.
Dots are individual counts of new distinct feasible points; orange bars
and red diamonds show the arithmetic mean and median of the original counts.
Both panels use the same $\log_{10}(1+z)$ ordinate and limits. All five PROTES
configurations and both EELS sampling backends are reported, not only the best observed configuration.
These are trained, time-limited runs, unlike the untrained full-indicator
sampling experiment in Figure~17 of the main article.}
\label{fig:protes_data_long}
\end{figure}

Longer initial fitting improves PROTES: at $n=8$, \texttt{fit20} finds
12--25 new points for linear objectives and 17--25 for quadratic objectives
in 50 seconds, while EELS with the legacy sampler finds 6.9K--7.3K and 9.3K--10K, respectively.
Thus, inability to generate any feasible point is not a robust conclusion;
the observed limitation is rare feasible generation under this data-based
initialization and update rule. New points need not improve the objective.
With the table sampler the corresponding EELS ranges are 6.9K--7.3K and 9.4K--10K.

These diagnostic instances are separate from the article's objective
files. For each $(n,j,o)$, with repetition index $j$ and objective family $o$,
NumPy's generator uses seed $20261008+1000n+10j+\mathbf{1}_{o=\mathrm{quadratic}}$.
It first draws 400 independent random permutations and converts them to
row-major binary assignments. It then draws either independent
$c_i\sim U[-50,50]$, with cost $c^Tx$, or a matrix with independent
$B_{ij}\sim U[-7,7]$, sets $Q=(B+B^T)/2$, and uses cost $x^TQx$.
Exactly the same arrays are exported to Julia.

Define novelty as the number of distinct generated feasible points outside
the set of initial observations. Initial duplicates are not new, repeated
generated points count once, and the initial best cost is retained in every
method. The 5-second screen uses $n=4,\ldots,10$ and $j=0$; the 50-second
comparison uses $n=8$ and $j=0,1,2$, for both objective families. All five
PROTES configurations in \cref{tab:protes_settings} are included at both
budgets. At $n=4$ the particular initial datasets cover all 24 solutions.

\begin{table}[!tbp]
\centering\small
\caption{PROTES settings. $k$ is the input/generated batch size, $k_{\rm top}$
the elite count, $g$ the number of Adam steps per accepted batch, and $e$
the number of passes through the initial 400 observations. The last column
is the number of initial Adam updates, $(400/k)eg$.}
\label{tab:protes_settings}
\begin{tabular}{lrrrrrrr}\toprule
Configuration & Rank & $k$ & $k_{\rm top}$ & $g$ & Learning rate & $e$ & Initial updates\\\midrule
default &5&100&10&1&0.05&1&4\\
rank20 &20&100&10&1&0.05&1&4\\
slow10 &5&100&10&10&0.01&1&40\\
batch400 &10&400&40&10&0.01&1&10\\
fit20 &5&100&10&1&0.05&20&80\\\bottomrule
\end{tabular}
\end{table}

PROTES starts from its ordinary random dense cores, receives the initial
observations through its external-sample interface in deterministically
shuffled passes, and applies its normal elite update. Thus only the best
10.0\% of each input batch contributes to the likelihood update, not all
400 observations simultaneously. After initialization
no more external observations are injected. The external evaluator assigns
valid points score $(2/\pi)\arctan C(x)$ and invalid points score 2, preserving
the ordering of feasible costs. An entirely infeasible batch uses the
upstream skip-update path. If fewer than $k_{\rm top}$ points are feasible,
the unchanged elite rule also includes penalized infeasible points.
No mask, repair, indicator, $A,b$, or permutation encoding is given to PROTES.

OhMyU1 receives $A,b$, deduplicates and cost-sorts the same initial data,
and rebuilds its TT from the retained data at each outer iteration, applying
EELS with $\gamma=1$ and sector degeneracy one before canonicalization.
Training uses one right--left sweep, learning rate 0.05, and normalized
weights proportional to $\exp[-(C(x)-C_{\min})/T]$, with $T$ the sample
standard deviation of the current costs; weights are uniform if $T=0$.
The normalized trained TT supplies up to 10K samples, processed in batches
of 100 for deadline checks. The next dataset contains up to 400 lowest-cost
distinct points from the union of retained and newly generated points.
No external feasible generator is called after initialization. The number
of iterations is set by the wall-time limit and reported below.

Initial model construction/fitting, subsequent training, sampling, novelty
bookkeeping, and objective evaluation all count toward the 5/50-second
budgets. Explicit warm-up and file loading do not: PROTES warms its kernels
on a synthetic binary objective at the same dimension, while OhMyU1 warms
one full outer iteration and discards its model. Residual compilation counts
toward the measured budget. Deadlines are cooperative: a current sample or
training operation can overrun, and a batch arriving after the deadline is
not accepted. All initial fits completed. There is no fixed number of
PROTES optimization iterations; generated counts and batch sizes determine
the sampling workload, whereas only batches containing feasible points
trigger further Adam steps. Total post-initialization Adam-update counts
were not separately logged. The result tables give generated counts rather than implying equal iteration counts; actual elapsed times remain available in the archived results.

All experiments ran locally on an Intel Core i5-10210U (four physical cores,
eight logical processors, approximately 8 GB RAM), Windows 11, with each
process pinned to logical CPU 0 and BLAS restricted to one thread. OhMyU1
uses Julia 1.12.5 with one thread; PROTES uses Python 3.12.14, NumPy 1.26.4,
JAX 0.4.38 (CPU/x64), and Optax 0.2.4. The upstream revision is
\texttt{06d447845dd26ffccbc8e95c7545fa1d2617940e}.
The scripts, compact results, hashes, and Julia plotting code accompany
the notebook \texttt{examples/PROTES\_assignment\_scaling.ipynb}.

\begin{table}[!tbp]
\centering\small
\caption{Recorded 50-second PROTES runs at $n=8$. Triples are new distinct feasible points in seed order $j=0,1,2$; ranges cover the three runs. Generated counts include infeasible draws.}
\label{tab:protes_long_counts}
\begin{tabular}{llrr}\toprule
Configuration & Objective & New points & Generated draws (range)\\\midrule
default & linear & 0, 0, 0 & 3.3M--3.4M\\
default & quadratic & 0, 0, 0 & 3.1M--3.4M\\
rank20 & linear & 0, 0, 0 & 1.2M--1.5M\\
rank20 & quadratic & 0, 0, 0 & 1.4M--1.6M\\
slow10 & linear & 13, 9, 4 & 3M--3.1M\\
slow10 & quadratic & 6, 3, 6 & 2.7M--2.9M\\
batch400 & linear & 0, 0, 0 & 3.4M--3.6M\\
batch400 & quadratic & 0, 0, 0 & 3.5M--3.5M\\
fit20 & linear & 22, 12, 25 & 2.9M--3.3M\\
fit20 & quadratic & 17, 18, 25 & 3.3M--3.4M\\
\bottomrule
\end{tabular}\end{table}
\begin{table}[!tbp]
\centering\small
\caption{Recorded EELS workloads with the legacy sampler. Iterations are completed/started outer iterations; the last started iteration may be interrupted. All generated draws counted here are feasible. Novelty excludes the initial observations.}
\label{tab:u1_timed_counts}
\begin{tabular}{rlrrrrr}\toprule
$n$ & Objective & $j$ & Budget (s) & Iterations & Draws & New points\\\midrule
4 & linear & 0 & 5 & 38/39 & 380K & 0\\
4 & quadratic & 0 & 5 & 39/40 & 390K & 0\\
5 & linear & 0 & 5 & 25/26 & 260K & 6\\
5 & quadratic & 0 & 5 & 26/27 & 270K & 3\\
6 & linear & 0 & 5 & 16/17 & 170K & 418\\
6 & quadratic & 0 & 5 & 15/16 & 160K & 406\\
7 & linear & 0 & 5 & 10/11 & 110K & 3K\\
7 & quadratic & 0 & 5 & 10/11 & 100K & 3.5K\\
8 & linear & 0 & 5 & 7/8 & 77K & 6.6K\\
8 & quadratic & 0 & 5 & 6/7 & 60K & 7.9K\\
9 & linear & 0 & 5 & 4/5 & 41K & 8K\\
9 & quadratic & 0 & 5 & 3/4 & 34K & 11K\\
10 & linear & 0 & 5 & 3/4 & 34K & 8.4K\\
10 & quadratic & 0 & 5 & 3/4 & 32K & 11K\\
8 & linear & 0 & 50 & 72/73 & 720K & 6.9K\\
8 & linear & 1 & 50 & 59/60 & 590K & 7.3K\\
8 & linear & 2 & 50 & 60/61 & 610K & 7.2K\\
8 & quadratic & 0 & 50 & 55/56 & 560K & 9.3K\\
8 & quadratic & 1 & 50 & 56/57 & 560K & 9.7K\\
8 & quadratic & 2 & 50 & 56/57 & 570K & 10K\\
\bottomrule
\end{tabular}\end{table}

\begin{table}[!tbp]\centering\small
\caption{EELS workloads with the table sampler under the same inputs and budgets as the legacy runs. Iterations are completed/started; novelty excludes initial observations. All counted draws are feasible.}
\label{tab:u1_timed_counts_table}
\begin{tabular}{rlrrrrr}\toprule
$n$ & Objective & $j$ & Budget (s) & Iterations & Draws & New points\\\midrule
4 & linear & 0 & 5 & 90/91 & 900K & 0\\
4 & quadratic & 0 & 5 & 74/74 & 740K & 0\\
5 & linear & 0 & 5 & 42/43 & 430K & 6\\
5 & quadratic & 0 & 5 & 43/44 & 440K & 3\\
6 & linear & 0 & 5 & 18/19 & 180K & 418\\
6 & quadratic & 0 & 5 & 20/21 & 210K & 406\\
7 & linear & 0 & 5 & 10/11 & 110K & 3K\\
7 & quadratic & 0 & 5 & 9/10 & 91K & 3.4K\\
8 & linear & 0 & 5 & 6/7 & 65K & 6.5K\\
8 & quadratic & 0 & 5 & 5/6 & 50K & 7.7K\\
9 & linear & 0 & 5 & 3/4 & 30K & 7.9K\\
9 & quadratic & 0 & 5 & 2/3 & 29K & 11K\\
10 & linear & 0 & 5 & 1/2 & 15K & 7.7K\\
10 & quadratic & 0 & 5 & 1/2 & 13K & 8K\\
8 & linear & 0 & 50 & 73/74 & 730K & 6.9K\\
8 & linear & 1 & 50 & 66/67 & 660K & 7.3K\\
8 & linear & 2 & 50 & 61/62 & 620K & 7.2K\\
8 & quadratic & 0 & 50 & 58/59 & 580K & 9.4K\\
8 & quadratic & 1 & 50 & 51/52 & 510K & 9.6K\\
8 & quadratic & 2 & 50 & 56/57 & 570K & 10K\\
\bottomrule
\end{tabular}\end{table}
\begin{table}[!tbp]\centering\small
\caption{PROTES 5-second screen: workload ranges over the 14 size/objective pairs per configuration. Each pair uses one seed. Exact per-instance counts and times are stored in the results snapshot.}
\label{tab:protes_short_workload}
\begin{tabular}{lr}\toprule Configuration & Generated draws (range)\\\midrule
default & 34K--390K\\
rank20 & 54K--140K\\
slow10 & 28K--260K\\
batch400 & 78K--470K\\
fit20 & 120K--340K\\
\bottomrule
\end{tabular}\end{table}

\end{document}